\documentclass[aop,seceqn,MSNbibl,citesort,dvips]{arximspdf}
\usepackage{graphicx}

\doi{10.1214/10-AOP546}
\volume{38}
\issue{6}
\pubyear{2010}
\firstpage{2322}
\lastpage{2378}

\makeatletter

\newproclaim{rem}{Remark}[section]
\newtheorem{theorem}[rem]{Theorem}
\newtheorem{proposition}[rem]{Proposition}
\newtheorem{lemma}[rem]{Lemma}
\newproclaim{definition}[rem]{Definition}

\newcommand{\sign}{\operatorname{sign}}

\newcommand{\ind}{\mathbf{1}}

\newcommand{\R}{\mathbb{R}}
\newcommand{\Z}{\mathbb{Z}}
\newcommand{\N}{\mathbb{N}}

\newcommand{\cF}{{\mathcal F}}
\newcommand{\cE}{{\mathcal E}}
\newcommand{\cH}{{\mathcal H}}
\newcommand{\cN}{{\mathcal N}}
\newcommand{\cL}{{\mathcal L}}
\newcommand{\cD}{{\mathcal D}}
\newcommand{\cS}{{\mathcal S}}
\newcommand{\cM}{{\mathcal M}}

\newcommand{\bP}{{\mathbf P}}
\newcommand{\bE}{{\mathbf E}}

\newcommand{\dd}{{d}}
\newcommand{\union}{\bigcup}

\newcommand{\bbE}{{\mathbb E}}
\newcommand{\bbP}{{\mathbb P}}

\newcommand{\ga}{\alpha}
\newcommand{\gb}{\beta}
\newcommand{\gga}{\gamma}
\newcommand{\gep}{\varepsilon}
\newcommand{\gD}{\Delta}
\newcommand{\go}{\omega}
\newcommand{\gl}{\lambda}
\newcommand{\gs}{\sigma}

\newcommand{\tf}{\mbox{\textsc{f}}}
\newcommand{\M}{\mbox{\textsc{M}}}
\newcommand{\slope}{\widetilde{\mathrm{m}}}

\newcommand{\myZ}{\mathcal{Z}}
\newcommand{\tg}{\mbox{\textsc{g}}}

\makeatother

\begin{document}
\begin{frontmatter}

\title{The weak coupling limit of disordered copolymer~models\thanksref{T1}}
\runtitle{Weak coupling limit for copolymer models}

\thankstext{T1}{Supported in part by ANR Grant POLINTBIO and by the University
of Padova Grant CPDA082105/08.}

\begin{aug}
\author[A]{\fnms{Francesco} \snm{Caravenna}\corref{}\ead[label=e1]{francesco.caravenna@math.unipd.it}} and
\author[B]{\fnms{Giambattista} \snm{Giacomin}\ead[label=e2]{giacomin@math.jussieu.fr}}
\runauthor{F. Caravenna and G. Giacomin}
\affiliation{Universit\`a degli Studi di Padova and Universit{\'e}
Paris Diderot}
\address[A]{Dipartimento di Matematica Pura e Applicata\\
Universit\`a degli Studi di Padova\\
via Trieste 63\\
35121 Padova\\
Italy\\
\printead{e1}}
\address[B]{Universit{\'e} Paris Diderot (Paris 7)\\
Laboratoire de Probabilit{\'e}s\\
\quad et Mod\`eles Al\'eatoires\\
(CNRS U.M.R. 7599)\\
U.F.R. Math\'ematiques\\
Case 7012 (Site Chevaleret)\\
75205 Paris Cedex 13\\
France\\
\printead{e2}}
\end{aug}

% HISTORY:
\received{\smonth{7} \syear{2009}}
\revised{\smonth{1} \syear{2010}}

% ABSTRACT
%
\begin{abstract}
A copolymer is a chain of repetitive units (\textit{monomers}) that
are almost identical, but they differ in their degree of affinity for
certain solvents. This difference leads to striking phenomena when the
polymer fluctuates in a nonhomogeneous medium, for example, made of two
solvents separated by an interface. One may observe, for instance, the
localization of the polymer at the interface between the two solvents.
A discrete model of such system, based on the simple symmetric random
walk on $\Z$, has been investigated in~\cite{cfBdH}, notably in the
weak polymer-solvent coupling limit, where the convergence of the
discrete model toward a continuum model, based on Brownian motion, has
been established. This result is remarkable because it strongly
suggests a universal feature of copolymer models. In this work, we
prove that this is indeed the case. More precisely, we determine the
weak coupling limit for a general class of discrete copolymer models,
obtaining as limits a one-parameter [$\ga\in(0,1)$] family of
continuum models, based on $\ga$-stable regenerative sets.
\end{abstract}

% KEYWORDS
%

\end{frontmatter}

%s1 ###
\section{Introduction}

%s1.1 ###
\subsection{The discrete model}
\label{sec:discmodel}
\begin{figure}

\includegraphics{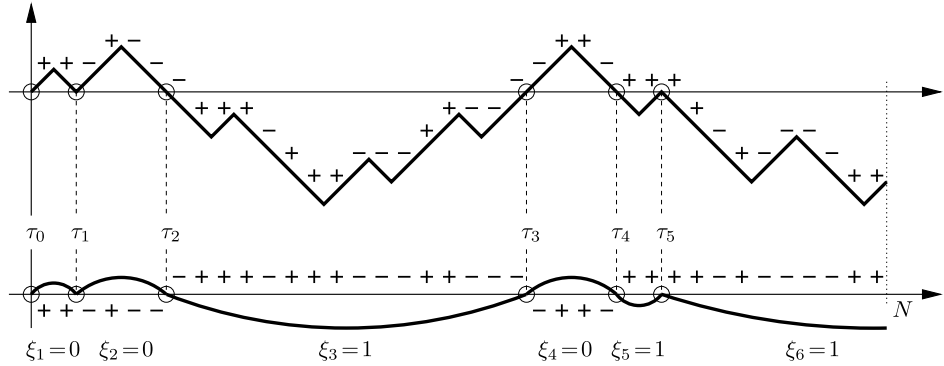}

\caption{The polymer model we deal with has been
introduced in the mathematical literature (see, e.g.,
\protect\cite{cfSinai}) as a modification of the law of the simple symmetric
random walk $\{S_n\}_{n \ge0}$ on $\Z$, with a density proportional to
$\exp[ \gl\sum_{n=1}^N (\go_n +h) \sign(S_n)]$ (\textit{Boltzmann
factor}). Each bond $(S_{n-1}, S_n)$ is interpreted as a
\textit{monomer} and by definition $\sign(S_n)$ is the sign of
$(S_{n-1}, S_n)$, that is, it is $+1$ (resp., $-1$) if the monomer
$(S_{n-1}, S_n)$ lies in the upper (resp., lower) half plane. In a
quicker way, $\sign(S_n)$ is just the sign of $S_{n-1}+S_{n}$. The
Boltzmann factor is somewhat different from the one appearing in
(\protect\ref{eq:SRWmodel}), but this is not a problem: in fact
$\gl\sum_{n=1}^N (\go_n +h) \sign(S_n)$ can be rewritten as $-2
\gl\sum_{n=1}^N \gD(S_{n-1}+S_{n}) (\go_n +h) + c_N$, where $c_N :=
\gl\sum_{n=1}^N (\go_n +h)$ does not depend on $S$, therefore the
quenched probability $\bP_{N, \go}$ is not affected by such a change.
It is clear that the trajectories of the walk, that are interpreted as
configurations of a \textit{polymer chain}, have an \textit{energetic
gain} (i.e., a~larger Boltzmann factor) if \textit{positively charged}
monomers [$(\go_n + h) > 0$] lie in the upper half plane and
\textit{negatively charged} ones [$(\go_n + h) < 0$] lie in the lower
one. The fulfillment of this requirement, even if only in a partial
way, entails however an \textit{entropic loss}: in fact the
trajectories have to stick very close to the horizontal axis (the
\textit{interface}) and there are only few such random walk
trajectories. The issue is precisely to understand who is the winner in
this \textit{energy-entropy competition}. The lower part of the figure
stresses the fact that the Boltzmann factor does not depend on the full
trajectory $S$, but only on the lengths and the signs of the successive
excursions, described by the variables $\tau, \xi$. In the figure, it
is also represented an example of the sequence of charges attached to
the copolymer, in the binary case ($\go_n \in\{ -1, +1\}$).}
\label{fig:cop-fig}
\end{figure}

Let $S:=\{ S_n\}_{n=0,1,\ldots}$ be the simple symmetric random walk
on $\Z$,
that is, $S_0=0$
and $\{ S_{n+1}-S_n\}_{n=0,1,\ldots}$ is an i.i.d. sequence
of random variables, each taking values $+1$ or $-1$ with probability $1/2$.
If $\bP$ is the law of $S$, we introduce
a new probability measure $\bP_{N, \go} = \bP_{N, \go, \gl, h}$
on the random walk trajectories
defined by
%
%e1.1 ###
%
\begin{equation}
\label{eq:SRWmodel}
\frac{\dd\bP_{N, \go}}{\dd\bP} (S) :=
\frac1{Z_{N, \go}}
\exp\Biggl(-2 \gl\sum_{n=1}^N \gD(S_{n-1}+S_n ) (\go_n +h)
\Biggr),
\end{equation}
where $N \in\N:= \{1, 2, \ldots\}$,
$\gl, h \in[0,\infty)$, we have set
$\gD(\cdot) := \ind_{(-\infty,0)} (\cdot)$ and
$\go:=\{\go_n \}_{n\in\N}$ is a sequence of real numbers.
Of course $Z_{N, \go}=Z_{N, \go, \gl, h}$ is the normalization constant,
called \textit{partition function} and given by
%
%e1.2 ###
%
\begin{equation} \label{eq:SRWZ}
Z_{N, \go} := \bE\Biggl[
\exp\Biggl(-2 \gl\sum_{n=1}^N \gD(S_{n-1}+S_n )
(\go_n +h) \Biggr) \Biggr] .
\end{equation}
We could have used $\gD(S_n)$ instead of $\gD(S_{n-1}+S_n)$, but
this apparently unnatural choice actually has a nice interpretation,
explained in the caption of Figure \ref{fig:cop-fig}.

%
%f1 ###
%

We are interested in the case when $\go$, called the sequence of \textit
{charges},
is chosen as a typical realization of
an i.i.d. sequence (call $\bbP$ its law). We assume that $\go$ and
$S$ are independent,
so that the relevant underlying law is $\bbP\otimes\bP$, but in
reality we are interested in \textit{quenched results}, that is,
we study $\bP_{N, \go}$ (in the limit $N \to\infty$)
for a \textit{fixed choice}
of $\go$. In the literature, the charge distribution
is often chosen Gaussian or of binary type, for example, $\bbP(\go_1=+1)=
\bbP(\go_1=-1)=1/2$.
We invite the reader to look at Figure \ref{fig:cop-fig} in order
to have a quick intuitive view of what this model describes (a polymer model).

Figure \ref{fig:cop-fig} also schematizes an aspect
of the model which is particularly relevant to us. Namely
that the \textit{Hamiltonian} of the model, that is,
the quantity appearing at the exponent in the right-hand side of
(\ref{eq:SRWmodel}), does not depend on the full
trajectories of $S$, but only on the
\textit{random set} $\tau:= \{n \in\N\cup\{0\}\dvtx
S_n=0\}$ (that we may also look at as an increasing \textit{random sequence}
$\tau=: \{\tau_0, \tau_1, \tau_2, \ldots\}$) and on the \textit{signs}
$\xi= \{\xi_j\}_{j\in\N}$, defined by $\xi_j := \Delta(S_n)$
for $n \in\{ \tau_j +1, \tau_{j+1}-1\}$ (i.e., $\xi_j = 0$
or~$1$ if the $j$th excursion of $S$ is positive or negative).
In fact, it is easily seen that $\Delta(S_{n-1}+S_n) =
\sum_{j=1}^\infty\xi_j \ind_{(\tau_{j-1}, \tau_j]}(n)$ is a function
of $\tau$ and $\xi$ only, and this suffices to reconstruct the Hamiltonian
[see (\ref{eq:SRWmodel})].
Note that we call the variables $\xi_n$ \textit{signs} even
if they take the values $\{0,1\}$ instead of $\{+1, -1\}$.

Under the simple random walk law $\bP$,
the two random sequences $\tau$ and $\xi$ are independent.
Moreover, $\xi$ is just an IID sequence of $B(1/2)$
(i.e., Bernoulli of parameter $1/2$) variables,
while $\tau$ is a \textit{renewal process}, that is, $\tau_0=0$ and
$\{ \tau_{j}-\tau_{j-1} \}_{j \in\N}$ is i.i.d.
Let us also point out that for every $j \in\N$,
%
%e1.3 ###
%
\begin{equation} \label{eq:tauSRW}
\bP( \tau_{j}-\tau_{j-1} =2n )
= \bP( \tau_{1}=2n ) \stackrel{n \to\infty}{\sim}
\frac1{2\sqrt{\pi} n^{3/2}},
\end{equation}
where we have introduced the notation
$f(x) \sim g(x)$ for \mbox{$\lim_{x\to\infty} f(x)/g(x)=1$}
[in the sequel, we will also use $\sim$ to denote
equality in law: e.g., $\go_1 \sim\go_2\sim\cN(0,1)$].

This discussion suggests a generalized framework
in which to work, that has been already introduced
in \cite{cfBGLT,cfBook}. We start from scratch: let us consider a
general renewal process $\tau= \{\tau_n\}_{n\ge0}$
on the nonnegative integers $\N\cup\{0\}$ such that
%
%e1.4 ###
%
\begin{equation} \label{eq:K}
K(n) := \bP( \tau_1 =n )
\stackrel{n\to\infty}{\sim} \frac{L(n)}{n^{1+\ga}},
\end{equation}
where $\ga\ge0$ and $L\dvtx(0,\infty) \to(0, \infty)$
a \textit{slowly varying function}, that is, a (strictly)
positive measurable function
such that $\lim_{x \to\infty}L(cx)/L(x)=1$, for every $c>0$
(see Remark \ref{rem:SVF} below for more details).
We assume that $\tau$ is a \textit{persistent} renewal,
that is,
%also that $K(n) > 0$ for every $n\in\N$ and
$\bP(\tau_1 < \infty) = \sum_{n \in\N} K(n)=1$, which is
equivalent to $\bP(\vert\tau\vert=\infty)=1$,
where $\vert\tau\vert$ denotes the cardinality of $\tau$,
viewed as a (random) subset of $\N\cup\{0\}$.
We will switch freely from looking at $\tau$
as a sequence of random variables or as a random set.

Let $\xi= \{\xi_n\}_{n\in\N}$ denote an i.i.d. sequence
of $B(1/2)$ variables, independent of~$\tau$, that we still
call signs.
With the couple $(\tau,\xi)$ in our hands, we
build a new sequence $\gD= \{\gD_n\}_{n\in\N}$ by setting
$\gD_n = \sum_{j=1}^\infty\xi_j \ind_{(\tau_{j-1}, \tau_j]}(n)$,
in analogy with the simple random walk case.
In words, the signs $\Delta_n$ are constant between the epochs of
$\tau$
and they are determined by $\xi$.

We are now ready to introduce the general \textit{discrete copolymer model},
as the probability law $\bP_{N, \go} = \bP_{N, \go}^{\gl, h}$ for
the sequence $\Delta$ defined by
%
%e1.5 ###
%
\begin{equation}
\label{eq:taumodel}
\frac{\dd\bP_{N, \go}}{\dd\bP} (\Delta) :=
\frac1{Z_{N, \go}}
\exp\Biggl(-2 \gl\sum_{n=1}^N \gD_n (\go_n +h)
\Biggr) ,
\end{equation}
where $N \in\N$, $\gl, h \in[0,\infty)$ and $\go= \{\go_n\}
_{n\in\N}$
is a sequence of real numbers (a~typical realization of
an i.i.d. sequence, see below).
The partition function
$Z_{N,\go} = Z_{N,\go}^{\gl, h}$ is given by
%
%e1.6 ###
%
\begin{equation} \label{eq:discZ}
Z_{N, \go} := \bE\Biggl[
\exp\Biggl(-2 \gl\sum_{n=1}^N \gD_n (\go_n +h) \Biggr) \Biggr] .
\end{equation}
In order to emphasize the value of $\ga$ in (\ref{eq:K}),
we will sometimes speak of a \textit{discrete $\ga$-copolymer model}, but
we stress
that $\bP_{N, \go}$ depends on the full law $K(\cdot)$.

Note that the new model (\ref{eq:taumodel}) only describes the
sequence of signs $\Delta$, while the
\textit{simple random walk model} (\ref{eq:SRWmodel})
records the full trajectory $S$.
However, once we project the probability law (\ref{eq:SRWmodel})
on the variables $\Delta_n := \Delta(S_{n-1} + S_n)$,
it is easy to check that the simple random walk model
becomes a particular case of (\ref{eq:taumodel})
and its partition function (\ref{eq:SRWZ}) coincides with
the general one given by (\ref{eq:discZ}), provided we choose $K(\cdot
)$ as
the law of the first return to zero of the simple random walk
[corresponding to $\alpha= \frac12$; see
(\ref{eq:tauSRW}) and (\ref{eq:K})]. As a matter of fact,
since we require that $K(n) > 0$ for all large
$n\in\N$ [cf. (\ref{eq:K})], strictly speaking
the case of the simple random walk is not covered.
We stress, however, that
our arguments can be adapted in a straightforward way to treat the
cases in
which there exists a positive
integer $T$ such that $K(n)=0$ if $n/T \notin\N$ and
relation (\ref{eq:K}) holds restricting $n \in T \N$
(of course $T=2$ for the simple random walk case).

To complete the definition of the discrete copolymer model,
let us state precisely our hypotheses on the \textit{disorder} variables
$\go= \{\go_n\}_{n\in\N}$. We assume that the sequence $\go$ is
i.i.d. and that $\go_1$ has locally finite exponential moments, that
is,
there exists $t_0 > 0$ such that
%
%e1.7 ###
%
\begin{equation}
\label{eq:Mgo}
\M(t) := \bbE[\exp(t \go_1)] < \infty\qquad
\mbox{for every } t \in[-t_0, t_0] .
\end{equation}
We also fix
%
%e1.8 ###
%
\begin{equation}
\label{eq:gonorm}
\bbE[\go_1 ] = 0 \quad\mbox{and}\quad \bbE[\go_1^2 ] = 1 ,
\end{equation}
which entails no loss
of generality (it suffices to shift $\gl$ and $h$).
In particular, these assumptions guarantee that there exists
$c_0>0$ such that
%
%e1.9 ###
%
\begin{equation}
\label{eq:omcond}
\max_{t \in[-t_0,t_0]} \M(t) \le
\exp(c_0 t^2 ).
\end{equation}

Although it only keeps track of the sequence of signs $\Delta$, we
still interpret
the probability law $\bP_{N,\go}$ defined in (\ref{eq:taumodel}) as
a model
for an \textit{inhomogeneous polymer} (this is the meaning
of \textit{copolymer}) that interacts with
two selective solvents (the upper and lower half planes) separated by
a flat interface (the horizontal axis), as it is explained in the caption
of Figure \ref{fig:cop-fig}. In particular, $\Delta_n = 0$ (resp.,
$1$) means that
the $n$th monomer of the chain lies above (resp., below) the interface.
To reinforce the intuition, we will sometimes describe the model in
terms of full
trajectories, like in the simple random walk case.
\begin{rem}
\label{rem:SVF}
We refer to \cite{cfBinGolTeu} for a full account on slowly
varying functions. Here, we just recall that the asymptotic behavior of
$L(\cdot)$ is \textit{weaker than any power}, in the sense that, as $x
\to\infty$,
$L(x) x^a$ tends to $\infty$
for $a>0$ and to zero if $a<0$. The most basic example of a slowly
varying function
is any positive measurable function that converges to
a positive constant at infinity (in this case, we say that the slowly
varying function is trivial). Other important examples are
positive measurable functions which behave asymptotically like
the power of a logarithm, that is,
$L(x) \sim\log(1+x)^{a}$, $a \in\R$.
%
% [FRA] Do we need this?
%
%Moreover, by \cite[Prop. 1.3.4]{cfBinGolTeu}, in our context
%there is no loss of generality in assuming $L(\cdot)$ to be $C^\infty$.
\end{rem}

%%%%%%%%%%%%%%%%%%%%%%%%%%%%%%%%%%%%%
%s1.2 ###
\subsection{The free energy: Localization and delocalization}

This work focuses on the properties of the free energy of the discrete
copolymer, defined by
%
%e1.10 ###
%
\begin{equation} \label{eq:fe-qa}\quad
\tf(\gl, h) := \lim_{N \to\infty} \tf_N(\gl, h)\qquad
\mbox{where } \tf_N(\gl, h) :=
\frac1N \bbE[ \log Z_{N, \go} ] .
\end{equation}
The existence of such a limit follows by a standard argument, see, for example,
\cite{cfBook}, Chapter 4, where it is also proven that for every $\gl$
and $h$
%
%e1.11 ###
%
\begin{equation}
\label{eq:fe-q}
\tf(\gl, h) = \lim_{N \to\infty}
\frac1N \log Z_{N, \go} ,\qquad
\mbox{$\bbP(\dd\go)$-a.s. and in $L^1(\bbP)$} .
\end{equation}
Equations (\ref{eq:fe-qa}) and (\ref{eq:fe-q})
are telling us that the limit in (\ref{eq:fe-q}) does not depend
on the (typical) realization of $\go$. Nonetheless, it is worthwhile
to stress that it does depend
on $\bbP$, that is, on the law of $\go_1$, as well as
on the renewal process on which the model
is built, namely on the inter-arrival law $K(\cdot)$.
This should be kept in mind, even if we omit
$\bbP$ and $K(\cdot)$ from the notation $\tf(\gl, h)$.

\begin{figure}

\includegraphics{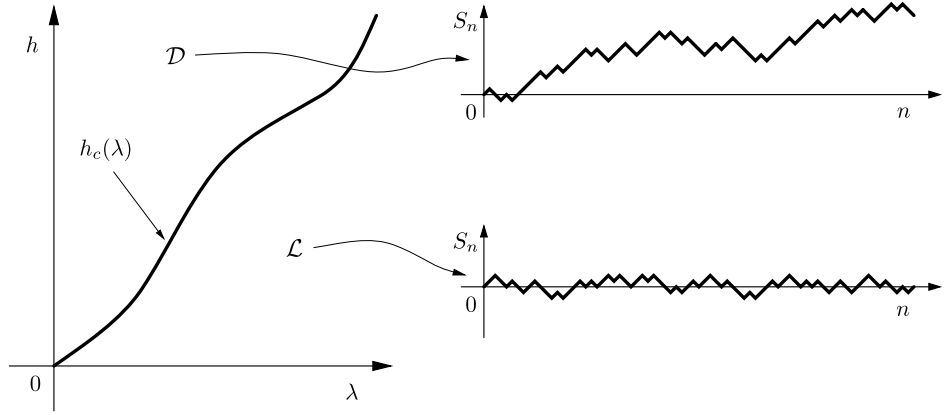}

\caption{In the figure, on the left, a sketch of the phase diagram of
the discrete
copolymer model. The critical curve $\gl\mapsto h_c(\gl)$ separates
the localized regime $\cL$ from the delocalized one $\cD$.
This is a free energy characterization of the notion of (de)localization,
but this characterization does correspond the to sharply different
path behaviors, sketched on the right side of the figure.
In a nutshell, if $(\gl,h)\in\cL$ then, for $N\to\infty$, the typical
paths intersect the interface (the horizontal axis) with a positive density,
while in the interior of $\cD$ the path strongly prefers not to
enter the lower half plane. In this work, we just
focus on properties of the free energy and for details on the link with path
properties, including a review of the literature and open problems;
we refer to \protect\cite{cfBook}, Chapters 7 and~8.}
\label{fig:phasediag}
\end{figure}

An elementary, but crucial observation is
%
%e1.12 ###
%
\begin{equation} \label{eq:ebutc}
\tf(\gl,h) \ge 0 \qquad \forall\gl, h \ge0 .
\end{equation}
This follows simply by restricting the expectation
in (\ref{eq:discZ}) to the event $\{\tau_1 > N, \xi_1 = 0\}$,
on which we have $\Delta_1=0, \ldots, \Delta_N = 0$, hence we obtain
$Z_{N, \go} \ge\frac12 \bP( \tau_1>N)$ and it suffices to
observe that $N^{-1}\log\bP( \tau_1>N) $ vanishes as $N \to\infty$,
thanks to (\ref{eq:K}). Notice that the event
$\{\tau_1 > N, \xi_1 = 0\}$ corresponds to the set of trajectories
that never visit the lower half plane, therefore the right-hand side
of (\ref{eq:ebutc}) may be viewed as the contribution to the free energy
given by these trajectories. Based on this observation,
it is customary to say that $(\gl,h) \in\cD$ (\textit{delocalized regime})
if $\tf(\gl,h)=0$, while $(\gl,h) \in\cL$ (\textit{localized regime})
if $\tf(\gl,h)>0$ (see also Figure \ref{fig:phasediag} and its caption).

%
%f2 ###
%

We have the following theorem.
\begin{theorem}
\label{th:phasediag}
If we set $h_c(\gl):= \sup\{h\dvtx \tf(\gl, h) >0\}$, then
$h_c(\gl)=\inf\{h\dvtx\break \tf(\gl, h) =0\}$ and the function
$h_c\dvtx [0, \infty) \to[0, \infty]$ is
strictly increasing and continuous
as long as it is finite. Moreover, we have the
explicit bounds
%
%e1.13 ###
%
\begin{equation}
\label{eq:LBUBhc}
\frac{1}{2\gl/(1+\ga)}\log\M\bigl( -2\gl/(1+\ga) \bigr)
\le h_c(\gl) \le \frac{1}{2\gl}\log\M( -2\gl
),
\end{equation}
where the left inequality is strict when $\ga\ge0.801$ (at least for
$\gl$ small)
and the right inequality is strict as soon as $\ga>0$ [for every
$\gl< \sup\{t\dvtx \log\M( -2t) < \infty\}$].
\end{theorem}

The first part of Theorem \ref{th:phasediag} is proven in
\cite{cfBdH} and \cite{cfBG} (see also \cite{cfBook}, Chapter~6). In
\cite{cfBG}, one also finds
the quantitative estimates (\ref{eq:LBUBhc}), except for the strict
inequalities proven in \cite{cfBGLT} (see also \cite{cfTAAP}).
From (\ref{eq:LBUBhc}), one directly extracts
%
%e1.14 ###
%
\begin{equation}
\frac1{1+\ga} \le \liminf_{\gl\searrow0}
\frac{h_c(\gl)}{\gl} \le \limsup_{\gl\searrow0}
\frac{h_c(\gl)}{\gl} \le 1 \qquad \forall\ga\ge0 .
\end{equation}
For $\ga> 0$, this result has been sharpened to
%
%e1.15 ###
%
\begin{equation}
\label{eq:slopebounds}
\max\biggl(\frac12, \frac{g(\ga)}{\sqrt{1+\ga}},\frac1{1+\ga
} \biggr) \le \liminf_{\gl\searrow0}
\frac{h_c(\gl)}{\gl} \le \limsup_{\gl\searrow0}
\frac{h_c(\gl)}{\gl} < 1 ,
\end{equation}
where $g(\cdot)$ is a continuous function such that
$g(\ga)=1$ for $\ga\ge1$ and for which one can show
that $g(\ga)/\sqrt{1+\ga} > 1/(1+\ga)$ for $\ga\ge0.801$
[by evaluating $g(\cdot)$ numerically one can go down to $\ga\ge0.65$].
In particular, the lower bound in (\ref{eq:slopebounds})
reduces to $1/2$ for $\ga\ge3$ and to
$1/\sqrt{1+\ga}$ for $\ga\in[1,3]$.
The bounds in (\ref{eq:slopebounds})
are proven in \cite{cfBGLT} and \cite{cfTcg}. We invite the reader
to look again at Figure \ref{fig:phasediag}.
We also point out that a numerical
study of the critical line in the simple random walk case ($\alpha =
\frac{1}{2}$) has been performed in \cite{cfCGG}, while the critical point of a
simplified copolymer model, originally introduced in \cite{cfBG}, has been
obtained in \cite{cfBCT}.

The focus on the behavior of the \textit{critical line} $h_c(\gl)$ for
$\gl$ small
has a reason, that is, at the heart of this paper:
our aim is to study the free energy $\tf(\gl,h)$ of discrete
copolymer models
in the \textit{weak coupling limit}, that is, when $\gl$
and $h$ are small.
We will show that the behavior of $\tf(\gl,h)$ in this regime
is captured by the exponent $\ga$ appearing in (\ref{eq:K}),
independently of the finer details of the inter-arrival law $K(\cdot)$.
In particular, we prove that $h_c^\prime(0)$
exists and that it depends only on $\ga$. In order
to state these results precisely, we need to introduce
a class of copolymer models in the continuum: in a suitable sense,
they capture the limit of discrete copolymer models as $\gl,h \searrow0$.

%%%%%%%%%%%%%%%%%%%%%%%%%%%%%%%%%%%%%%%%%%%%%%%%%%%%%%%%%%%%%%%%%%%%

%s1.3 ###
\subsection{The continuum model: Brownian case}

Bolthausen and den Hollander introduced in \cite{cfBdH}
the \textit{Brownian copolymer model}, whose partition function is
given by
%
%e1.16 ###
%
\begin{equation}
\label{eq:ZBrownian}
\widetilde Z_{t,\gb}^{\mathrm{BM}} := \bE
\biggl[ \exp\biggl(
- 2 \gl\int_0^t \Delta( \widetilde B (u) ) \bigl(\dd\gb
(u) + h \,\dd u \bigr)
\biggr) \biggr],
\end{equation}
where once again $\gl, h\ge0$, $\Delta(x) := \ind_{(-\infty
,0)}(x)$ and
$\widetilde B(\cdot)$ (the polymer), $\gb(\cdot)$ (the medium)
are independent standard Brownian motions
with laws $\bP$ and $\bbP$, respectively.
%The remarkable result in \cite{cfBdH} is the link with
%the small $\gl$ behavior of the discrete (random walk based)
%copolymer model: we are going to recall it here in a rather sketchy way
%(more details will come later, since the result in \cite{cfBdH}
%is a particular case
%of our main result).

The corresponding free energy $\widetilde\tf_{\mathrm{BM}} (\gl, h)$
is defined as the limit as $t \to\infty$ of
$\frac1t \bbE[ \log\widetilde Z_{t,\gb}^{\mathrm{BM}}]$
and one has $\widetilde\tf_{\mathrm{BM}}(\gl, h)\ge0$ for every
$\gl, h \ge0$,
in analogy with the discrete case.
Therefore, by looking at the positivity of $\widetilde\tf_{\mathrm
{BM}}$, one can
define also for the Brownian copolymer model
the localized and delocalized regimes, that are separated
by the critical line $\widetilde h_c(\gl):=
\sup\{ h\dvtx \widetilde\tf_{\mathrm{BM}}(\gl, h)>0\}$. Now a real
novelty comes
into the game: the scaling properties of the two Brownian motions
yield easily that for every $a>0$
%
%e1.17 ###
%
\begin{equation}
\label{eq:scalingBM}
\frac1{a^2}
\widetilde\tf_{\mathrm{BM}} (a\gl, ah) = \widetilde\tf
_{\mathrm{BM}}(\gl, h) .
\end{equation}
In particular, the critical line is a straight line:
$\widetilde h_c (\gl)= \slope_{\mathrm{BM}} \gl$, for every $\gl
\ge0$, with
%
%e1.18 ###
%
\begin{equation}
\label{eq:mBM}
\slope_{\mathrm{BM}} := \sup\{c \ge0\dvtx \widetilde\tf
_{\mathrm{BM}}(1,c)>0 \}.
\end{equation}

We are now ready to state the main result in \cite{cfBdH}.
\begin{theorem}
\label{th:BdH}
For the simple random walk model (\ref{eq:SRWmodel}),
with $\go_1$ such that $\bbP( \go_1=+1)=\bbP( \go_1=-1)=1/2$,
we have
%
%e1.19 ###
%
\begin{equation}
\label{eq:BdH0}
\lim_{a \searrow0} \frac1{a^2}\tf(a \gl, ah) =
\widetilde\tf_{\mathrm{BM}} (\gl, h) \qquad \forall\gl,h \ge
0 ,
\end{equation}
and
%
%e1.20 ###
%
\begin{equation}
\label{eq:BdH1}
\lim_{\gl\searrow0} \frac{h_c(\gl)}{\gl} =
\slope_{\mathrm{BM}} \in(0,1].
\end{equation}
\end{theorem}

The great interest of this result is that
it provides a precise formulation for the fact that the Brownian
copolymer model is
the weak coupling scaling limit of the
simple random walk copolymer model (\ref{eq:SRWmodel}).
At the same time, the fact that such a result is proven
only for the simple random walk model and only for a single choice
of the charges distribution appears to be a limitation.
In fact, since Brownian motion is the scaling limit of many
discrete processes, it is natural to guess that
Theorem \ref{th:BdH} should hold for a large class of discrete
copolymer models
and for a vast choice of charge distributions (remaining of course
in the domain of attraction of the Gaussian law
and adding some technical assumptions). This
would show that the Brownian copolymer model has indeed
a \textit{universal} character.

In fact, Theorem \ref{th:BdH} has been generalized in \cite{cfGT} to
a large class of
disorder random variables (including all bounded random variables).
A further generalization has been obtained in \cite{cfNicolas},
in the case when, added to the copolymer
interaction, there is also a \textit{pinning} interaction at the
interface, that is,
an energy reward in touching the interface. We stress, however,
that these generalizations are always for the copolymer model
built over the simple random walk: going beyond the simple random walk case
appears indeed to be a very delicate (albeit natural) step.

The main result of this paper is that Theorem \ref{th:BdH} can be generalized
to any discrete $\ga$-copolymer model with
$\ga\in(0,1)$ and to any disorder distribution
satisfying (\ref{eq:Mgo}) and (\ref{eq:gonorm}) (see Theorem \ref
{th:main} below).
For $\ga= \frac12$, the scaling limit
is precisely the Brownian copolymer model (\ref{eq:ZBrownian}),
like in the simple random walk case,
while for $\ga\ne\frac12$
the continuum copolymer model is defined in the next
subsection. We stress from now that the scaling limit
\textit{depends only on $\ga$}: in particular,
there is no dependence on the slowly varying function $L(\cdot)$
appearing in (\ref{eq:K})
and no dependence on $\bP(\tau_1=n)$ for any finite $n$.

%%%%%%%%%%%%%%%%%%%%%%%%%%%%%%%%%%%%%%%%%%%%%%%%%%%%%%%%%%%%%%%%%

%s1.4 ###
\subsection{The continuum $\ga$-copolymer model}
\label{sec:genC0}

Let us start by recalling that, for \mbox{$\delta\ge0$}, \textit{the square of
$\delta$-dimensional
Bessel process} (started at $0$) is the process $X = \{X_t\}_{t\ge0}$
with values
in $[0,\infty)$, that is, the unique strong solution of the following equation:
%
%e1.21 ###
%
\begin{equation}
X_t = 2 \int_0^t \sqrt{X_s} \,\dd w_s + \delta t ,
\end{equation}
where $\{w_t\}_{t \ge0}$ is a standard Brownian motion.
The \textit{$\delta$-dimensional Bessel
process} is by definition the process
$Y = \{Y_t := \sqrt{X_t}\}_{t\ge0}$:
it is a Markov process on $[0,\infty)$
that enjoys the standard Brownian scaling (\cite{cfRevYor}, Chapter XI,
Proposition (1.10)).
We focus on the case $\delta\in(0,2)$,
when a.s. the process $Y$ visits the origin infinitely
many times \cite{cfRevYor}, Chapter XI, Proposition (1.5).
We actually use
the parametrization $\delta= 2(1-\ga)$ and we then restrict to $\ga
\in(0,1)$.

It is easily checked using It\^o's formula that for $\ga= \frac12$
(i.e., $\delta= 1$) the process $Y$ has the same law as
the absolute value of Brownian motion on $\R$. Since
to define the Brownian copolymer model (\ref{eq:ZBrownian})
we have used the full Brownian motion process, not only its absolute value,
we need a modification of the Bessel process in which
each excursion from zero may be either positive or negative,
with the sign chosen \textit{by fair coin tossing}. Such a process, that
we denote by $\widetilde B^\ga:=\{\widetilde B^\ga(t)\}_{t\ge0}$,
has been
considered in the literature, for example, in \cite{cfBPY}, and
is called \textit{Walsh process of index} $\ga$
(in \cite{cfBPY} a more general case is actually considered:
in their notation, our process corresponds to the choices $k=2$,
$E_1=[0, \infty)$,
$E_2=(-\infty,0]$ and $p_1=p_2=1/2$). It is easy to see that
the process $\widetilde B^\ga$ inherits the Brownian scaling. We
denote by $\bP$ its law.

We are now ready to generalize the Brownian copolymer model (\ref
{eq:ZBrownian}):
given $\ga\in(0,1)$, we define the partition function of the
\textit{continuum $\ga$-copolymer model} through the formula
%
%e1.22 ###
%
\begin{equation}
\label{eq:ZWalsh}
\widetilde Z_{t,\gb}^{\ga} := \bE\exp\biggl(
- 2 \gl\int_0^t \Delta( \widetilde B^\ga(u) )
\bigl(\dd\gb(u) + h \,\dd u \bigr) \biggr) ,
\end{equation}
where $\gb= \{\gb(t)\}_{t\ge0}$ always denotes a standard Brownian motion
with law $\bbP$, independent of $\widetilde B^\ga$, and $\Delta(x) =
\ind_{(-\infty,0)}(x)$.
Since for $\ga= \frac12$ the process $\widetilde B^{1/2}$ is just a
standard Brownian motion, $\widetilde Z_{t,\gb}^{1/2}$ coincides with
$\widetilde Z_{t,\gb}^{\mathrm{BM}}$ defined in (\ref{eq:ZBrownian}).
For the sake of simplicity, in (\ref{eq:ZWalsh})
we have only defined the partition function
of the continuum $\ga$-copolymer model: of course, one can easily introduce
the corresponding probability measure $\bP_{t,\gb}$ on the
paths of $\widetilde B^\ga$, in analogy with the discrete case, but we
will not need it.

Let us stress that the integral in (\ref{eq:ZWalsh}),
as well as the one in (\ref{eq:ZBrownian}), does not really depend
on the full path of the process $\widetilde B^\ga$; in fact, being a function
of $\Delta(\widetilde B^\ga(\cdot))$, it only matters to know, for
every $u \in[0,t]$,
whether $B(u) < 0$ or $B(u) \ge0$.
For this reason, it is natural to introduce (much like in the discrete case)
the \textit{zero level set} $\widetilde\tau^\ga$ of $\widetilde B^\ga
(\cdot)$:
%
%e1.23 ###
%
\begin{equation}
\label{eq:tautilde}
\widetilde\tau^\ga := \{ s \in[0,\infty)\dvtx \widetilde
B^\ga(s)=0 \} .
\end{equation}
The set $\widetilde\tau^\ga$ contains \textit{almost} all the
information we need,
because, conditionally on~$\widetilde\tau^\ga$, the sign of
$\widetilde B^\ga$ inside each excursion
is chosen just by tossing an independent fair coin.
Moreover, the random set $\widetilde\tau^\ga$ is a much studied object:
it is, in fact, the \textit{$\ga$-stable regenerative set}
(\cite{cfRevYor}, Chapter XI, Exercise (1.25)). Regenerative sets may be viewed
as the continuum analogues of renewal processes: we discuss them
in some detail in Section \ref{sec:C0},
also because it will come very handy to restate the model
in terms of regenerative sets for the proofs.

The free energy for the continuum $\ga$-copolymer model
is defined in close analogy to the discrete case, but proving its existence
turns out to be a highly nontrivial task. For this reason, we state
it as a result in its own.
\begin{theorem}
\label{th:existence}
The limit of $\frac1t
\bbE[ \log\widetilde Z_{t,\gb}^{\ga} ]$ as $t\to\infty$
exists and we call it $\widetilde\tf_\ga(\gl, h)$.
For all $\ga\in(0,1)$
and $\gl, h \in[0,\infty)$ we have $0 \le\widetilde\tf_\ga(\gl,
h) < \infty$
and furthermore
%
%e1.24 ###
%
\begin{equation}\label{eq:existence}
\lim_{t \to\infty} \frac1t \log\widetilde Z_{t,\gb}^{\ga}
= \widetilde\tf_\ga(\gl, h) ,
\end{equation}
both $\bbP(\dd\gb)$-a.s. and in $L^1(\bbP)$.
The function $(\gl, h) \mapsto\tf_\ga(\gl, h)$ is continuous.
\end{theorem}

Like before, the nonnegativity of the free energy
leads to exploiting the dichotomy $ \widetilde\tf_\ga(\gl, h)=0$
and $\widetilde\tf_\ga(\gl, h)>0$ in order to define, respectively, the
delocalized and localized regimes of the continuum
$\ga$-copolymer model. The monotonicity
of $\widetilde\tf_\ga( \gl, \cdot)$ guarantees that if we set
$\widetilde h_c^{\ga}(\gl):= \sup\{ h \ge0\dvtx \widetilde\tf_\ga
(\gl, h)>0\}$, then
we also have
$\widetilde h_c^{\ga}(\gl):= \inf\{ h \ge0\dvtx \widetilde\tf_\ga
(\gl, h)=0\}$.
Moreover, the scaling properties of $\gb$ and $\widetilde B^\ga$
imply that
(\ref{eq:scalingBM}) holds unchanged for $\widetilde\tf_\ga(\cdot,
\cdot)$
so that the critical line is again a straight line:
$\widetilde h_c^{\ga} (\gl)= \slope_\ga\gl$ for every $\gl\ge0$, with
%
%e1.25 ###
%
\begin{equation}
\label{eq:mga}
\slope_\ga := \sup\{c \ge0\dvtx \widetilde\tf_\ga
(1,c)>0 \},
\end{equation}
in direct analogy with (\ref{eq:mBM}). Plainly, $\slope_{1/2} =
\slope_{\mathrm{BM}}$.

%%%%%%%%%%%%%%%%%%%%%%%%%%%%%%%%%%%%%%%%%%%%%%%%%%%%%%%%%%%%%%%%%

%s1.5 ###
\subsection{The main result}

We can finally state the main result of this paper.

\begin{theorem}
\label{th:main}
Consider an arbitrary discrete $\ga$-copolymer model
satisfying the hypotheses (\ref{eq:K}), (\ref{eq:Mgo}) and (\ref{eq:gonorm}),
with $\ga\in(0,1)$. For all $\gl,h \ge0$, we have
%
%e1.26 ###
%
\begin{equation}
\label{eq:mainfe}
\lim_{a \searrow0} \frac1{a^2} \tf( a\gl, a h)
= \widetilde\tf_\ga(\gl, h) .
\end{equation}
Moreover,
%
%e1.27 ###
%
\begin{equation}
\label{eq:mainslope}
\lim_{\gl\searrow0} \frac{h_{c}(\gl)}{\gl} = \slope_\ga.
\end{equation}
\end{theorem}

Theorem \ref{th:main} shows that the continuum $\ga$-copolymer is the
universal
weak interaction limit of arbitrary discrete $\ga$-copolymer models.
Although the phase diagram of a discrete copolymer model does depend
on the details of the inter-arrival law~$K(\cdot)$, it nevertheless exhibits
universal features for weak coupling. In particular, the critical line
close to the origin
is, to leading order, a straight line of slope~$\slope_\ga$.
It is therefore clear that
computing $\slope_\ga$ or, at least, being able of improving the
known bounds on $\slope_\ga$
would mean a substantial progress in understanding the phase transition
in this model. Note that, of course, given (\ref{eq:mainslope}), the
bounds in
(\ref{eq:slopebounds}) are actually bounds on $\slope_\ga$
(and they represent the state of the art on this important issue, to
the the authors'
knowledge).

It is remarkable that in the
physical literature there is, on the one hand, quite some attention on
the slope
at the origin of the critical curve (see, e.g., \cite{cfGHLO})
but, on the
other hand, its \textit{universal} aspect has not been appreciated
(some of the physical predictions
are even in contradiction with the universality of the slope).
We refer to \cite{cfBGLT,cfBook,cfdH} for overviews of the extensive
physical literature on copolymer models.

We do not expect a generalization of Theorem \ref{th:main} to $\ga
\notin
(0,1)$. To be more precise, the
case $\ga=0$ is rather particular: the critical curve is known explicitly
by Theorem \ref{th:phasediag},
the slope at the origin is universal and its value is one.
The case $\ga=1$ with $\bE[ \tau_1]=\infty$ may still be treatable,
but the associated regenerative set is the full line, so Theorem \ref{th:main}
cannot hold as stated. An even more substantial problem
arises whenever $\bE[ \tau_1]<\infty$ (in particular, for every $\ga>1$):
apart from the fact that the regenerative set becomes trivial,
there is a priori no reason why universality should hold.
The \textit{rationale} behind Theorem \ref{th:main} is that at small coupling
the renewal trajectories are not much perturbed by the
interaction with the charges. If $\bE[ \tau_1]=\infty$, one may then
hope that
\textit{long} inter-arrival gaps dominate, as they do when there is
no interaction with the charges:
since the statistics of long gaps depends only on the tail of $K(\cdot)$
and within long gaps the disorder can be replaced by Gaussian disorder,
%and $K(\cdot)$ by its asymptotic behavior,
Theorem \ref{th:main} is plausible.
This is, of course, not at all the case if $\bE[ \tau_1]<\infty$.
\begin{rem}
One may imagine that (\ref{eq:mainslope}) is a consequence
of (\ref{eq:mainfe}), but this is not true. In fact, it is easy
to check that (\ref{eq:mainfe}) directly implies
%
%e1.28 ###
%
\begin{equation}
\label{eq:mainslopeinf}
\liminf_{\gl\searrow0} \frac{h_{c}(\gl)}{\gl} \ge \slope
_\ga,
\end{equation}
but the opposite bound (for the superior limit) does not follow
automatically. We obtain it as a corollary of
our main technical result (Theorem \ref{th:main-tech}).
\end{rem}

%%%%%%%%%%%%%%%%%%%%%%%%%%%%%%%%%%%%%%%%%%%%%%%%%%%%%%%%

%s1.6 ###
\subsection{Outline of the paper}
\label{sec:overview}

We start, in
Section \ref{sec:C0}, by taking a closer look at the continuum model
and by giving a proof of the existence of the free energy (Theorem \ref
{th:existence}).
Such an existence result had been overlooked in \cite{cfBdH}.
A proof was proposed in \cite{cfGB}, in the Brownian context,
giving for granted a suitable uniform boundedness property, that is not
straightforward (this is the issue addressed in Appendix \ref
{sec:contbound}).
The proof that we give here therefore generalizes [from $\ga=1/2$
to $\ga\in(0,1)$] and completes the proof in \cite{cfGB}.
We follow the general scheme of the proof in \cite{cfGB},
that is, we first define a suitably modified partition function,
that falls in the realm of Kingman's super-additive ergodic theorem
\cite{cfKin},
and then we show that such a modified partition function has the
same Laplace asymptotic behavior as the original one.
Roughly speaking, the modified partition function is obtained by
relaxing the condition that $\widetilde B^\ga(0)=0$; one takes, rather,
the infimum
over a finite interval of starting points.
If introducing such a modified partition function is a standard procedure,
a straightforward application of this idea does not seem to lead far.
Such an infimum procedure has to be set up in a careful
%and rather atypical
way
in order to be able to perform the second step of the proof,
that is, stepping back to the original partition function. With respect
to the
proof in \cite{cfGB}, that exploits the full path of
the Brownian motion $B(\cdot)$, the
one we present here is fully based on the regenerative set.
Overall, establishing the existence of the continuum free energy
is very much harder than the discrete counterpart case and it appears
to be remarkably subtle and complex when compared to the analogous statement
for \textit{close relatives} of our model (see, e.g., \cite{cfCY}).

In Section \ref{sec:main_result_proof}, we give
the proof of our main result, Theorem \ref{th:main},
following the scheme set forth in \cite{cfBdH}
(we refer to it as the \textit{original approach}), which
is based on a four step procedure.
We outline it here, in order to
give an overview of the proof and to stress the points
at which our arguments are more substantially novel.

\begin{enumerate}[(1)]
\item[(1)] \textit{Coarse graining of the renewal process.}
In this step, we replace the Boltzmann factor by a new, \textit{coarser}
one, which
does not depend on the \textit{short} excursions of the renewal
process (in the sense
that these excursions inherit the sign of a neighbor long gap).
This step is technically, but not
substantially different
from the one in the original approach.

\item[(2)] \textit{Switching to Gaussian charges.}
The original approach exploits the well known, and
highly nontrivial, coupling result due to Koml\'os, Major
and Tusn\'ady \cite{cfKMT}. We take
instead a more direct, and more elementary, approach.
In doing so, we get rid of any assumption, beyond local exponential
integrability, on the disorder.

\item[(3)] \textit{From the renewal process to the regenerative set.} This is probably
the crucial step. The original approach exploits heavily
the underlying simple random walk and the exact formulas
available for such a process. Our approach necessarily
sticks to the renewal process and, in a sense, the point is showing that
suitable local limit theorems (crucial here are results by Doney \cite{cfDon})
suffice to perform this step. There is, however, another issue that makes our
general case different from
the simple random walk case. In fact this step, in the original
approach, is based
on showing that a suitable Radon--Nikodym derivative,
comparing the renewal process and the regenerative set, is uniformly bounded.
In our general set-up, this Radon--Nikodym derivative is not bounded
and a more careful estimate has to be carried out.

\item[(4)] \textit{Inverse coarse graining of the regenerative set.}
We are now left with a model
based on the regenerative set, but depending only on the
\textit{large} excursions. We have therefore to show that
putting back the dependence on the small excursions does not
modify substantially the quantity we are dealing with. This is
parallel to the first step: it involves estimates that are different from
the ones in the original approach, because we are sticking
to the regenerative set formulation and because $\ga$ is not necessarily
equal to $1/2$, but the difference is, essentially, just technical.
\end{enumerate}

Let us finally mention that our
choice of focusing on discrete copolymer models built over renewal
processes leaves out another possible (and perhaps more natural)
generalization of the simple random walk copolymer model (\ref{eq:SRWmodel}):
namely, the one obtained by replacing the simple random walk
with a more general random walk. A general random walk
crosses the interface without necessarily touching it,
therefore the associated point process is a Markov renewal process
\cite{cfAsm}, because one has to carry along
not only the switching-sign times, but also the height of the walk
at these times (sometimes called the \textit{overshoot}).
This is definitely an interesting
and nontrivial problem that goes in a direction which is complementary
to the one
we have taken. However two remarks are in order:
\begin{enumerate}[(1)]
\item[(1)] Symmetric random walks with i.i.d. increments in $\{-1,0,1\}$
touch the interface when they cross it, hence, they
are covered by our analysis: their weak coupling limit
is the continuum $1/2$-copolymer, because $K(n) \stackrel{n\to\infty
}{\sim}
(\mbox{const.}) n^{-3/2}$ (e.g., \cite{cfBook}, Appendix A.5).

\item[(2)] While one definitely expects an analog of Theorem \ref{th:main}
to hold
for \textit{rather general} random walks with increments in the domain
of attraction of the normal law (with the continuum $1/2$-copolymer as
weak coupling limit), it is less clear what to expect
when the increments of the walk are in the domain of attraction of a
non-Gaussian stable law. In our view, working with generalized
copolymer models has, in any case, a considerable flexibility
with respect to focusing on the random walk set-up.
\end{enumerate}

\section{A closer look at the continuum model}
\label{sec:C0}

In this section, we prove the existence of the continuum free energy
$\widetilde\tf_\ga(\gl,h)$, that is, we prove Theorem \ref{th:existence}.
In Section \ref{sec:modified}, we define a modified partition function,
to which Kingman's super-additive ergodic theorem can be applied,
and then in Section \ref{sec:existence_proof},
we show that this modified partition function yields the
same free energy as the original one.
Before starting with the proof, in Section \ref{sec:redefine} we redefine
the partition function $\widetilde Z_{t,\gb}^{\ga}$
more directly in terms of the $\ga$-stable
regenerative set $\widetilde\tau^\ga$, whose basic properties
are recalled in Section \ref{sec:reg} (cf. also
Appendix \ref{sec:regapp}).
We are going to drop some dependence on $\ga$ for short,
writing, for example, $\widetilde\tf(\gl,h)$.

%%%%%%%%%%%%%%%%%%%%%%%%%%%%%%%%%%%%%%%%%%%%%%%%%%%%%%%%%%%%%%%%%%%%%%

%s2.1 ###
\subsection{Preliminary considerations}
\label{sec:redefine}

As explained in Section \ref{sec:genC0}, the process $\widetilde B^\ga
$ is introduced
just to help visualizing the copolymer, but the underlying relevant process
is $\gD(\widetilde B^\ga) := \ind_{(-\infty,0)}(\widetilde B^\ga)$.
So let us reintroduce $\widetilde Z_{t,\gb}$ more explicitly,
in terms of the random set $\widetilde\tau^\ga$
[cf. (\ref{eq:tautilde})] and of the
signs of the excursions, that are sufficient to determine $\gD
(\widetilde B^\ga)$.

There is no need to pass through
the process $\widetilde B^\ga$ to introduce $\widetilde\tau^\ga$:
we can define it directly as
the \textit{stable regenerative set}
of index $\ga$, that is, the closure of the image
of the stable subordinator of index $\ga$; cf. \cite{cfFFM}.
Some basic properties of regenerative sets are recalled in Section \ref
{sec:reg}
and Appendix \ref{sec:regapp}; we mention in
particular the scale invariance property:
$\widetilde\tau^\ga\sim c \widetilde\tau^\ga$, for every $c>0$.
Since $\widetilde\tau^\ga$
is a closed set, we can write the open
set $(\widetilde\tau^\ga)^\complement= \union_{n\in\N} I_n$
as the disjoint union of countably many (random) open intervals $I_n$,
the connected components (i.e.,
maximal open intervals) of $(\widetilde\tau^\ga)^\complement$.
Although there is no canonical way of numbering these intervals,
any reasonable rule is equivalent for our purpose.
As an example, one first numbers the intervals that start
(i.e., whose left endpoint lies) in $[0,1)$ in decreasing order
of width,
obtaining $\{J^1_n\}_{n\in\N}$; then one does the same with the
intervals that start
in $[1,2)$, getting $\{J^2_n\}_{n\in\N}$; and so on.
Finally, one sets $I_n := J^{a_n}_{b_n}$, where $n \mapsto(a_n, b_n)$
is any
fixed bijection from $\N$ to $\N\times\N$.

Let $\widetilde\xi= \{\widetilde\xi_n\}_{n\in\N}$ be an i.i.d.
sequence of
Bernoulli random variables of parameter $1/2$,
defined on the same probability space as $\widetilde\tau^\ga$ and
independent of~$\widetilde\tau^\ga$, that represent
the signs of the excursions of $\widetilde B^\ga$. We then define the
process $\widetilde\gD^\ga(u):= \sum_n \widetilde\xi_n \ind_{I_n}(u)$,
which takes values in $\{0,1\}$
and is a continuum analogue of the discrete process
$\{\gD_n\}_{n\in\N}$ introduced in Section \ref{sec:discmodel}:
$\widetilde\gD^\ga(u) = 1$ (resp., $0$) means that the continuum
copolymer in $u$ is below
(resp., above or on) the interface.
With this definition, we have the equality in law
%
%e2.1 ###
%
\begin{equation}
\{ \widetilde\gD^\ga(u) \}_{u\ge0} \sim
\{ \gD( \widetilde B^\ga(u) ) \}_{u\ge0} ,
\end{equation}
so that we can use $\widetilde\gD^\ga(\cdot)$ instead
of $\gD(\widetilde B^\ga(\cdot) )$.
More precisely, for $0 \le s \le t < \infty$ we set
%
%e2.2 ###
%
\begin{eqnarray}
\label{eq:Zcont}
\widetilde\myZ_{s,t;\gb} &=& \widetilde\myZ^{\gl, h}_{s,t;\gb}
:= \bE [ \exp( \cH_{s,t;\gb}(\widetilde\gD^\ga
) ) ] ,
\nonumber\\[-8pt]\\[-8pt]
\cH_{s,t;\gb}(\widetilde\gD^\ga) &=&
\cH^{\gl, h}_{s,t;\gb}(\widetilde\gD^\ga)  := - 2 \gl
\int_s^t \widetilde\gD^\ga(u) \bigl(\dd\gb(u) + h \,\dd u \bigr)
,\nonumber
\end{eqnarray}
so that\vspace*{2pt} the partition function $\widetilde Z_{t,\gb}^\ga$ defined
in (\ref{eq:ZWalsh}) coincides with $\widetilde\myZ_{0,t;\gb}$.
For later convenience, we introduce
the finite-volume free energy
%
%e2.3 ###
%
\begin{equation}
\label{eq:tildeFt}
\widetilde\tf_t(\gl,h) := \frac1t \bbE[
\log\widetilde\myZ_{0,t;\gb} ].
\end{equation}
To be precise,\vspace*{2pt} for
$\widetilde\myZ_{s,t;\gb}$ and $\widetilde\tf_t(\gl,h)$ to be
well defined
we need to use a measurable version of $\cH_{s,t;\gb}(\widetilde\gD
^\ga)$
(we build one in Remark \ref{rem:well-posedness} below).

Notice that we have the following additivity property:
%
%e2.4 ###
%
\begin{equation} \label{eq:addit}
\cH_{r,t;\gb}(\widetilde\gD^\ga) =
\cH_{r,s;\gb}(\widetilde\gD^\ga) + \cH_{s,t;\gb}(\widetilde
\gD^\ga) ,
\end{equation}
for every $r < s < t$ and $\bP\otimes\bbP$-a.e. $(\widetilde\gD
^\ga, \gb)$.
Another important observation is that, for any fixed realization
of $\widetilde\Delta^\ga(\cdot)$, the process
$\{ \cH_{s,t;\gb}(\widetilde\gD^\ga)\}_{s,t}$ under $\bbP$ is Gaussian.
\begin{rem}
\label{rem:well-posedness}
Some care is needed for definition (\ref{eq:Zcont})
to make sense. The problem is that $\cH_{s,t;\gb}(\widetilde\gD^\ga)$,
being a stochastic (Wiener) integral, is defined
(for every fixed realization of $\widetilde\gD^\ga$) through an
$L^2$ limit,
hence it is not canonically defined
for every $\gb$, but only $\bbP(\dd\gb)$-a.s.
However, in order to define $\widetilde\myZ_{s,t;\gb}$,
we need $\cH_{s,t;\gb}(\widetilde\gD^\ga)$ to be a measurable
function of $\widetilde\gD^\ga$,
for every (or at least $\bbP$-almost every) fixed $\gb$.
For this reason, we now show that it is possible to define
a \textit{version} of $\cH_{s,t;\gb}(\widetilde\gD^\ga)$, that is, a
measurable function of
$(\gb, \widetilde\gD^\ga,s,t,\gl, h)$.

Let us fix a realization of the process $\{\widetilde\gD^\ga(u)\}_{u
\in[0,\infty)}$.
We build a sequence of approximating functions as follows: for $k\in\N
$ we set
%
%e2.5 ###
%
\begin{equation}
\widetilde\gD^\ga_k(u) :=
\sum_{n \in\N\dvtx |I_n| \ge\frac1k} \widetilde\xi_n \ind
_{I_n}(u) ,
\end{equation}
that is, we only keep the excursion intervals of width at least $\frac1k$.
Note that $\widetilde\gD^\ga_k(u) \to\widetilde\gD^\ga(u)$ as $k
\to\infty$,
for every $u \in\R^+$.
By dominated convergence, we then have $\widetilde\gD^\ga_k \to
\widetilde\gD^\ga$
in $L^2((s,t), \dd x)$,
for all $0 \le s \le t < \infty$, hence by the theory of Wiener integration
it follows that $\lim_{k\to\infty} \cH_{s,t;\gb}(\widetilde\gD
^\ga_k)
= \cH_{s,t;\gb}(\widetilde\gD^\ga)$ in $L^2(\dd\bbP)$.
Note that, for any $k\in\N$, we have
%
%e2.6 ###
%
\begin{equation} \label{eq:HBk}
\cH_{s,t;\gb}(\widetilde\gD^\ga_k) = -2\gl\sum_{n \in\N\dvtx
|I_n| \ge1/k} \widetilde\xi_n
\bigl( \gb_{I_{n} \cap(s,t)} + h |I_{n} \cap(s,t)| \bigr) ,
\end{equation}
where we have set $\gb_{(a,b)} := \gb_{b} - \gb_{a}$ and $\gb
_{\varnothing} := 0$
[note that the right-hand side of (\ref{eq:HBk}) is a finite sum].
This shows that $\cH_{s,t;\gb}(\widetilde\gD^\ga_k)$
is a measurable function of $(\gb, \widetilde\gD^\ga, s, t, \gl,
h)$. Therefore,\vspace*{2pt}
if we prove that $\lim_{k\to\infty} \cH_{s,t;\gb}(\widetilde\gD
^\ga_k) =
\cH_{s,t;\gb}(\widetilde\gD^\ga)$ $\bbP(\dd\gb)$-a.s., we can redefine
$\cH_{s,t;\gb}(\widetilde\gD^\ga) := \liminf_{k\to\infty} \cH
_{s,t;\gb}(\widetilde\gD^\ga_k)$
and get the measurable version we are aiming at.
However,\vspace*{2pt} for any fixed realization of $\widetilde\gD^\ga$, it is
clear from
(\ref{eq:HBk}) that $(\{\cH_{s,t;\gb}(\widetilde\gD^\ga_k)\}
_{k\in\N}, \bbP)$ is a
supermartingale (it is a process with independent
Gaussian increments of negative mean) bounded in $L^2$, hence it
converges $\bbP(\dd\gb)$-a.s.
\end{rem}
%%%%%%%%%%%%%%%%%%%%%%%%%%%%%%%%%%%%%%%%%%%%%%%%%%%%%%%%%%%%%%%%%%%%%%

%s2.2 ###
\subsection{On the $\ga$-stable regenerative set}
\label{sec:reg}

%We recall that $\widetilde\tau$ under $\bP$ is the stable regenerative
%set of index $\ga\in(0,1)$, that is, the closure of the image of the
%$\ga$-stable subordinator. It is a random closed subset of $[0,
%scale-invariant: for every $c >0$, the set $c \widetilde\tau$ has
%the same law of $\widetilde\tau$.

We collect here a few basic formulas on~$\widetilde\tau^\ga$.

For $x \in\R$, we denote by $\bP_x$ the law of the regenerative
set started at $x$, that is,
$\bP_x ( \widetilde\tau^\ga\in\cdot) := \bP( \widetilde\tau
^\ga+ x \in\cdot)$.
Analogously, the process $\{\widetilde\gD^\ga(u)\}_{u\ge x}$ under
$\bP_x$ is distributed
like the process
$\{\widetilde\gD^\ga(u - x)\}_{u\ge x}$ under $\bP=:\bP_0$.
Two variables of basic interest are
%
%e2.7 ###
%
\begin{eqnarray} \label{eq:gtds}
g_t &=& g_t(\widetilde\tau^\ga) := \sup\{ x \in
\widetilde\tau^\ga\cap(-\infty,t] \} ,\nonumber\\[-8pt]\\[-8pt]
d_t &=& d_t(\widetilde\tau^\ga) := \inf\{ x \in
\widetilde\tau^\ga\cap(t, \infty) \} .\nonumber
\end{eqnarray}
%
%or equivalently $g_t = \sup\{ u \le t: B^\ga(u) = 0\}$
%and $d_t = \inf\{ u > t: B^\ga(u) = 0\}$, because
%$\widetilde\tau^\ga= \{u \in[0,\infty): B^\ga(u) = 0\}$.
The joint density of $(g_t,d_t)$ under $\bP_x$ is
%
%e2.8 ###
%
\begin{equation} \label{eq:joint}
\frac{\bP_x ( g_t \in\dd a , d_t \in\dd b )}
{\dd a \, \dd b}
= \frac{\ga \sin(\pi\ga)}{\pi}
\frac{\ind_{(x,t)}(a) \ind_{(t,\infty)}(b)}{(a-x)^{1-\ga}
(b-a)^{1+\ga}} ,
\end{equation}
from which we easily obtain the marginal distribution
of $g_t$: for $y \in[x,t]$
%
%e2.9 ###
%
\begin{equation} \label{eq:G}
G_{x,t}(y) := \bP_x ( g_t \le y ) =
\frac{\sin(\pi\ga)}{\pi} \int_x^{y}
\frac{1}{(a-x)^{1-\ga} (t-a)^\ga} \,\dd a .
\end{equation}
Observing that
$\frac{\dd}{\dd x} (x^\ga/(1-x)^\ga) = \ga (x^{1-\ga
}(1-x)^{1+\ga})^{-1}$,
%for every $x \in(0,1)$,
one obtains also the distribution of $d_t$: for $y \in[t,\infty)$
%
%e2.10 ###
%
\begin{equation} \label{eq:D}
D_{x,t}(y) := \bP_x ( d_t \le y ) =
\frac{\sin(\pi\ga)}{\pi} \int_t^{y}
\frac{(t-x)^\ga}{(b-t)^\ga (b-x)} \,\dd b .
\end{equation}

Let us denote by $\cF_u$ the $\gs$-field generated
by $\widetilde\tau^\ga\cap[0,u]$. The set $\widetilde\tau^\ga$
enjoys the remarkable
\textit{regenerative property}, the continuum analogue of the
well-known renewal property, that can be stated as follows:
for every $\{\cF_u\}_{u \ge0}$-stopping time $\gamma$
such that $\bP(\gamma\in\widetilde\tau^\ga) = 1$,
the portion of $\widetilde\tau^\ga$ after $\gamma$,
that is, the
set $(\widetilde\tau^\ga-\gamma) \cap[0, \infty)$, under $\bP$
is independent of $\cF_{\gamma}$ and distributed like
the original set $\widetilde\tau^\ga$. Analogously,
the translated process $\{\widetilde\gD^\ga(\gamma+ u)\}_{u \ge0}$ is
independent of $\cF_\gamma$ and distributed like the original
process $\widetilde\gD^\ga$.

%%%%%%%%%%%%%%%%%%%%%%%%%%%%%%%%%%%%%%%%%%%%%%%%%%%%%%%%%%%%%%%%%%%%%%

%s2.3 ###
\subsection{A modified partition function}
\label{sec:modified}

In order to apply super-additivity arguments, we
introduce a modification of the partition function.
We extend the Brownian motion $\gb(t)$
to negative times, setting $\gb(t) := \gb'(-t)$ for $t < 0$,
where $\gb'(\cdot)$ is another standard Brownian motion independent
of $\gb$, so that $\gb(t) - \gb(s) \sim\cN(0,t-s)$ for all
$s,t \in\R$ with $s \le t$.

Observe that $\{d_a < b\} = \{\widetilde\tau^\ga\cap(a,b) \ne
\varnothing\}$,
where the random variable $d_t$ has been defined in (\ref{eq:gtds}).
Then for $0 \le s < t$ we set
%
%e2.11 ###
%
\begin{equation}
\widetilde\myZ^{*}_{s,t;\gb} := \inf_{x \in[s-1,s]}
\bE_x [ \exp( \cH_{x,d_{t-1};\gb}(\widetilde\gD^\ga)
) ,
d_{t-1} < t ] .
\end{equation}
In words: we take the smallest free energy
among polymers starting at $x \in[s-1,s]$
and coming back to the interface at some point in $(t-1,t)$.
Notice that the Hamiltonian looks at the polymer only in the
interval $(x, d_{t-1})$. Also notice that for $t < s+1$ the
expression is somewhat degenerate, because for $x > t-1$
we have $d_{t-1} = x$ and therefore $\cH_{x,d_{t-1};\gb} (\widetilde
\gD^\ga)
= \cH_{x,x;\gb} (\widetilde\gD^\ga) = 0$. Therefore,
we may restrict the infimum over $x \in[s-1,\min\{s, t - 1\}]$,
and for clarity we state it explicitly:
%
%e2.12 ###
%
\begin{equation} \label{eq:Zmod}
\widetilde\myZ^{*}_{s,t;\gb}
:= \inf_{x \in[s-1,\min\{s, t-1\}]}
\bE_x [ \exp( \cH_{x,d_{t-1};\gb}(\widetilde\gD^\ga)
) ,
d_{t-1} < t ] .
\end{equation}
Let us stress again that
$\{d_{t-1} < t\} = \{\widetilde\tau^\ga\cap(t-1, t) \ne\varnothing\}$.

It is sometimes more convenient to use $\bE= \bE_0$
instead of $\bE_x$. To this purpose, by a simple change of variables
we have $\cH_{x,a;\gb}(\widetilde\gD^\ga)
= \cH_{0,a-x;\theta_x\gb}(\theta_x \widetilde\gD^\ga)$, where
$\theta_x f(\cdot) := f(x+\cdot)$, as it follows easily from
the definition (\ref{eq:Zcont}) of the Hamiltonian.
Since by definition the process $\theta_x \widetilde\gD^\ga$
under $\bP_x$ is distributed like
$\widetilde\gD^\ga$ under $\bP= \bP_0$, we can write
%
%e2.13 ###
%
\begin{equation} \label{eq:bahlbound}
\bE_x [ \exp( \cH_{x,y;\gb}(\widetilde\gD^\ga) )
] =
\bE[ \exp( \cH_{0,y-x; \theta_x \gb}(\widetilde\gD^\ga
) ) ] .
\end{equation}
Analogously, since the random variable $d_{t-1}$ under $\bP_x$ is distributed
like $x+d_{t-1-x}$ under $\bP$, we can rewrite the term appearing in
(\ref{eq:Zmod}) as
%
%e2.14 ###
%
\begin{eqnarray} \label{eq:bahubound}
&&\bE_x [ \exp( \cH_{x,d_{t-1};\gb}(\widetilde\gD^\ga)
) ,
d_{t-1} < t ] \nonumber\\[-8pt]\\[-8pt]
&&\qquad=
\bE[ \exp( \cH_{0,d_{t-1-x};\theta_x\gb}(\widetilde\gD
^\ga) ) ,
d_{t-1-x} < t - x ] .\nonumber
\end{eqnarray}
These alternative expressions are very useful to get uniform
bounds. In fact, if we set
%
%e2.15 ###
%
\begin{equation}\label{eq:Theta}
\Theta_T (\gb, \widetilde\gD^\ga)
:= {\sup_{-1 \le x \le T , 0 \le y \le T + 1}}
| \cH_{0,y;\theta_x \gb}(\widetilde\gD^\ga) | ,
\end{equation}
from (\ref{eq:Zmod}) and (\ref{eq:bahubound}) we have the following
upper bound:
%
%e2.16 ###
%
\begin{equation} \label{eq:uub}
\sup_{0 \le s < t \le T} \widetilde\myZ^{*}_{s,t;\gb} \le
\bE
[ \exp( \Theta_T (\gb, \widetilde\gD^\ga) )
] .
\end{equation}
In a similar fashion, from relation (\ref{eq:bahlbound}) we obtain
the lower bound
%
%e2.17 ###
%
\begin{equation} \label{eq:llb}
\inf_{-1 \le x \le T , 0 \le y \le T+1} \bE_x [
\exp( \cH_{x,y;\gb}(\widetilde\gD^\ga) ) ] \ge
\bE[ \exp( - \Theta_T(\gb, \widetilde\gD^\ga) )
] .
\end{equation}
We finally state a very useful result which we prove in Appendix \ref
{sec:contbound}:
for every $\eta\in(0,\infty)$ there exists $D(\eta) \in(0,\infty
)$ such that
%
%e2.18 ###
%
\begin{equation} \label{eq:cru}\quad
\bE[ \bbE[ \exp( \eta
\Theta_T (\gb, \widetilde\gD^\ga) ) ] ]
\le D(\eta) e^{D(\eta) T} < \infty\qquad
\mbox{for every } T > 0 .
\end{equation}

%%%%%%%%%%%%%%%%%%%%%%%%%%%%%%%%%%%%%%%%%%%%%%%%%%%%%%%%%%%%%%%%%%%%%%
%s2.4 ###
\subsection{\texorpdfstring{Proof of Theorem
\protect\ref{th:existence}}{Proof of Theorem 1.4}}
\label{sec:existence_proof}

We start by\vspace*{1pt} proving the existence of the limit in (\ref{eq:existence}) when
the partition function $\widetilde Z_{t,\gb}^\ga= \widetilde\myZ
_{0,t;\gb}$
is replaced by $\widetilde\myZ^*_{0,t;\gb}$.
%Then we are going to show that the same result can be transferred to
%the original partition function.
%
\begin{proposition} \label{th:Kin}
For all $\gl, h \ge0$,
the following limit exists $\bbP(\dd\gb)$-a.s. and in $L^1(\dd\bbP)$:
%
%e2.19 ###
%
\begin{equation} \label{eq:Kin}
\lim_{t \to\infty} \frac1t \log \widetilde\myZ
^{*}_{0,t;\gb}
=: \widehat\tf(\gl,h) ,
\end{equation}
where $\widehat\tf(\gl,h)$ is finite and nonrandom.
\end{proposition}
\begin{pf}
We are going to check that, for all fixed $\gl, h \ge0$,
the process $\{\log\widetilde\myZ_{s,t;\gb}^{*}\}_{0 \le s < t <
\infty}$ under $\bbP$
satisfies the four hypotheses of
\textit{Kingman's super-additive ergodic theorem}; cf. \cite{cfKin}.
This entails the existence of the limit in the left-hand side of (\ref{eq:Kin}),
both $\bbP$-a.s. and in $L^1(\dd\bbP)$, as well as the fact that the
limit is a function of $\gb$ which is invariant under time translation
$\gb(\cdot) \mapsto\theta_t \gb(\cdot) := \gb(t+\cdot)$, for
every $t \ge0$.
Therefore, the limit must be measurable w.r.t. the tail $\gs$-field
of $\gb(\cdot)$, hence nonrandom by Kolmogorov 0--1 law for Brownian motion.

The first of Kingman's conditions is that for every $k \in\N$,
any choice of $\{(s_j, t_j)\}_{k\in\N}$,
with $0 \le s_j < t_j$, and for every $a > 0$ we have
%
%e2.20 ###
%
\begin{equation}
( \widetilde\myZ^{*}_{s_1, t_1;\gb} , \ldots,
\widetilde\myZ^{*}_{s_k, t_k;\gb} ) \stackrel{d}{=}
( \widetilde\myZ^{*}_{s_1+a, t_1+a;\gb} , \ldots,
\widetilde\myZ^{*}_{s_k+a, t_k+a;\gb} ) .
\end{equation}
However this is trivially true, because
$\widetilde\myZ^{*}_{s+a, t+a;\gb} =
\widetilde\myZ^{*}_{s, t;\theta_a\gb}$, as it follows from (\ref{eq:Zmod}),
recalling the definition of the Hamiltonian in (\ref{eq:Zcont}).

The second condition is the super-additivity property:
for all $0 \le r < s < t$
%
%e2.21 ###
%
\begin{equation} \label{eq:super-ad}
\widetilde\myZ^{*}_{r, t; \gb} \ge \widetilde\myZ^{*}_{r, s;
\gb} \cdot
\widetilde\myZ^{*}_{s, t; \gb} .
\end{equation}
To this purpose, for any fixed $x \in[r-1, r]$ the inclusion bound yields
%
%e2.22 ###
%
\begin{eqnarray} \label{eq:mmm}\quad
&& \bE_x \bigl( \exp( \cH_{x,d_{t-1};\gb} ) ,
d_{t-1} < t \bigr)\nonumber\\[-8pt]\\[-8pt]
&&\qquad \ge
\bE_x \bigl( \exp( \cH_{x,d_{s-1};\gb} )
\exp( \cH_{d_{s-1},d_{t-1};\gb} ) ,
d_{s-1} < s , d_{t-1} < t \bigr) ,\nonumber
\end{eqnarray}
where we have used the additivity of the Hamiltonian, see (\ref{eq:addit}).
We integrate over the possible values of $d_{s-1}$ and, using
the regenerative property, we rewrite the right-hand side of (\ref{eq:mmm})
as follows:
%
%e2.23 ###
%
\begin{eqnarray} \label{eq:mmm2}\quad
&& \int_{y \in(s-1, s)}
\bE_x \bigl( \exp( \cH_{x,y;\gb} ) ,
d_{s-1} \in\dd y \bigr)
\bE_y \bigl( \exp( \cH_{y,d_{t-1};\gb} ) ,
d_{t-1} < t \bigr)\nonumber\\[-8pt]\\[-8pt]
&&\qquad \ge \bE_x \bigl( \exp( \cH_{x,d_{s-1};\gb} ),
d_{s-1} < s \bigr)
\cdot\widetilde\myZ^{*}_{s,t;\gb} ,\nonumber
\end{eqnarray}
where the inequality is just a consequence of taking the infimum over
$y \in[s-1, s]$
and recalling the definition (\ref{eq:Zmod}) of $\widetilde\myZ
^{*}_{s,t;\gb}$.
Putting together the relation (\ref{eq:mmm}) and
(\ref{eq:mmm2}) and taking the infimum over $x \in[r-1, r]$, we have proven
(\ref{eq:super-ad}).

The third condition to check is
%
%e2.24 ###
%
\begin{equation} \label{eq:easybound}
\sup_{t > 0} \frac1t \bbE( \log\widetilde\myZ
^{*}_{0,t;\gb} )
< \infty.
\end{equation}
Recalling (\ref{eq:Zmod}) and applying Jensen's inequality
and Fubini's theorem, we can write
%
%e2.25 ###
%
\begin{equation} \label{eq:JenFub}
\bbE( \log\widetilde\myZ^{*}_{0,t;\gb} ) \le
\log \bE\bigl( \bbE[
\exp( \cH_{0,d_{t-1};\gb}(\widetilde\gD^\ga) ) ],
d_{t-1} < t \bigr) .
\end{equation}
Since the Hamiltonian is a stochastic integral [cf. (\ref{eq:Zcont})]
for fixed $a < b$ and $\widetilde\gD^\ga$ we have
$\cH_{a,b;\gb}(\widetilde\gD^\ga) \sim\cN(\mu, \gs^2)$, where
$\mu= -2 \gl h \int_a^b \widetilde\Delta^\ga(u) \,\dd u$ and
$\gs^2 = 4 \gl^2 \int_a^b |\widetilde\Delta^\ga(u)|^2 \,\dd u$.
In particular, $|\mu| \le2 \gl h (b-a)$ and
$\gs^2 \le4 \gl^2 (b-a)$, hence, on the
event $\{d_{t-1} < t\}$, we have
$\bbE[ \exp( \cH_{0,d_{t-1};\gb}(\widetilde\gD^\ga)
) ]
\le\exp(2\gl h t + 2 \gl^2 t)$, and (\ref{eq:easybound}) follows.

Finally, the fourth and last condition is that for some
(hence every) $T > 0$,
%
%e2.26 ###
%
\begin{equation} \label{eq:sup}
\bbE\Bigl(
{\sup_{0 \le s < t \le T}} | {\log\widetilde\myZ^{*}_{s,t;\gb
} }|
\Bigr) < \infty.
\end{equation}
We need both a lower and an upper bound on $\widetilde\myZ
^{*}_{s,t;\gb}$.
For the upper bound, directly from (\ref{eq:uub}) we have
%
%e2.27 ###
%
\begin{equation} \label{eq:uubb}
\sup_{0 \le s< t \le T} \log\widetilde\myZ^{*}_{s,t;\gb} \le
\log\bE( \exp( \Theta_T(\gb, \widetilde\gD^\ga)
) ) .
\end{equation}
The lower bound is slightly more involved.
The additivity of the Hamiltonian yields
$\cH_{x,d_{t-1};\gb}(\widetilde\gD^\ga) =
\cH_{x,t-1;\gb}(\widetilde\gD^\ga) + \cH_{t-1,d_{t-1};\gb
}(\widetilde\gD^\ga)$.
Observing that $\widetilde\Delta^\ga(s)$ is
constant for $s \in(t-1, d_{t-1}(\widetilde\tau^\ga))$, from the definition
(\ref{eq:Zcont}) of the Hamiltonian, we can write
%
%e2.28 ###
%
\begin{eqnarray} \label{eq:Cbeta}
\cH_{t-1, d_{t-1};\gb}(\widetilde\gD^\ga) &\ge&
- 2\gl |\gb_{d_{t-1}} - \gb_{t-1}| - 2\gl h
\bigl(d_{t-1} - (t-1)\bigr) \nonumber\\[-8pt]\\[-8pt]
&\ge& - {2\gl \sup_{0 \le s < t \le T}} |\gb_t - \gb_s|
- 2 \gl h T =: -C_T(\gb) .\nonumber
\end{eqnarray}
Recalling (\ref{eq:Zmod}), we can therefore
bound $\widetilde\myZ^{*}_{s,t;\gb}$ from below by
%
%e2.29 ###
%
\begin{eqnarray}\qquad
\widetilde\myZ^{*}_{s,t;\gb} &\ge& e^{-C_T(\gb)}\Bigl(
\inf_{x \in[s-1, \min\{s, t-1\}]} \bE_x
\bigl( \exp(\cH_{x,t-1;\gb}(\widetilde\gD^\ga) )
|
d_{t-1} < t \bigr) \Bigr)\nonumber\\[-8pt]\\[-8pt]
&&{}\times \bP_x ( d_{t-1} < t ) .\nonumber
\end{eqnarray}
From (\ref{eq:D}) it follows easily that, for fixed $T$,
%
%e2.30 ###
%
\begin{equation}
\inf_{0 \le s < t \le T}
\inf_{x \in[s-1, \min\{s, t-1\}]}
\bP_x ( d_{t-1} < t ) > 0 .
\end{equation}
Furthermore, we now show that we can replace the law $\bP_x( \cdot
| d_{t-1} < t )$
with $\bP_x( \cdot)$ by paying a positive constant.
In fact, the laws of the set $\widetilde\tau^\ga\cap[x,t-1]$
under these two probability measures are mutually absolutely continuous.
The Radon--Nikodym derivative, which depends only on $g_{t-1}$, is computed
with the help of (\ref{eq:joint}), (\ref{eq:G}), (\ref{eq:D}) and equals
%
%e2.31 ###
%
\begin{eqnarray} \label{eq:R-N}
&&\frac{\dd\bP_x( \cdot| d_{t-1} < t )}{\dd\bP_x( \cdot)}
(\widetilde\tau^\ga\cap[x,t-1]) \nonumber\\
&&\qquad =
\frac{\bP_x(g_{t-1} \in\dd y, d_{t-1} < t)}
{\bP_{x}(g_{t-1} \in\dd y) \bP_x(d_{t-1} < t)} \bigg|_{y =
g_{t-1}}\\
&&\qquad = \biggl( 1 -
\frac{(t-1-g_{t-1})^\ga}{(t-g_{t-1})^\ga} \biggr) \cdot
\frac{1}{D_{x,t-1}(t)} .\nonumber
\end{eqnarray}
Using (\ref{eq:D}), it is straightforward to check that,
for every fixed $T$, the infimum of this expression
over $0 \le s < t \le T$ and $x \in[s-1, \min\{s, t-1\}]$ is strictly
positive. Therefore, uniformly in the range of parameters, we have
%
%e2.32 ###
%
\begin{eqnarray}\quad
\widetilde\myZ^{*}_{s,t;\gb} &\ge& (\mbox{const.}) e^{-C_T(\gb)}
\inf_{x \in[s-1, \min\{s, t-1\}]} \bE_x
( \exp( \cH_{x,t-1;\gb}(\widetilde\gD^\ga) )
)\nonumber\\[-8pt]\\[-8pt]
&\ge& (\mbox{const.}) e^{-C_T(\gb)} \bE(
\exp( -\Theta_T(\gb, \widetilde\gD^\ga) )
) ,\nonumber
\end{eqnarray}
where we have applied (\ref{eq:llb}). By Jensen's inequality, we then obtain
%
%e2.33 ###
%
\begin{equation} \label{eq:llbb}
\inf_{0 \le s < t \le T} \log\widetilde\myZ^{*}_{s,t;\gb}
\ge - \bE( \Theta_T(\gb, \widetilde\gD^\ga) )
- C_T(\gb) + (\mbox{const.}') .
\end{equation}
Putting together (\ref{eq:uubb}) and (\ref{eq:llbb}), we then get
%
%e2.34 ###
%
\begin{eqnarray}
{\sup_{0 \le s < t \le T}} |{ \log\widetilde\myZ^{*}_{s,t;\gb}}|
&\le& \log\bE( \exp( \Theta_T(\gb, \widetilde\gD
^\ga) ) )
+ \bE( \Theta_T(\gb, \widetilde\gD^\ga) )\nonumber\\[-8pt]\\[-8pt]
&&{} + C_T(\gb) + (\mbox{const.}) .\nonumber
\end{eqnarray}
It is clear from (\ref{eq:Cbeta}) that $\bbE(C_T(\gb)) < \infty$,
for every $T > 0$.
Moreover,\break by Jensen's inequality and (\ref{eq:cru}) we have
$\bbE\log\bE[ \exp(\Theta_T(\gb, \widetilde\gD^\ga)) ] \le\break
\log\bE[ \bbE[ \exp(\Theta_T(\gb, \widetilde\gD^\ga))] ] <
\infty$, so that
$\bE[ \bbE[\Theta_T(\gb, \widetilde\gD^\ga)] ] < \infty$.
Therefore, (\ref{eq:sup}) is proven.
\end{pf}

We finally show that Proposition \ref{th:Kin} still holds
if we replace the modified partition function $\widetilde\myZ
^{*}_{0,t;\gb}$
with the original partition function
$\widetilde\myZ_{0,t;\gb}$; in particular, the free energy
$\widetilde\tf(\gl, h)$
is well defined and coincides with $\widehat\tf(\gl,h)$.
We first need a technical lemma.
\begin{lemma} \label{th:barabao}
For every fixed $h \ge0$,
the function $\widehat\tf(\gl,h)$ appearing in Proposition \ref{th:Kin}
is a nondecreasing and continuous function of $\gl$.
\end{lemma}
\begin{pf}
Note that sending $\gl\to c\gl$ is the same as multiplying
the Hamiltonian by $c$. By Jensen's inequality, for every $\gep> 0$ we have
%
%e2.35 ###
%
\begin{equation}\hspace*{28pt}
\bE_x \bigl( \exp( \cH_{x,d_{t-1};\gb} )
\ind_{\{d_{t-1}<t\}} \bigr)^{1+\gep} \le
\bE_x \bigl( \exp\bigl( (1+\gep) \cH_{x,d_{t-1};\gb} \bigr)
\ind_{\{d_{t-1}<t\}} \bigr) ,
\end{equation}
hence, taking the infimum over $x \in[-1,0]$, then $\frac1t \bbE\log
(\cdot)$
and letting $t \to\infty$, we obtain
$\widehat\tf((1+\gep)\gl, h) \ge(1+\gep) \widehat\tf(\gl, h)$.
In particular,
$\gl\mapsto\widehat\tf(\gl,h)$ is nondecreasing for fixed $h$.

To prove the continuity, we use H\"older's inequality
with $p=\frac1{1-\gep}$ and $q=\frac{1}{\gep}$, getting
%
%e2.36 ###
%
\begin{eqnarray}\label{eq:bilibao}
&&
\bE_x \bigl( e^{(1+\gep)\cH_{x,d_{t-1};\gb}}
\ind_{\{d_{t-1}<t\}} \bigr) \nonumber\\
&&\qquad=
\bE_x \bigl( e^{(1-\gep)\cH_{x,d_{t-1};\gb}}
e^{2 \gep\cH_{x,d_{t-1};\gb}} \ind_{\{d_{t-1}<t\}} \bigr)
\\
&&\qquad\le \bE_x \bigl( e^{\cH_{x,d_{t-1};\gb}}
\ind_{\{d_{t-1}<t\}} \bigr)^{1-\gep} \bE_x \bigl( e^{2 \cH
_{x,d_{t-1};\gb}}
\ind_{\{d_{t-1}<t\}} \bigr)^{\gep} .\nonumber
\end{eqnarray}
Now observe that by (\ref{eq:bahubound}) and (\ref{eq:Theta}) we can write
%
%e2.37 ###
%
\begin{equation}
\bE_x \bigl( e^{2 \cH_{x,d_{t-1};\gb}}
\ind_{\{d_{t-1}<t\}} \bigr)^{\gep} \le
\bE\bigl( e^{2 \Theta_{t+1}(\gb, \widetilde\gD^\ga)} \bigr)^{\gep
} .
\end{equation}
Taking $\frac1t \bbE \inf_{x \in[-1,0]} \log(\cdot)$ in
(\ref{eq:bilibao}),
applying Jensen's inequality to the last term, using (\ref{eq:cru})
and letting $t \to\infty$ then yields
%
%e2.38 ###
%
\begin{eqnarray}
\widehat\tf\bigl((1+\gep)\gl, h\bigr) \le (1-\gep) \widehat\tf(\gl
, h)
+ \gep D(2) \nonumber\\[-8pt]\\[-8pt]
\eqntext{\mbox{for every } \gl, h \ge0
\mbox{ and every } \gep> 0 .}
\end{eqnarray}
Since $\gl\mapsto\widehat\tf(\gl,h)$ is nondecreasing,
this shows that $\gl\mapsto\widehat\tf(\gl,h)$ is continuous.
%\rightqed
\end{pf}

We now pass from $\widetilde\myZ^{*}_{0,t;\gb}$ to the original partition
function $\widetilde\myZ_{0,t;\gb}$ in three steps: first, we remove the
infimum over $x \in[-1,0]$, then we replace $\cH_{0,d_{t-1};\gb}$
with $\cH_{0,t-1;\gb}$ and finally we remove the event $\{d_{t-1}<t\}$.
From now till the end of the proof, we assume $t\ge1$.

%%%%%%%%%%%%%%%%%%%%%%%%%%%%%%%%%

%s2.4.1 ###
\subsubsection*{Step 1}

It follows from the regenerative property of $\widetilde\tau^\ga$
that the laws of the random set $\widetilde\tau^\ga\cap[1,\infty)$
under the probabilities $\bP= \bP_0$ and $\bP_x$, with $x \in[-1,0]$,
are mutually absolutely continuous, with Radon--Nikodym derivative
depending only on $d_1$, given by
%
%e2.39 ###
%
\begin{equation}\quad
\frac{\dd\bP(\widetilde\tau^\ga\cap[1,\infty) \in\cdot)}
{\dd\bP_x(\widetilde\tau^\ga\cap[1,\infty) \in\cdot)} =
\frac{\bP(d_1 \in\dd z)}{\bP_x(d_1 \in\dd z)} \bigg|_{z
= d_1} =
\frac{1}{(1-x)^\ga} \frac{d_1}{d_1 - x} .
\end{equation}
It is clear that, uniformly on $x \in[-1,0]$,
this expression is bounded from above by
some constant $0 < C < \infty$. Therefore, for every $\gep> 0$,
by the H\"older inequality with $p = \frac{1+\gep}{\gep}$
and $q=1+\gep$ we can write
%
%e2.40 ###
%
\begin{eqnarray}
&&
\bE\bigl( e^{\cH_{0,d_{t-1};\gb}} \ind_{\{d_{t-1} < t\}} \bigr)
\nonumber\\
&&\qquad= \bE\bigl( e^{\cH_{0,1;\gb} + \cH_{1,d_{t-1};\gb}}
\ind_{\{d_{t-1} < t\}} \bigr) \nonumber\\
&&\qquad \le \bE\bigl( e^{({1+\gep})/{\gep} \cH_{0,1;\gb}}
\bigr)^{{\gep}/({1+\gep})}
\bE\bigl( e^{(1+\gep)\cH_{1,d_{t-1};\gb}}
\ind_{\{d_{t-1} < t\}} \bigr)^{1/(1+\gep)} \\
&&\qquad \le \bE\bigl( e^{({1+\gep})/{\gep} \cH_{0,1;\gb}}
\bigr)^{{\gep}/({1+\gep})}
C^{1/(1+\gep)}\nonumber\\
&&\qquad\quad{}\times
\inf_{x \in[-1,0]} \bE_x \bigl( e^{(1+\gep)\cH_{1,d_{t-1};\gb}}
\ind_{\{d_{t-1} < t\}} \bigr)^{1/(1+\gep)} .\nonumber
\end{eqnarray}
Analogously, again by the H\"older inequality, we have
%
%e2.41 ###
%
\begin{eqnarray}\qquad
&&
\bE_x \bigl( e^{(1+\gep)\cH_{1,d_{t-1};\gb}} \ind_{\{d_{t-1}
< t\}} \bigr) \nonumber\\
&&\qquad=
\bE_x \bigl( e^{(1+\gep)(\cH_{x,d_{t-1};\gb} - \cH_{x,1;\gb})}
\ind_{\{d_{t-1} < t\}} \bigr) \\
&&\qquad \le \bE_x \bigl( e^{-{(1+\gep)^2}/{\gep} \cH_{x,1;\gb}}
\bigr)^{{\gep}/{(1+\gep)}}
\bE_x \bigl( e^{(1+\gep)^2 \cH_{x,d_{t-1};\gb} }
\ind_{\{d_{t-1} < t\}} \bigr)^{1/(1+\gep)} .\nonumber
\end{eqnarray}

However, $\bE_x ( e^{-{(1+\gep)^2}/{\gep} \cH_{x,1;\gb}}
)
\le\bE( e^{{(1+\gep)^2}/{\gep} \Theta_2(\gb, \widetilde
\gD^\ga)} )$,
by (\ref{eq:bahubound}) and (\ref{eq:Theta}).
Putting together these relations, Proposition \ref{th:Kin} and (\ref
{eq:cru}), we get
$\bbP(\dd\gb)$-a.s.
%
%e2.42 ###
%
\begin{eqnarray}
\label{eq:tomltln}
&&\limsup_{t\to\infty} \frac1 t \log
\bE\bigl( e^{\cH_{0,d_{t-1};\gb}} \ind_{\{d_{t-1} < t\}} \bigr)
\nonumber\\
&&\qquad\le \frac{1}{(1+\gep)^2}
\limsup_{t\to\infty} \frac1 t \log
\inf_{x \in[-1,0]} \bE_x \bigl( e^{(1+\gep)^2 \cH
_{x,d_{t-1};\gb}}
\ind_{\{d_{t-1} < t\}} \bigr) \\
&&\qquad=
\frac{\widehat\tf((1+\gep)^2 \gl, h)}{(1+\gep)^2} ,\nonumber
\end{eqnarray}
and since $\gep> 0$ is arbitrary, by Lemma \ref{th:barabao}
the left-hand side in (\ref{eq:tomltln}) does not exceed
$\widehat\tf(\gl, h)$.
By the definition (\ref{eq:Zmod}) of $\widetilde\myZ^{*}_{0,t;\gb
}$, we have immediately
an analogous lower bound for the $\liminf$, hence we have proven
that $\bbP(\dd\gb)$-a.s.
%
%e2.43 ###
%
\begin{equation} \label{eq:quasife}
\lim_{t\to\infty} \frac1 t \log
\bE\bigl( e^{\cH_{0,d_{t-1};\gb}} \ind_{\{d_{t-1} < t\}} \bigr)
= \widehat\tf(\gl, h) .
\end{equation}
Furthermore, the convergence holds also in $L^1(\bbP)$, because the sequence
in the left-hand side is uniformly integrable, as it follows from the
bounds we
have obtained.

%%%%%%%%%%%%%%%%%%%%%%%%%%%%%%%%%

%s2.4.2 ###
\subsubsection*{Step 2}

With analogous arguments, we now show that we can replace $\cH
_{0,d_{t-1};\gb}$
with $\cH_{0,t-1;\gb}$ in (\ref{eq:quasife}), that is,
the following limit holds, $\bbP(\dd\gb)$-a.s. and in $L^1(\dd\bbP)$:
%
%e2.44 ###
%
\begin{equation} \label{eq:quasife1}
\lim_{t\to\infty} \frac1 t \log
\bE\bigl( e^{\cH_{0,t-1;\gb}} \ind_{\{d_{t-1} < t\}} \bigr)
= \widehat\tf(\gl, h) .
\end{equation}
Since $\cH_{0,d_{t-1};\gb} = \cH_{0,t-1;\gb} + \cH
_{t-1,d_{t-1};\gb}$,
for every $\gep> 0$, we can write
%
%e2.45 ###
%
\begin{eqnarray}\label{eq:samebound1}
&&
\bE\bigl( e^{(1-\gep)\cH_{0,d_{t-1};\gb}} \ind_{\{d_{t-1} < t\}
} \bigr)
\nonumber\\[-8pt]\\[-8pt]
&&\qquad\le \bE\bigl( e^{({1-\gep})/{\gep} \cH_{t-1,d_{t-1};\gb}}
\ind_{\{d_{t-1} < t\}} \bigr)^{\gep}
\bE\bigl( e^{\cH_{0,t-1;\gb}}
\ind_{\{d_{t-1} < t\}} \bigr)^{1-\gep} ,\nonumber
\end{eqnarray}
and analogously
%
%e2.46 ###
%
\begin{eqnarray}\label{eq:samebound2}
&& \bE\bigl( e^{\cH_{0,t-1;\gb}} \ind_{\{d_{t-1} < t\}} \bigr)
\nonumber\\
&&\qquad \le \bE\bigl( e^{- ({1+\gep})/{\gep} \cH
_{t-1,d_{t-1};\gb}}
\ind_{\{d_{t-1} < t\}} \bigr)^{{\gep}/({1+\gep})}\\
&&\qquad\quad{}\times
\bE\bigl( e^{(1+\gep)\cH_{0,d_{t-1};\gb}}
\ind_{\{d_{t-1} < t\}} \bigr)^{{1}/({1+\gep})} .\nonumber
\end{eqnarray}
Now notice that, by definition (\ref{eq:Zcont}),
since $\widetilde\tau^\ga\cap(t-1, d_{t-1}) = \varnothing$, we can write
%
%e2.47 ###
%
\begin{equation} \label{eq:useboh}
| \cH_{t-1,d_{t-1};\gb} |
\le 2\gl \bigl( |\gb_{d_{t-1}} - \gb_{t-1}|
+ h \bigl(d_{t-1} - (t-1)\bigr) \bigr) ,
\end{equation}
from which it follows easily that $\bbP(\dd\gb)$-a.s. and in
$L^1(\dd\bbP)$
%
%e2.48 ###
%
\begin{equation}
\lim_{t\to\infty}
\frac1t \log \bE\bigl( e^{ \gamma |\cH_{t-1,d_{t-1};\gb}|}
\ind_{\{d_{t-1} < t\}} \bigr) = 0 \qquad \forall\gamma\ge
0 .
\end{equation}
From (\ref{eq:samebound1}), (\ref{eq:samebound2}) and (\ref{eq:quasife})
we then have $\bbP(\dd\gb)$-a.s.
\begin{eqnarray*}
\frac{\widehat\tf((1-\gep)\gl, h)}{1-\gep} &\le&
\liminf_{t\to\infty} \frac1 t \log
\bE\bigl( e^{\cH_{0,t-1;\gb}} \ind_{\{d_{t-1} < t\}} \bigr) \\
&\le& \limsup_{t\to\infty} \frac1 t \log
\bE\bigl( e^{\cH_{0,t-1;\gb}} \ind_{\{d_{t-1} < t\}} \bigr) \\
&\le& \frac{\widehat\tf((1+\gep)\gl, h)}{1+\gep} .
\end{eqnarray*}
Letting $\gep\to0$ and using Lemma \ref{th:barabao}, we see that
(\ref{eq:quasife1}) holds $\bbP(\dd\gb)$-a.s. and also in $L^1(\dd
\bbP)$,
thanks to the bounds (\ref{eq:samebound1}), (\ref{eq:samebound2}) and
(\ref{eq:useboh}) that ensure the uniform integrability.

%%%%%%%%%%%%%%%%%%%%%%%%%%%%%%%%%

%s2.4.3 ###
\subsubsection*{Step 3}

We finally show that we can remove the indicator function $\ind_{\{
d_{t-1} < t\}}$
from (\ref{eq:quasife1}).
We have already observed that the laws of $\widetilde\tau^\ga\cap
[0,t-1]$ under
the two probabilities $\bP( \cdot| d_{t-1} < t)$ and $\bP$ are mutually
absolutely continuous: the corresponding Radon--Nikodym derivative $f_t
= f_t(g_{t-1})$
is given by (\ref{eq:R-N}), from which we extract the bound
%
%e2.49 ###
%
\begin{equation}
f_t(g_{t-1}) \ge 1 - \frac{(t-g_{t-1}-1)^\ga}{(t-g_{t-1})^\ga}
\ge 1 - \frac{(t-1)^\ga}{t^\ga} \ge \frac{\ga}{t} ,
\end{equation}
where the last inequality holds for large $t$.
Therefore, for large $t$,
%
%e2.50 ###
%
\begin{eqnarray}
\bE\bigl( e^{\cH_{0,t-1;\gb}} \ind_{\{d_{t-1} < t\}} \bigr)
&=& \bE( e^{\cH_{0,t-1;\gb}} | d_{t-1} < t )
\bP(d_{t-1} < t) \nonumber\\[-8pt]\\[-8pt]
&\ge&
\frac{\ga}{t} \bE( e^{\cH_{0,t-1;\gb}} )
\bP(d_{t-1} < t) ,\nonumber
\end{eqnarray}
and note that $\bP(d_{t-1} < t) = G_{0,t-1}(t) \sim(\mbox{const.})/t^{1-\ga}$
as $t\to\infty$, by (\ref{eq:D}). Therefore,
%
%e2.51 ###
%
\begin{equation}\hspace*{32pt}
\bE\bigl( e^{\cH_{0,t-1;\gb}} \ind_{\{d_{t-1} < t\}} \bigr)
\le
\bE( e^{\cH_{0,t-1;\gb}} ) \le
(\mbox{const.}) t^{2-\ga}
\bE\bigl( e^{\cH_{0,t-1;\gb}} \ind_{\{d_{t-1} < t\}} \bigr) ,
\end{equation}
for large $t$, hence by (\ref{eq:quasife1}) it follows that,
$\bbP(\dd\gb)$-a.s. and in $L^1(\dd\bbP)$, we have
%
%e2.52 ###
%
\begin{equation}
\lim_{t\to\infty} \frac1 t \log
\bE( e^{\cH_{0,t-1;\gb}} )
= \widehat\tf(\gl, h) .
\end{equation}
Replacing $\frac1 t$ with $\frac{1}{t-1}$ in the left-hand side
shows that the free energy $\widetilde\tf(\gl,h)$,
defined as the limit in (\ref{eq:existence}), does exist and coincides
with $\widehat\tf(\gl,h)$ (we recall that $\widetilde Z_{t,\gb}^\ga
= \widetilde\myZ_{0,t;\gb}$).

To complete the proof of Theorem \ref{th:existence}, it only
remains to show that the free energy $\widetilde\tf(\gl, h)$
is nonnegative and continuous.
By restricting, for $t>1$, the expectation that
defines $\widetilde Z_{0,t;\gb}$ to the event $E_t:= \{ d_1 > t,
\widetilde\gD^\ga(\frac{t+1}{2}) = 0 \} = \{ d_1 > t,
\widetilde B^\ga(\frac{t+1}{2})>0\}$ and by using
Jensen inequality, we have
%
%e2.53 ###
%
\begin{eqnarray}
\label{eq:quspr}
\frac1t \bbE\log\widetilde\myZ_{0,t;\gb} &\ge& \frac1t
\bE[ \bbE[\cH_{0,1;\gb} (\widetilde\gD^\ga)]
\vert E_t ] + \frac1t \log\bP(E_t) \nonumber\\[-8pt]\\[-8pt]
&\ge&
-\frac{2\gl h}{t} + \frac1t \log\bP(E_t).\nonumber
\end{eqnarray}
By (\ref{eq:D}), we have $\bP(E_t) \stackrel{t \to\infty}\sim
(\mbox{const.}) t^{-\ga}$
so that the right-most side in (\ref{eq:quspr}) vanishes as $t \to
\infty$
and therefore $\widetilde\tf(\gl, h)\ge0$.

For the continuity, it is convenient to use a different parametrization.
For $t > 0$ and $a,b \in\R$, we set
%
%e2.54 ###
%
\begin{equation} \label{eq:tgt}
\tg_t(a,b) := \frac1t \bbE\biggl[ \log\bE
\biggl[ \exp\biggl( - 2 \int_0^t \widetilde\gD^\ga(u)
\bigl(a \,\dd\gb(u) + b \,\dd u \bigr) \biggr)
\biggr] \biggr] .
\end{equation}
Since the argument of the exponential is a bilinear function of $(a,b)$,
it is easily checked, using H\"older's inequality, that
for every fixed $t>0$ the function $(a,b) \mapsto\tg_t(a,b)$
is convex on $\R^2$. By a straightforward adaptation of the results
proven in this section, the limit
%
%e2.55 ###
%
\begin{equation} \label{eq:tg}
\tg(a, b) := \lim_{t\to\infty} \tg_t (a,b)
\end{equation}
exists and is finite, for all $a,b \in\R$.
For instance, for $a > 0$ and $b \ge0$,
by (\ref{eq:tildeFt}) and (\ref{eq:Zcont}) we have
$\tg_t(a, b) = \widetilde\tf_t(a, b/a)$, therefore
the limit in (\ref{eq:tg}) exists and equals $\widetilde\tf(a, b/a)$;
the restriction to $a > 0$ and $b \ge0$ is however not necessary
for the existence of such a limit.

Being the pointwise limit of convex functions, $\tg(a,b)$ is convex
too on $\R^2$,
hence continuous (because finite).
Therefore, $\widetilde\tf(\gl, h) = \tg(\gl, \gl h)$ is continuous
too on $[0,\infty) \times[0,\infty)$.

%%%%%%%%%%%%%%%%%%%%%%%%%%%%%%%%%%%%%%%%%%%%%%%%%%%%%%%%%%%%%%%%%%%%%%
%%%%%%%%%%%%%%%%%%%%%%%%%%%%%%%%%%%%%%%%%%%%%%%%%%%%%%%%%%%%%%%%%%%%%%
%%%%%%%%%%%%%%%%%%%%%%%%%%%%%%%%%%%%%%%%%%%%%%%%%%%%%%%%%%%%%%%%%%%%%%
%s3 ###
\section{The proof of the main result}
\label{sec:main_result_proof}

We fix an arbitrary value of $\ga\in(0,1)$
and an arbitrary discrete $\ga$-copolymer model
(and we omit
$\ga$ in most of the notations of this section).
We aim at proving an analogue of Theorem 6 in \cite{cfBdH}.
More precisely, we want to show the following theorem.
\begin{theorem}
\label{th:main-tech}
For every choice of
$\gl>0$ and $h>0$, and for every choice of $\rho\in(0,1)$
there exists $a_0>0$ such that for every $a \in(0, a_0]$, we have
%
%e3.1 ###
%
\begin{equation}
\label{eq:main}
\widetilde\tf\biggl( \frac{\gl}{1+\rho}, \frac{h}{1-\rho} \biggr)
\le \frac{1}{a^2} \tf(a \gl, a h) \le
\widetilde\tf\bigl( (1+\rho) \gl, (1-\rho) h \bigr) .
\end{equation}
\end{theorem}

Theorem \ref{th:main-tech} implies Theorem \ref{th:main}.
In fact, notice that it directly yields (\ref{eq:mainfe})
when both $\gl$ and $h$ are positive [by continuity of
$\widetilde\tf(\cdot, \cdot)$].
If $\gl=0$, there is nothing to prove, because
$\tf(0,h) = \widetilde\tf(0,h) = 0$. If $\gl>0$ and $h=0$
instead, (\ref{eq:mainfe}) follows because
for $h \ge0$ we have
$\tf(\gl, 0) -2\gl h \le\tf(\gl, h) \le\tf(\gl, 0)$
by (\ref{eq:taumodel}) and (\ref{eq:fe-q}), hence for every $h>0$,
%
%e3.2 ###
%
\begin{eqnarray}
\widetilde\tf(\gl, h) &=&
\lim_{a \searrow0} \frac1{a^2} \tf(a\gl, ah)
\le
\liminf_{a \searrow0} \frac1{a^2} \tf(a\gl, 0)\nonumber\\
&\le&
\limsup_{a \searrow0} \frac1{a^2} \tf(a\gl, 0)
\le
\lim_{a \searrow0} \frac1{a^2} \tf(a\gl, ah) + 2\gl h \\
&=&
\widetilde\tf(\gl, h) + 2 \gl h \nonumber
\end{eqnarray}
so that (\ref{eq:mainfe}) for $h=0$
follows by continuity of $\widetilde\tf(\gl, \cdot)$.
For (\ref{eq:mainslope}), in view of (\ref{eq:mainslopeinf})
it suffices to show that
%
%e3.3 ###
%
\begin{equation}
\label{eq:mainslopesup}
\limsup_{\gl\searrow0} \frac{h_c(\gl)}{\gl} \le \slope_\ga,
\end{equation}
and Theorem \ref{th:main-tech} does yield (\ref{eq:mainslopesup}).
In fact if $c> \slope_\ga$, then $\widetilde\tf( (1+\rho) \gl,
(1-\rho) c \gl)=0$
for $\rho$ sufficiently small and every $\gl\ge0$;
the upper bound in (\ref{eq:main}) then
yields $\tf(a \gl, ac\gl)=0$ for $a$ small,
that is, $h_c(\gl) \le c \gl$ for $\gl$ small, which implies (\ref
{eq:mainslopesup}).

In order to carry out the proof Theorem \ref{th:main-tech}, it is
convenient to introduce
the following basic order relation.
\begin{definition} \label{def:main}
Let $f_{t,\gep,\delta}(a, \gl, h)$ and $g_{t,\gep,\delta}(a, \gl,
h)$ be two real
functions.
We write $f \prec g$ if for all fixed $\gl, h > 0$ and $\rho\in(0,1)$
there exists $\delta_0 > 0$ such that for every $0 < \delta< \delta
_0$ there exists
$\gep_0 = \gep_0(\delta) > 0$ such that for every $0 < \gep< \gep
_0$ there exists
$a_0 = a_0(\delta, \gep) > 0$ such that for every $0 < a < a_0$
%
%e3.4 ###
%
\begin{equation} \label{eq:order}
\limsup_{t \to\infty} f_{t,\gep,\delta}(a, \gl, h) \le
\limsup_{t \to\infty} g_{t,\gep,\delta} \bigl( a, (1+\rho)\gl,
(1-\rho)h \bigr).
\end{equation}
The values $\delta_0, \gep_0, a_0$ may also depend on $\gl, h ,\rho$.
If both $f \prec g$ and $g \prec f$, we write $f \simeq g$.
\end{definition}

Recalling the definitions (\ref{eq:fe-qa}) and (\ref{eq:existence})
of the discrete and continum finite-volume
free energies $\tf_N(\gl, h)$ and $\widetilde\tf_t(\gl, h)$, we set
%
%e3.5 ###
%
\begin{equation} \label{eq:psipsitilde}\quad
f^0_{t,\gep, \delta}(a, \gl, h) := \frac{1}{a^2}
\tf_{\lfloor t/a^2 \rfloor} (a\gl, a h) ,\qquad
f^4_{t,\gep, \delta}(a, \gl, h) := \widetilde\tf
_{t} (\gl, h) ,
\end{equation}
(that in fact do not depend on $\gep, \delta$ and on $\gep, \delta, a$).
Thanks to Definition \ref{def:main}, we see immediately that
proving Theorem \ref{th:main-tech} is equivalent to showing
that \mbox{$f^0 \simeq f^4$}.
Since the relation $\simeq$ is symmetric and transitive, we proceed
by successive approximations: more precisely, we are going to prove that
%
%e3.6 ###
%
\begin{equation} \label{eq:4steps}
f^0 \simeq f^1 \simeq f^2 \simeq f^3 \simeq f^4 ,
\end{equation}
where $f^{i} = f^i_{t,\gep, \delta}(a, \gl, h)$
for $i=1,2,3$ are suitable intermediate quantities.

The proof is divided into \textit{four steps}, corresponding to the equivalences
in (\ref{eq:4steps}). In each step,
%to prove that $f^{i-1} \simeq f^{i}$,
we will make statements that hold when $\delta$, $\gep$ and $a$ are small
in the sense prescribed by Definition \ref{def:main},
that is,
when $0 < \delta< \delta_0$, $0 < \gep< \gep_0(\delta)$ and $0 < a
< a_0(\delta, \gep)$,
for a suitable choice of $\delta_0$, $\gep_{0}(\cdot)$ and
$a_0(\cdot, \cdot)$.
For brevity, we will refer to this notion of smallness
by saying that $\gep$, $\delta$, $a$ are \textit{small in the usual sense}.
It is important to keep in mind that
%
%e3.7 ###
%
\begin{equation}
\label{eq:magn}
t^{-1} \ll a \ll
\gep \ll \delta \ll 1 .
\end{equation}
At times, we will
commit abuse of notation by writing $a_0(\gep)$
or $a_0(\delta)$ to stress the parameter that enters the specific computation.
In order to simplify notationally the proof,
we also assume that all the large numbers built with $\delta, \gep,
a, t$
that we encounter, such as $\gep/a^2$, $\delta/\gep$, $t/\delta$
(hence $\delta/a^2$, $t/\gep$, $t/a^2, \ldots$), are integers.

Before starting with the proof,
let us describe a general scheme, that is common to all the four steps.
The functions $f^{i}$ that we consider will always be of the form
%
%e3.8 ###
%
\begin{equation}
\label{eq:fform}
f^{i}_{t,\gep,\delta}(a, \gl, h) = \frac1 t \bbE \log
\bE
[ \exp(
-2 a \gl H^{i}_{t,\gep,\delta}(a, h) ) ] ,
\end{equation}
for a suitable Hamiltonian $H^{i}_{t,\gep,\delta}(a, h)$. Now, for
$\rho\in(0,1)$,
let us write
%
%e3.9 ###
%
\begin{equation}
\label{eq:DH}
H^{i}_{t,\gep,\delta}(a, h) = H^{j}_{t,\gep,\delta} \bigl(a,
(1-\rho) h \bigr) +
\gD H^{(i,j)}_{t,\gep,\delta}(a, h, \rho)
\end{equation}
(this relation is the definition of $\gD H$). Applying H\"older, Jensen and
Fubini, we get
%
%e3.10 ###
%
\begin{eqnarray}
\label{eq:Holder}\quad
&&
f^{i}_{t,\gep, \delta} (a, \gl, h) \nonumber\\
&&\qquad \le \frac{1}{1+\rho}
f^{j}_{t,\gep, \delta} \bigl(a, (1+\rho)\gl, (1-\rho)h \bigr)
\\
&&\qquad\quad{} + \frac{1}{(1+\rho^{-1}) t} \log\bE\bbE\exp\bigl(
-2 a (1+\rho^{-1}) \gl \gD H^{(i,j)}_{t,\gep, \delta}(a, h, \rho)
\bigr) .\nonumber
\end{eqnarray}
Therefore, to prove $f^{i} \prec f^{j}$ it suffices to show that for
every positive constant $A$ we can choose the parameters $\delta,\gep
, a$ small
in the usual sense such that
%
%e3.11 ###
%
\begin{equation}
\label{eq:ctrl}
\limsup_{t\to\infty} \frac1 t \log\bE\bbE\exp\bigl(
-a A \gD H^{(i,j)}_{t,\gep,\delta}(a, h,\rho) \bigr) \le 0 .
\end{equation}
Replacing $\gD H^{(i,j)}$ by $\gD H^{(j,i)}$ in this relation,
we prove that $f^{j} \prec f^{i}$ and therefore
$f^{i} \simeq f^{j}$.

%%%%%%%%%%%%%%%%%%%%%%%%%%%%%%%%%%%%%%%%%%%%%%%%%%%%%%%%%%%%
%%%%%%%%%%%%%%%%%%%%%%%%%%%%%%%%%%%%%%%%%%%%%%%%%%%%%%%%%%%%
%s3.1 ###
\subsection{Step 1: Coarse-graining of the renewal process}
\label{sec:step1}

We recall that by definition [see (\ref{eq:psipsitilde}), (\ref{eq:fe-qa})
and (\ref{eq:discZ})] $f^{0}$ is given by
%
%e3.12 ###
%
\begin{equation} \label{eq:f0}
f^{0}_{t,\gep, \delta}(a, \gl, h) := \frac{1}{a^2} \tf_{t/a^2}
(a\gl, a h)
= \frac1 t \bbE \log \bE[ \exp(
-2 a \gl H^{0}_{t,\gep,\delta}(a, h) ) ] ,\hspace*{-28pt}
\end{equation}
where $H^{0}$ is defined by
%
%e3.13 ###
%
\begin{equation} \label{eq:H0}
H^{0}_{t,\gep,\delta}(a, h) = \sum_{i=1}^{t/a^{2}} (\go_{i} +
ah) \gD_{i} .
\end{equation}
The purpose of this section is to define a first intermediate approximation
$f^{1}$ and to show that $f^{0} \simeq f^{1}$, in the sense of
Definition \ref{def:main}, following the general scheme
(\ref{eq:fform})--(\ref{eq:ctrl}).

We recall that the sequence $\gD_i \in\{0,1\}$ is constant
for $i \in\{\tau_j+1, \tau_j+2, \ldots, \tau_{j+1}\}$ and it is
chosen by flipping a fair coin.
%(in the simple random walk case $\gD_{i} = \ind_{\{S_{i} + S_{i-1} < 0
%than in $\{+1,-1\}$.
We start by defining,
for $j \in\N\cup\{0\}$, the basic coarse-grained blocks
%
%e3.14 ###
%
\begin{equation}
I_j :=
\bigl( (j-1) \gep/a^2, j \gep/a^2 \bigr].
\end{equation}
Then we set $\gs_0:=0$ and for $k \ge1$
%
%e3.15 ###
%
\begin{equation} \label{eq:gsk}
\gs_k := \inf\{j \ge\gs_{k-1}+( \delta/\gep):
\tau\cap I_j \neq\varnothing
\},
\end{equation}
thus introducing a coarse-grained version $\gs$
of the underlying renewal $\tau$
that has a resolution of $\gep/a^2\gg1$.
We say that the block $I_j$ is \textit{visited} if there exists $k$
such that
$\gs_k =j$. We stress that $\gs$ is built in such a way
that if $I_j$ is visited, we disregard the
content of the next $(\delta/\gep)-1 \gg1$ blocks, that is, we dub them
as \textit{not visited} (even if they may contain renewal points).
Since we are interested only in the blocks that
fall inside the interval $[0, t/a^2]$, we set
$m_{t/ a^2}:= \min\{ k\dvtx \gs_k \ge t/\gep\}$. Moreover,
for $k \in\N$, we give a notation for the union of blocks
between visited sites (that should be interpreted as \textit
{coarse-grained excursions}):
%
%e3.16 ###
%
\begin{equation}
\bar I_k := \Biggl( \bigcup_{j = \gs_{k-1} + 1}^{\gs_{k}} I_{j}
\Biggr) \cap( 0, t/a^2 ].
\end{equation}
Note that $\bar I_k \neq\varnothing$
if and only if $k \le\gs_{m_{t/a^2}}$; furthermore
$(0, t/a^{2}] = \bigcup_{k=1}^{m_{t/a^{2}}} \bar I_{k}$.
Each coarse-grained excursion $\bar I_{k}$ with $1 \le k < m_{t/a^{2}}$
contains exactly one visited block, namely $I_{\gs_{k}}$, at its
right extremity. The last
coarse-grained excursion $\bar I_{m_{t/a^{2}}}$ may or may not
end with a visited block, depending on whether $\gs_{m_{t/a^{2}}} =
t/\gep$
or $\gs_{m_{t/a^{2}}} > t/\gep$.

\begin{figure}

\includegraphics{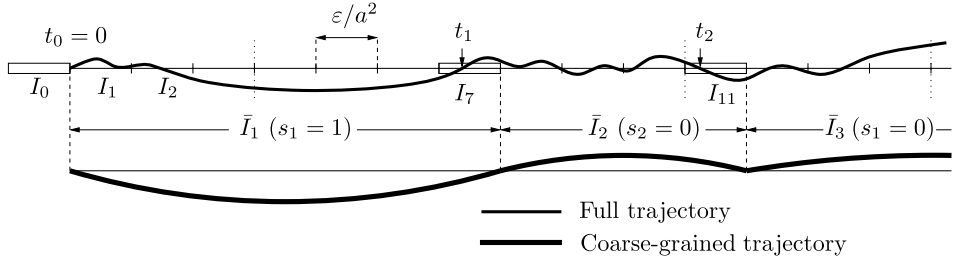}

\caption{A full trajectory, on top, and the corresponding
coarse-grained trajectory, below.
The visited blocks are surrounded by a box and the first renewal point
inside such
blocks is marked by a vertical arrow: a coarse-grained excursion is
everything that lies between
visited blocks.
One stipulates that there is a visited block to the left of the origin,
containing the origin.
The visited block on the right belongs to the coarse-grained excursion,
while the
one on the left does not. The sign of the excursion is just the sign
of the full trajectory just before the vertical arrow
(except possibly for the last excursion).
In this example $\delta/\gep=4$, so the first three blocks
to the right of a visited block
(i.e., up to the vertical dotted lines) cannot be visited
blocks.}\vspace*{-4pt}
\label{fig:cg}
\end{figure}

For $1 \le k < m_{t/a^2}$, we assign a
\textit{sign} $s_k$ to the $k$th coarse-grained excursion
by stipulating that it coincides with the \textit{sign}
just before the first renewal point in $I_{\gs_k}$ (that we call
$t_k$, and $t_0:=0$),
that is, we set $s_{k} := \gD_{t_{k}}$. When $k=m_{t/a^2}$, we need to make
a distinction: if the coarse-grained excursion $\bar I_{k}$
ends with a visited block ($\gs_{m_{t/a^{2}}} = t/\gep$)
we set $s_{k} := \gD_{t_{k}}$ as before; if the coarse-grained
excursion $\bar I_{k}$ is truncated ($\gs_{m_{t/a^{2}}} > t/\gep$)
we set $s_k= \gD_{t/a^2}$.
We refer to Figure \ref{fig:cg}
for a graphical description of the quantities introduced so far.

%
%f3 ###
%

We are now ready to introduce the first intermediate approximation $f^{1}$.
According to (\ref{eq:fform}), it suffices to define the corresponding
Hamiltonian
%
%e3.17 ###
%
\begin{equation} \label{eq:H1}\quad
H^1_{t, \gep, \delta}( a,h) :=
\sum_{k=1}^{m_{t/a^2}} \sum_{i \in\bar I _k} (\go_i +ah) s_{k}
= \sum_{k=1}^{m_{t/a^2}} s_k \bigl( Z_k(\go) +ah \vert\bar I_k
\vert\bigr) ,
\end{equation}
where $Z_k(\go):= \sum_{i \in\bar I_k} \go_i$.
Note that we may rewrite $H^{0}$ [see (\ref{eq:H0})] as
%
%e3.18 ###
%
\begin{equation}
H^0_{t, \gep, \delta}( a,h)
% := \sum_{i=1}^{t/a^2} (\go_i +ah) \gD_i
= \sum_{k=1}^{m_{t/a^2}} \sum_{i \in\bar I _k} (\go_i +ah)\gD
_i .
\end{equation}
Passing from $H^{0}$ to $H^{1}$, we are thus replacing the renewal
$\tau$ by its coarse-grained version.
Applying the\vspace*{1pt} general scheme (\ref{eq:fform})--(\ref{eq:Holder}),
to prove that $f^0 \simeq f^1$ we have to establish (\ref{eq:ctrl}) for
$\gD H^{(0,1)}$ and $\gD H^{(1,0)}$, defined by
%
%e3.19 ###
%
\begin{eqnarray}
\label{eq:wtonw}\hspace*{32pt}
\gD H^{(0,1)}_{t, \gep, \delta}( a,h, \rho)  :\!&=&
H^0_{t, \gep, \delta}( a,h) - H^1_{t, \gep, \delta}\bigl( a,(1-\rho
)h\bigr) \nonumber\\[-8pt]\\[-8pt]
&=& a \rho h \sum_{i=1}^{t/a^2}\gD_i +
\sum_{k=1}^{m_{t/a^2}} \sum_{i \in\bar I_k}
\bigl( \go_i +a(1-\rho)h \bigr)
( \gD_i -s_k ) ,\nonumber
\end{eqnarray}
and
%
%e3.20 ###
%
\begin{eqnarray}
\label{eq:wtonw2}\quad
\gD H^{(1,0)} _{t, \gep, \delta}( a,h, \rho) :\!&=&
H^1_{t, \gep, \delta}( a,h) - H^0_{t, \gep, \delta}\bigl( a,(1-\rho
)h\bigr) \nonumber\\[-8pt]\\[-8pt]
&=& a \rho h \sum_{i=1}^{t/a^2}\gD_i +
\sum_{k=1}^{m_{t/a^2}} \sum_{i \in\bar I_k}
( \go_i +ah ) ( s_k-\gD_i ) .\nonumber
\end{eqnarray}
Formulas (\ref{eq:wtonw}) and (\ref{eq:wtonw2}) are minimally different:
in particular we are going to estimate the second term in the
right-hand side
by taking the absolute value. For this reason, we detail
only the case of (\ref{eq:wtonw}).

In order to establish (\ref{eq:ctrl}) for $\gD H^{(0,1)}$,
we observe that for $a \le t_0/A^2$ [$t_0$ is the constant
in (\ref{eq:omcond})]
%
%e3.21 ###
%
\begin{eqnarray}
\label{eq:4.10}\qquad
&&\bbE e^{-A a \gD H^{(0,1)}}\nonumber\\
&&\qquad= \bbE\exp\Biggl( -Aa^2 \rho h \sum_{i=1}^{t/a^2}\gD_i -Aa
\sum_{k=1}^{m_{t/a^2}} \sum_{i \in\bar I_k}
\bigl( \go_i +a(1-\rho)h \bigr) ( \gD_i -s_k )
\Biggr) \nonumber\\
&&\qquad= \exp\Biggl(
-Aa^2 \rho h \sum_{i=1}^{t/a^2}\gD_i -Aa^2(1-\rho)h
\sum_{k=1}^{m_{t/a^2}} \sum_{i \in\bar I_k}
( \gD_i -s_k ) \Biggr)\\
&&\qquad\quad{}\times \prod_{k=1}^{m_{t/a^2}} \prod_{i
\in\bar I_k}\M\bigl( Aa ( \gD_i -s_k ) \bigr) \nonumber\\
&&\qquad= \exp\Biggl( -C a^2 \sum_{i=1}^{t/a^2}\gD_i + B a^2
\sum_{k=1}^{m_{t/a^2}} \sum_{i \in\bar I_k} \vert\gD_i
-s_k \vert\Biggr),\nonumber
\end{eqnarray}
where $C := A \rho h$ and $B:=A(1-\rho)h+ c_0 A^2$. Here, we have used
(\ref{eq:omcond}) and the fact that $|\gD_{i} - s_{k}|^{2} = |\gD
_{i} - s_{k}|$
because $|\gD_{i} - s_{k}| \in\{0,1\}$.
This shows that (\ref{eq:ctrl}) is proven if we can show that for any given
$B,C>0$ we have
%
%e3.22 ###
%
\begin{equation}
\label{eq:ctrl1}\quad
\limsup_{t \to\infty}
\frac1t \log\bE\exp\Biggl(
-Ca^2 \sum_{i=1}^{t/a^2}\gD_i + B a^2
\sum_{k=1}^{m_{t/a^2}} \sum_{i \in\bar I_k}
\vert\gD_i -s_k \vert
\Biggr) \le 0 ,
\end{equation}
for $\delta$, $\gep$, $a$ small in the usual sense
[recall the discussion before (\ref{eq:magn})].

Let us re-express (\ref{eq:ctrl1}) explicitly in terms of
the renewal process $\tau$ and of the signs $\xi= \{\xi_j\}_{j\in\N}$,
where $\xi_j = \gD_{\tau_j}$. This notation\vspace*{2pt} has been already
introduced in Section \ref{sec:discmodel}:
here we need also $\cN_s := \vert\tau\cap[0,s] \vert
= \min\{k \ge1\dvtx \tau_k > s \}$ ($s \in\N$).
Observe that $\xi$ is an i.i.d. sequence, as well as the
sequence of the inter-arrivals $\{\eta_{j} := \tau_{j} - \tau_{j-1}\}
_{j\in\N}$.
First of all,
%
%e3.23 ###
%
\begin{equation}\quad
\sum_{i=1}^{t/a^2}\gD_i = \sum_{j=1}^{\cN_{t/a^2}-1} \xi_j\eta
_j +
\xi_{\cN_{t/a^2}} \bigl( (t/a^2)- \tau_{\cN_{t/a^2} -1} \bigr)
\ge \sum_{j=1}^{\cN_{t/a^2} -1} \xi_j\eta_j .
\end{equation}
Concerning the second addendum in the exponent in (\ref{eq:ctrl1}),
we use the fact that if $\eta_j= \tau_j-\tau_{j-1}
\ge(\delta/\gep)(\gep/a^2)= \delta/a^2$, then necessarily the
inter-arrival $\eta_j$
determines a coarse-grained excursion (say, $\bar I_k$).
We can then distinguish two cases: either $\tau_{j-1} \in\bar I_k$, or
$\tau_{j-1} \in\bar I_{k-1}$. If $\tau_{j-1} \in\bar I_k$,
we know that $\gD_i =s_k$ for every $i\in\{ \tau_{j-1}+1,\ldots,
\tau_j\}$,
by our definition of the sign of the coarse-grained excursions.
If, on the other hand, $\tau_{j-1} \in\bar I_{k-1}$, which happens
if and only if $\tau_{j-1} \in I_{\sigma_{k-1}}$, we can only be sure that
$\gD_i =s_k$ for every $i\in\{ \tau_{j-1}+1,\ldots, \tau_j\}
\setminus I_{\gs_{k-1}}$.
Since $|I_{\gs_{k-1}}| = \frac{\gep}{a^2}$ and there are $m_{t/a^2}$
visited blocks,
we are lead to the bound
%
%e3.24 ###
%
\begin{equation} \label{eq:hihi}
\sum_{k=1}^{m_{t/a^2}} \sum_{i \in\bar I_k} \vert\gD_i
-s_k \vert
\le \sum_{j=1}^{\cN_{t/a^2}-1} \eta_j \ind_{\eta_j < \delta
/a^2} +
\frac{\gep}{a^2} m_{t/a^2} .
\end{equation}
This step of the proof is therefore completed by applying the following lemma.
\begin{lemma}
\label{th:lem1}
For every $B,C>0$, we have
%
%e3.25 ###
%
\begin{eqnarray}
&&\limsup_{t \to\infty} \frac1t
\log\bE\exp\Biggl( Ba^2
\sum_{j=1}^{\cN_{t/a^2}-1} \eta_j \ind_{\eta_j <
\delta/a^2}\nonumber\\[-8pt]\\[-8pt]
&&\qquad\hspace*{67.1pt}{}+ B \gep m_{t/a^2}
- C a^2 \sum_{j=1}^{\cN_{t/a^2}-1} \xi_j\eta_j
\Biggr) \le 0,\nonumber
\end{eqnarray}
for $\delta$, $\gep$ and $a$ small in the usual sense.
\end{lemma}
\begin{pf}
%Let us start by observing that
%we can safely drop $j = \cN_{t/a^2} $ in the first sum in
%since it gives an additive contribution that is smaller than $B\gd/t$.
Since $\xi$ and $\eta$ are independent and since $\xi$ is an i.i.d.
sequence of $B(1/2)$ variables,
%
%e3.26 ###
%
\begin{eqnarray}\quad
&&\bE\exp\Biggl( Ba^2
\sum_{j=1}^{\cN_{t/a^2}-1} \eta_j \ind_{\eta_j < \delta/a^2}
+ B \gep m_{t/a^2}
- C a^2 \sum_{j=1}^{\cN_{t/a^2}-1} \xi_j\eta_j
\Biggr)
\nonumber\\
&&\qquad=\bE\exp\Biggl( Ba^2
\sum_{j=1}^{\cN_{t/a^2}-1} \eta_j \ind_{\eta_j < \delta/a^2}
\\
&&\qquad\quad\hspace*{29.5pt}{}+ B \gep m_{t/a^2}
+ \sum_{j=1}^{\cN_{t/a^2}-1}\log\biggl( \frac12 + \frac12
\exp(- C a^2 \eta_j ) \biggr)
\Biggr).\nonumber
\end{eqnarray}
The proof now proceeds in two steps: first, we will show that if
$\delta$, $\gep$ and
$a$ are small in the usual sense,
%
%e3.27 ###
%
\begin{equation}
\label{eq:step1.1}
B \gep m_{t/a^2} + \frac12 \sum_{j=1}^{\cN_{t/a^2}
-1}\log\biggl( \frac12 + \frac12 \exp(- C a^2 \eta_j )
\biggr) \le B \gep,
\end{equation}
uniformly in $\eta$, and then that
%
%e3.28 ###
%
\begin{eqnarray}
\label{eq:step1.2}
&&
\limsup_{t \to\infty} \frac1t\log
\bE\Biggl( Ba^2
\sum_{j=1}^{\cN_{t/a^2}-1} \eta_j \ind_{\eta_j < \delta/a^2}
\nonumber\\[-8pt]\\[-8pt]
&&\qquad\hspace*{49pt}{} + \frac12
\sum_{j=1}^{\cN_{t/a^2}-1}\log\biggl( \frac12 + \frac12 \exp
(- C a^2 \eta_j ) \biggr)
\Biggr) \le 0.\nonumber
\end{eqnarray}

For the proof of (\ref{eq:step1.1}), recall first that $t_k$ is the
first contact in $I_{\gs_k}$ for
$k < m_{t/a^2}$,
that is, $t_k:= \min\{ n \in I_{\gs_k}\dvtx n \in\tau\}$. Now, let
us consider
the intervals $(t_{k-1}, t_k]$ for $k=1, \ldots, m_{t/a^2} -1$
($t_0:=0$). Given a value of $k$:
\begin{enumerate}[(1)]
\item[(1)] either in $(t_{k-1}, t_k]$ there is a \textit{long excursion},
that is,
there exists $j^{*}$ such that $(\tau_{j^{*}-1},\tau_{j^{*}}]\subset
(t_{k-1}, t_k]$
with $\tau_{j^{*}}-\tau_{j^{*}-1}\ge\delta/a^2$,
so that
%
%e3.29 ###
%
\begin{eqnarray}\label{eq:parapa0}
&& B\gep+\frac12 \sum_{j\dvtx (\tau_{j-1},\tau_j]\subset(t_{k-1},
t_k] }
\log\biggl( \frac12 + \frac12 \exp(- C a^2 \eta_j )
\biggr) \nonumber\\
&&\qquad \le B\gep+\frac12
\log\biggl( \frac12 + \frac12 \exp(- C a^2 \eta
_{j^{*}} ) \biggr)\\
&&\qquad \le
B\gep+\frac12
\log\biggl( \frac12 + \frac12 e^{-C \delta} \biggr)
\le 0 ,\nonumber
\end{eqnarray}
where the last inequality holds for $\gep\le\gep_0(\delta)$;

\item[(2)] or in $(t_{k-1}, t_k]$ there are only \textit{short excursions},
that is,
$\eta_{j} := \tau_j-\tau_{j-1}< \delta/a^2$
for every $j$ such that $(\tau_{j-1},\tau_j]\subset(t_{k-1}, t_k]$.
In this case, we bound from above
$\log( \frac12 + \frac12 \exp(- C a^2 \eta_j )
)$ by
$-\frac14 Ca^2 \eta_j$ for $\delta\le\delta_0$, so that
%
%e3.30 ###
%
\begin{eqnarray} \label{eq:parapa1}
&&
B\gep+\frac12 \sum_{j\dvtx (\tau_{j-1},\tau_j]\subset(t_{k-1}, t_k] }
\log\biggl( \frac12 + \frac12 \exp(- C a^2 \eta_j )
\biggr) \nonumber\\[-8pt]\\[-8pt]
&&\qquad\le B\gep-\frac18 Ca^2(t_k-t_{k-1})
\le 0 ,\nonumber
\end{eqnarray}
where the last inequality holds for $\gep\le\gep_0(\delta)$ and it
follows by observing that
$t_k-t_{k-1}> ((\delta/\gep)-1)(\gep/a^2)= (\delta- \gep)/a^2$.
\end{enumerate}
Summing (\ref{eq:parapa0}) and (\ref{eq:parapa1}) from $k=1$
to $k=m_{t/a^{2}} - 1$, we see that (\ref{eq:step1.1}) holds true.

Let us therefore turn to (\ref{eq:step1.2}): note that we need to estimate
%
%e3.31 ###
%
\begin{eqnarray}
\label{eq:step1.2.1}
\frac1t\log\bE\exp\Biggl(
\sum_{j=1}^{\cN_{t/a^2}-1} g(a^{2} \eta_j) \Biggr)\nonumber\\[-8pt]\\[-8pt]
\eqntext{\mbox{ with }
g(x) := B x \ind_{x< \delta} + \frac12
\log\bigl( \frac12 + \frac12 e^{- C x} \bigr) .}
\end{eqnarray}
Since $g(\cdot) \ge- \frac12 \log2$, we can add the term $j=\cN_{t/a^2}$
by paying at most $\sqrt{2}$, that is,
%
%e3.32 ###
%
\begin{eqnarray} \label{eq:truc}
\bE\exp\Biggl(
\sum_{j=1}^{\cN_{t/a^2} -1 } g(a^{2} \eta_j) \Biggr) \le
\sqrt{2} \bE[ G_{\cN_{t/a^{2}}} ] \nonumber\\[-8pt]\\[-8pt]
\eqntext{\mbox{where }
\displaystyle G_{n} := \exp\Biggl(\sum_{j=1}^{n} g(a^{2} \eta_j) \Biggr).}
\end{eqnarray}
Let us set $G_{0} := 1$ and $\gga:= \bE[\exp(g(a^{2} \eta_{1}))]$
for convenience. Since $G_{n}$ is the product of $n$ i.i.d. random variables,
the process $\{G_{n}/\gga^{n}\}_{n \ge0}$ is a martingale
(with respect to the natural filtration of the sequence $\{\tau_{n}\}
_{n \ge0}$).
Assume now that $\gga\le1$: the process $\{G_{n}\}_{n \ge0}$ is a
supermartingale
and, since $\cN_{t/a^2}$ is a bounded stopping time, the optional sampling
theorem yields $\bE[ G_{\cN_{t/a^{2}}} ] \le1$. Then from (\ref
{eq:truc}) it
follows immediately that (\ref{eq:step1.2}) holds, thus completing
the proof.

We are left with showing that $\gga\le1$, that is,
$\bE[\exp(g(a^{2} \eta_1)] \le1$, when $\delta$, $\gep$ and $a$ are
small in the usual sense (actually $\gep$ does not appear in
this quantity). Note that
%
%e3.33 ###
%
\begin{equation}
\bE[\exp(g(a^{2 }\eta_1))]-1 = \sum_{n\in\N} [\exp
(g(a^{2}n)) - 1 ]
K(n) ,
\end{equation}
and recall that $K(n) \sim L(n) / n^{1+\ga}$ as $n\to\infty$,
with $L(\cdot)$ slowly varying at infinity.
Then it follows by Riemann sum approximation that
%
%e3.34 ###
%
\begin{eqnarray}
\label{eq:Riemann1}
&&
\lim_{a \searrow0} \frac{\bE[\exp(g(a^{2}\eta_1))]-1}{a^{2\ga}L(1/a^2)}
\nonumber\\[-8pt]\\[-8pt]
&&\qquad= \int_0^\infty\biggl[ \exp\biggl( B x \ind_{x <\delta} +
\frac12
\log\biggl( \frac12 + \frac12 e^{-C x} \biggr) \biggr) - 1\biggr]\,
\frac{\dd x}{x^{1+\ga}} .\nonumber
\end{eqnarray}
The Riemann sum approximation is justified since
$L(cn)/L(n) \to1$ as $n \to\infty$
uniformly for $c$ in compact sets of $(0,\infty)$ \cite
{cfBinGolTeu}, Theorem 1.2.1,
and since
for every $\epsilon>0$ there exists $b>0$ such that $L(n) \le b
n^\epsilon$ for every $n$
(the latter property is used to deal with very large and small values
of $n$).
A simple look at (\ref{eq:Riemann1}) suffices to see that
the right-hand side is negative if $\delta\le\delta_0$.
%Let us now choose $\gd\le\gd_0$ and $a\le a(\gd)$ such that
%$\bE[\exp(g(\eta_1)]\le1$ and let us go back to \eqref{eq:step1.2.1}
%and observe that,
%since $g(\cdot)$ is bounded,
%$\cN_{t/a^2}-1$ can be replaced by $\cN_{t/a^2}$. Moreover there
%exists $c>0$,
%depending uniquely on $K(\cdot)$,
%such that
% \log\bE\exp(
% ) \le
% ); t/ a^2 \in\tau] + c \log( t/a^2 ),
%in strict analogy with ADD. Therefore it is sufficient to
%estimate the first term in the right-hand side:
% ); N \in\tau] =
%where $\widetilde\tau$ is a renewal process with inter-arrival law
%given by $\widetilde K(n) := K(n)\exp(g(n))$, and $\sum_n
%This completes the proof of Lemma \ref{th:lem1}.
\end{pf}

%%%%%%%%%%%%%%%%%%%%%%%%%%%%%%%%%%%%%%%%%%%%%%%%%%%%%%%%%%%%
%%%%%%%%%%%%%%%%%%%%%%%%%%%%%%%%%%%%%%%%%%%%%%%%%%%%%%%%%%%%
%s3.2 ###
\subsection{Step 2: Switching to Gaussian charges}

In this step, we introduce the second intermediate approximation $f^{2}$:
following (\ref{eq:fform}), we define the corresponding Hamiltonian
$H^{2}$ by
%
%e3.35 ###
%
\begin{equation}
\label{eq:H2}
H^2_{t, \gep, \delta}( a,h) := \sum_{k=1}^{m_{t/a^2}} s_k
\bigl( Z_k(\widehat\go) +ah \vert\bar I_k \vert\bigr),
\end{equation}
where $\widehat\go=\{ \widehat\go_i\}_{i \in\N}$ is an i.i.d.
sequence of standard Gaussian
random variables and we recall that $Z_k(\widehat\go):= \sum_{i \in
\bar I_k} \widehat\go_i$.
We stress that,
with respect to the preceding Hamiltonian $H^{1}$ [cf. (\ref
{eq:H1})] we have
just changed the charges $\go_{i} \to\widehat\go_{i}$.\vadjust{\goodbreak}

In order to apply the general scheme (\ref{eq:fform})--(\ref{eq:ctrl}),
we build the two sequences of disorder variables
$\go=\{ \go_i\}_{i \in\N}$ and $\widehat\go=\{ \widehat\go_i\}
_{i \in\N}$
on the same probability space $(\Omega, \cF, \bbP)$, that is, we
define a
\textit{coupling}. Actually, the disorder does not appear any longer
in terms of the individual charges $\go_i$, but it is by now
summed over the
coarse-grained blocks $I_j = ((j-1) \frac{\gep}{a^{2}}, j \frac{\gep
}{a^{2}}]$, so we\vspace*{2pt} just
need to
couple the two i.i.d. sequences $\{ \sum_{i\in I_j} \go_i\}_{j\in\N
}$ and
$\{ \sum_{i\in I_j} \widehat\go_i\}_{j\in\N}$. The coupling is
achieved via
the standard Skorohod representation in the following way:
given the i.i.d. sequence $\{ \widehat\go_i\}_{i \in\N}$
of $\cN(0,1)$ variables, if we set $ \widehat F(t) := \bbP(\widehat
\go_1 \le t)$
and $n := \vert I_1\vert$,
then $ \widehat F ( \sum_{i\in I_j} \widehat\go_i / \sqrt{n}
) = : U_j$
is uniformly distributed over $(0,1)$. Therefore, if we set
$F_n(t):= \bbP( \sum_{i\in I_j} \go_i / \sqrt{n} \le t)$ and
$F_n^{-1}(s) := \inf\{t \in\R\dvtx F_{n}(t) > s\}$, that is, $F_n^{-1}$
is the generalized inverse of $F_n$,
then the sequence $\{F_n^{-1} (U_j)\}_{j\in\N}$
has the same law as $\{ \sum_{i\in I_j} \go_i / \sqrt{n}\}_{j\in\N}$
and we have built a coupling.
In short, we set $X^{(n)}_j := F_n^{-1}(U_j)$ and $Y_j := \widehat
F^{-1} (U_j)=
\sum_{i\in I_j} \widehat\go_i / \sqrt{n}$.
Moreover, we observe that, by the central limit theorem, $\lim_{n\to
\infty} F_{n}(t) = \widehat F(t)$
for every $t\in\R$ and
therefore $\lim_{n\to\infty} X^{(n)}_{j} = Y_{j}$, in $\bbP$-probability.
%we assume that on $\Omega$ is defined a sequence of IID variables $\{
%U_j\}_{j \in\N}$ uniformly distributed on $(0,1)$. For $n := |I_{1}|
%= \gep/a^2$, we set $F_{n} (t):= \bbP( \sum_{i\in I_1} \go_i \le t
%quantity for $\widehat\go$, which is of course the distribution
%function of a standard normal variable.
%Then we define $X^{(n)}_j := F_n^{-1}(U_j)$, where $F_n^{-1}(s) := \inf
%$Y_j := \widehat F^{-1} (U_j)$, so that $\{X^{(n)}_j\}_j$ has the same
%law as $\{ \sum_{i\in I_j} \go_i/\sqrt{\vert I_j \vert}\}_{j}$ and$\{
%Y_j\}_j$ has the same law as $\{ \sum_{i\in I_j} \widehat\go_i/\sqrt{
%
\begin{lemma}
\label{th:coupling}
For every $C>0$,
%
%e3.36 ###
%
\begin{equation}
\lim_{n \to\infty}
\bbE\bigl[ \exp\bigl( C \bigl\vert
X^{(n)}_1
-Y_1 \bigr\vert\bigr) \bigr] = 1.
\end{equation}
\end{lemma}
\begin{pf}
Since $\lim_{n\to\infty} X^{(n)}_1 = Y_1$ in probability it suffices
to prove
that the sequence of random variables
$\{\exp( C \vert X^{(n)}_1 - Y_1 \vert)\}_{n\in n_0+\N}$
is bounded in $L^{2}$ (hence, uniformly integrable) for a given $n_0\in
\N$.
We choose $n_0$ to be the smallest integer number larger than
$ 16C^2/ t_0^2$, with $t_0$ the constant in (\ref{eq:omcond}).
By the triangle and Cauchy--Schwarz inequalities, we get
%
%e3.37 ###
%
\begin{eqnarray}
&&
\sup_{n> n_0}
\bbE\bigl[ \exp\bigl( 2C \bigl\vert X^{(n)}_1 - Y_1 \bigr\vert
\bigr) \bigr]
\nonumber\\[-8pt]\\[-8pt]
&&\qquad\le\sqrt{ \Bigl( \sup_{n> n_0} \bbE\bigl[ \exp\bigl( 4C
\bigl\vert X^{(n)}_1
\bigr\vert\bigr) \bigr] \Bigr)
\bbE[ \exp( 4C \vert Y_1 \vert) ] }
< \infty,\nonumber
\end{eqnarray}
where the second inequality follows from (\ref{eq:omcond}) and the
choice of $n_0$,
recalling that $X_1^{(n)} \sim\sum_{i=1}^n \go_i / \sqrt{n}$ and
$Y_1 \sim\cN(0,1)$.
\end{pf}

Let us see why Lemma \ref{th:coupling} implies
$f^1 \simeq f^2$. First of all,
\begin{eqnarray*}
&&\min\bigl(
H^1_{t, \gep, \delta}( a,h) - H^2_{t, \gep, \delta}\bigl( a,(1-\rho)h\bigr)
,
H^2_{t, \gep, \delta}( a,h) - H^1_{t, \gep, \delta}\bigl( a,(1-\rho)h\bigr)
\bigr) \\
&&\qquad \ge - \sum_{k=1}^{m_{t/a^2}} s_k
\vert Z_k(\go)-Z_k(\widehat\go) \vert +
a\rho h \sum_{k=1}^{m_{t/a^2}} s_k \vert\bar I_k \vert
\\
&&\qquad
\ge - \sum_{k=1}^{m_{t/a^2}} s_k \sum_{j = \gs_{k-1}+1}^{\gs_k}
\biggl\vert\sum_{i \in I_j} \go_i - \sum_{i \in I_j} \widehat
\go_i
\biggr\vert + a\rho h \sum_{k=1}^{m_{t/a^2}} s_k \vert
\bar I_k
\vert,
\end{eqnarray*}
where we redefine $\gs_{m_{t/a^{2}}} := t/\gep$ for
notational convenience (otherwise, we should treat the
last term $j=m_{t/a^2}$ separately).
In view of (\ref{eq:DH})--(\ref{eq:ctrl}), it suffices to show that
for $a, \gep$ and $\delta$ small in the usual sense
[recall the discussion before (\ref{eq:magn})] we have
%
%e3.38 ###
%
\begin{eqnarray} \label{eq:fkeg2}\hspace*{32pt}
&&\limsup_{t \to\infty} \frac1t \log\bE\Biggl[
\exp\Biggl(
-
Aa^2\rho h \sum_{k=1}^{m_{t/a^2}} s_k \vert
\bar I_k \vert\Biggr)
\nonumber\\
&&\qquad\hspace*{46.46pt}{}\times\bbE\Biggl( \exp\Biggl( A a
\sum_{k=1}^{m_{t/a^2}} s_k \sum_{j = \gs_{k-1}+1}^{\gs_k}
\biggl( \frac{\sqrt{\gep}}{a} \biggr)
\bigl\vert X_j^{(\gep/a^2)}- Y_j \bigr\vert
\Biggr) \Biggr)
\Biggr] \\
&&\qquad\le 0 .\nonumber
\end{eqnarray}
By independence,
%
%e3.39 ###
%
\begin{eqnarray}
&&\bbE\Biggl[ \exp\Biggl( A a
\sum_{k=1}^{m_{t/a^2}} s_k \sum_{j = \gs_{k-1}+1}^{\gs_k}
\biggl( \frac{\sqrt{\gep}}{a} \biggr)
\bigl\vert X_j^{(\gep/a^2)}- Y_j \bigr\vert\Biggr) \Biggr]
\nonumber\\[-8pt]\\[-8pt]
&&\qquad= \prod_{k=1}^{m_{t/a^2}} \bbE\bigl[ \exp\bigl( A\sqrt{\gep} s_k
\bigl\vert X^{(\gep/a^2)}_1-Y_1 \bigr\vert\bigr) \bigr] ^{\gs
_k -\gs_{k-1}}
,\nonumber
\end{eqnarray}
and since $a^2 \vert\bar I_k\vert= \gep(\gs_k -\gs_{k-1}
)$ the
term between square brackets in (\ref{eq:fkeg2}) is equal to
%
%e3.40 ###
%
\begin{equation}
\prod_{k=1}^{m_{t/a^2}} \bigl( \exp( -A \rho h s_k \gep)
\bbE\bigl[ \exp\bigl( A\sqrt{\gep} s_k
\bigl\vert X^{(\gep/a^2)}_1-Y_1 \bigr\vert
\bigr) \bigr] \bigr)^{\gs_k -\gs_{k-1}} .
\end{equation}
Since $s_{k} \in\{0,1\}$, (\ref{eq:fkeg2}) is implied by
%
%e3.41 ###
%
\begin{equation}
\exp( -A \rho h \gep) \bbE\bigl[ \exp\bigl(
A\sqrt{\gep}
\bigl\vert X^{(\gep/a^2)}_1-Y_1 \bigr\vert\bigr) \bigr] \le
1 ,
\end{equation}
which holds for $a \le a_0(\gep)$ by Lemma \ref{th:coupling}.
The proof of $f^1 \simeq f^2$ is complete.

%%%%%%%%%%%%%%%%%%%%%%%%%%%%%%%%%%%%%%%%%%%%%%%%%%%%%%%%%%%%%%%
%%%%%%%%%%%%%%%%%%%%%%%%%%%%%%%%%%%%%%%%%%%%%%%%%%%%%%%%%%%%%%%

%s3.3 ###
\subsection{Step 3: From the renewal process to the regenerative set}
\label{sec:step3}

In this crucial step, we replace the discrete renewal process
$\tau= \{\tau_{n}\}_{n\in\N}$ with
the continuum regenerative set $\widetilde\tau^\ga$
(both processes are defined under the law $\bP$).
Recall that for the renewal process $\tau$ we have defined
the coarse-grained returns $\{\gs_{k}\}_{k\in\N}$ as well as the
coarse-grained signs $s_{k}$, and $m_{t/a^{2}} :=
\inf\{k\dvtx \gs_{k} \ge t/\gep\}$. Henceforth,
we set $m := m_{t/a^{2}}$ for short\vspace*{1pt} and we redefine
for notational convenience $\gs_{m} := t/\gep$ (as in the previous step).

Since $\bar{I}_{k} = ( \frac{\gep}{a^{2}} \gs_{k-1}, \frac{\gep
}{a^{2}} \gs_{k}]$,
the second intermediate Hamiltonian $H^{2}$ [cf. (\ref{eq:H2})]
can be rewritten as
%
%e3.42 ###
%
\begin{equation} \label{eq:H2b}\qquad\quad
H^{2}_{t,\gep,\delta}(a,h) = \frac{1}{a} \sum_{k=1}^{m} s_{k}
\biggl[
\biggl( \sum_{({\gep\gs_{k-1}})/{a^{2}} < i \le{\gep\gs
_{k}}/{a^{2}}}
a \widehat\go_{i} \biggr) + h \gep(\gs_{k} - \gs_{k-1}) \biggr].
\end{equation}
We now introduce the rescaled returns $\underline\gs_{k} := \gep\gs
_{k}$ and
we let $\gb= \{\gb_{t}\}_{t \ge0}$ be a standard Brownian motion,
defined on the disorder probability space $(\Omega, \cF, \bbP)$.
With some abuse of notation, we can redefine $H^{2}$ as
%
%e3.43 ###
%
\begin{equation} \label{eq:H2mod}
H^{2}_{t,\gep,\delta}(a,h) = \frac{1}{a} \sum_{k=1}^{m} s_{k}
\bigl( \gb_{\underline\gs_{k}} - \gb_{\underline\gs_{k-1}}
+ h (\underline\gs_{k} - \underline\gs_{k-1}) \bigr),
\end{equation}
which has the same law as the quantity in (\ref{eq:H2b}), hence
through formula (\ref{eq:fform}) it yields \textit{the same} $f^{2}$.
It is clear that $H^{2}$ depends on the renewal
process $\tau= \{\tau_{n}\}_{n\in\N}$ only through the vector
%
%e3.44 ###
%
\begin{equation} \label{eq:Sigma}
\Sigma:= (m; s_{1}, \ldots, s_{m}; \underline\gs_{1},
\ldots, \underline\gs_{m}) ,
\end{equation}
whose definition depends of course on $t, a, \gep, \delta$.

One can define an analogous vector $\widetilde\Sigma$
in terms of the regenerative set $\widetilde\tau^\ga$, by looking
at the returns on blocks of width $\gep$, skipping
$(\delta/\gep)$ blocks between successive returns. More precisely, we set
$\widetilde I_{j} := ((j-1)\gep, j\gep]$ for $j \in\N$ and define
%
%e3.45 ###
%
\begin{eqnarray} \label{eq:tildegs}
\widetilde{\underline\gs}_{0} &:=& 0 ,\nonumber\\[-8pt]\\[-8pt]
\widetilde{\underline\gs}_{k} &:=&
\gep\cdot\inf\{ j \ge(\widetilde{\underline\gs}_{k-1}/\gep
) + (\delta/\gep) \dvtx
\widetilde\tau^\ga\cap\widetilde I_{j} \ne\varnothing\} ,\qquad
n\in\N.\nonumber
\end{eqnarray}
We then\vspace*{2pt} set $\widetilde m := \inf\{k \in\N\dvtx \widetilde{\underline
\gs}_{k} \ge t \}$
and redefine $\widetilde{\underline\gs}_{\widetilde m} := t$. The signs
$\{\widetilde s_{k}\}_{1 \le k \le m}$ are defined in complete analogy
with the discrete case, by looking at the sign $\widetilde\gD^{\ga}$
at the
beginning of each visited block $\widetilde I_{\widetilde\gs_{k}}$.
We have thus completed the
definition of
%
%e3.46 ###
%
\begin{equation} \label{eq:Sigmatilde}
\widetilde\Sigma:= (\widetilde m; \widetilde s_{1}, \ldots,
\widetilde s_{\widetilde m};
\widetilde{\underline\gs}_{1}, \ldots, \widetilde{\underline\gs
}_{\widetilde m}) .
\end{equation}

We are ready to introduce the third intermediate quantity $f^{3}$,
which, in agreement with (\ref{eq:fform}),
will be defined by the corresponding Hamiltonian $H^{3}$.
We replace in the right-hand side of (\ref{eq:H2mod}) the quantities
$m, s_{k}, \underline\gs_{k}$ with their continuum analogues
$\widetilde m, \widetilde s_{k}, \widetilde{\underline\gs}_{k}$,
that is, we set
%
%e3.47 ###
%
\begin{equation} \label{eq:H3a}
H^{3}_{t,\gep,\delta}(a,h) := \frac{1}{a} \sum
_{k=1}^{\widetilde m} \widetilde s_{k}
\bigl( \gb_{\widetilde{\underline\gs}_{k}} - \gb_{\widetilde
{\underline\gs}_{k-1}}
+ h (\widetilde{\underline\gs}_{k} - \widetilde{\underline\gs
}_{k-1}) \bigr) .
\end{equation}

It is now convenient to modify slightly the definition (\ref
{eq:H2mod}) of $H^2$.
The laws of the vectors $\Sigma$ and $\widetilde\Sigma$
are mutually absolutely continuous (note that they are probability laws on
the same finite set) and we denote by
$\frac{\dd\Sigma}{\dd\widetilde\Sigma} =
\frac{\dd\Sigma}{\dd\widetilde\Sigma} (\widetilde m;
\widetilde{\underline\gs}_{1}, \ldots, \widetilde{\underline\gs
}_{\widetilde m})$
the corresponding Radon--Nikodym derivative, which does not depend
on $(\widetilde s_{1}, \ldots,\widetilde s_{\widetilde m})$: in fact,
conditionally on
$\widetilde m$, $\widetilde{\underline\gs}_{1}, \ldots, \widetilde
{\underline\gs}_{\widetilde m}$,
the signs $\widetilde s_{1}, \ldots, \widetilde s_{m}$ are i.i.d.
variables that take the values
$\{0,1\}$ with equal probability, and an analogous statement
holds for $s_{1}, \ldots, s_{m}$.
We then redefine
%
%e3.48 ###
%
\begin{equation} \label{eq:H2modbis}
H^{2}_{t,\gep,\delta}(a,h) := H^{3}_{t,\gep,\delta}(a,h) -
\frac{1}{2 a \gl} \log\frac{\dd\Sigma}{\dd\widetilde\Sigma}
.
\end{equation}
Note that this definition of $H^{2}$ yields \textit{the same} $f^{2}$
as (\ref{eq:H2mod}), according to (\ref{eq:fform}).

To prove that $f^{2} \simeq f^{3}$, we can now apply the general scheme
(\ref{eq:fform})--(\ref{eq:ctrl}). Plainly,
%
%e3.49 ###
%
\begin{eqnarray}
&& \min\bigl\{ H^{2}_{t,\gep,\delta}(a,h)
- H^{3}_{t,\gep,\delta}\bigl(a,(1-\rho)h\bigr),\nonumber\\
&&\qquad\hspace*{2.7pt} H^{3}_{t,\gep,\delta}(a,h)
- H^{2}_{t,\gep,\delta}\bigl(a,(1-\rho)h\bigr) \bigr\}
\\
&&\qquad \ge - \frac{1}{2 a \gl} \biggl|
\log\frac{\dd\Sigma}{\dd\widetilde\Sigma} \biggr|
+ \frac{\rho h}{a} \sum_{k=1}^{\widetilde m}
\widetilde s_{k} (\widetilde{\underline\gs}_{k} -
\widetilde{\underline\gs}_{k-1}) ,\nonumber
\end{eqnarray}
therefore, in view of (\ref{eq:ctrl}), we are left with showing that
for all $A, B > 0$ and for $\delta, \gep, a$ small in the usual sense
we have
%
%e3.50 ###
%
\begin{equation} \label{eq:aimstep3}\quad
\limsup_{t\to\infty} \frac1t \log \bE\Biggl[
\exp\Biggl( - A \sum_{k=1}^{\widetilde m}
\widetilde s_{k} (\widetilde{\underline\gs}_{k} - \widetilde
{\underline\gs}_{k-1})
+ B \biggl| \log\frac{\dd\Sigma}{\dd\widetilde\Sigma}
\biggr|
\Biggr) \Biggr] \le 0 .
\end{equation}

We have already observed that,
conditionally on $\widetilde m, \widetilde{\underline\gs}_{1},
\ldots,
\widetilde{\underline\gs}_{\widetilde m}$,
the variables $\widetilde s_{1}, \ldots, \widetilde s_{\widetilde m}$
are i.i.d., taking the values
$\{0, 1\}$ with probability $\frac12$ each, hence
$\frac{\dd\Sigma}{\dd\widetilde\Sigma}$ does not depend on
these variables. Integrating over $\widetilde s_{1}, \ldots,
\widetilde s_{\widetilde m}$, we can
rewrite the expectation in (\ref{eq:aimstep3}) as
%
%e3.51 ###
%
\begin{equation} \label{eq:aimstep3b}
\bE\Biggl[ \Biggl( \prod_{k=1}^{\widetilde m} \biggl( \frac12
+ \frac12 \exp\bigl(-A (\widetilde{\underline\gs}_{k} -
\widetilde{\underline\gs}_{k-1})\bigr) \biggr)
\Biggr)
\exp\biggl( B \biggl| \log\frac{\dd\Sigma}{\dd\widetilde
\Sigma} \biggr|
\biggr) \Biggr] .
\end{equation}
We need some bounds on $\frac{\dd\Sigma}{\dd\widetilde\Sigma}$,
that are given
in the following lemma (whose proof is deferred
to Appendix \ref{sec:crucialk}). Since the result
we are after at this stage is for fixed $\delta>0$, for the sake of simplicity
we are going to fix $\delta=1$: arbitrary values of $\delta$
lead to very similar estimates.
\begin{lemma} \label{th:crucialk}
Fix $\delta= 1$. There exists $\kappa(\gep, a) > 0$ with the
property that
%
%e3.52 ###
%
\begin{equation} \label{eq:crucialk2}
\lim_{\gep\to0} \limsup_{a \to0} \kappa(\gep, a) = 0 ,
\end{equation}
such that, for all values of
$\widetilde m$, $\widetilde{\underline\gs}_{1}, \ldots,
\widetilde{\underline\gs}_{\widetilde m}$, the following bound holds:
%
%e3.53 ###
%
\begin{equation} \label{eq:crucialk1}
\biggl| \log\frac{\dd\Sigma}{\dd\widetilde\Sigma}
(\widetilde m; \widetilde{\underline\gs}_{1}, \ldots,
\widetilde{\underline\gs}_{\widetilde m}) \biggr|
\le \kappa(\gep, a) \sum_{i=1}^{\widetilde m}
\bigl( \log(\widetilde{\underline\gs}_{i} -
\widetilde{\underline\gs}_{i-1}) + 1 \bigr) .
\end{equation}
\end{lemma}

Note that by definition
$(\widetilde{\underline\gs}_{i} - \widetilde{\underline\gs
}_{i-1}) \ge\delta=1$ and
therefore the right-hand side of (\ref{eq:crucialk1}) is positive.
By applying (\ref{eq:crucialk1}), we now see
that the expression in (\ref{eq:aimstep3b})
is bounded above by $\bE[ G_{\widetilde m} ]$, where for $n\in\N$ we set
%
%e3.54 ###
%
\begin{equation} \label{eq:aimstep3c}
G_n := \prod_{i=1}^{n} \frac12 \bigl( 1
+ e^{-A (\widetilde{\underline\gs}_{i} - \widetilde{\underline\gs
}_{i-1})} \bigr)
e^{B\kappa(\gep, a)}
(\widetilde{\underline\gs}_{i} - \widetilde{\underline\gs
}_{i-1})^{B \kappa(\gep, a)} .
\end{equation}
To prove (\ref{eq:aimstep3}), thus completing the proof that
$f^2 \simeq f^3$, it therefore suffices to show that
%
%e3.55 ###
%
\begin{equation} \label{eq:aimstep3d}
\limsup_{t\to\infty}
\frac1t \log \bE[ G_{\widetilde m} ] \le 0 .
\end{equation}

We recall that $\widetilde m = \inf\{k \in\N\dvtx \widetilde
{\underline\gs}_{k} \ge t \}$
and that we had redefined $\widetilde{\underline\gs}_{\widetilde m}
:= t$
for notational convenience.
It is now convenient to switch back to the natural definition (\ref
{eq:tildegs}) of
$\widetilde{\underline\gs}_{\widetilde m}$.
This produces a minor change in $G_{\widetilde m}$, see (\ref
{eq:aimstep3c}): in fact,
only the last factor in the product is modified, and
since $(1+e^{-x}) \le2(1+e^{-y})$ for all $x,y \ge0$,
the \textit{new} $G_{\widetilde m}$ is at most twice the \textit{old} one.
The change is therefore immaterial for the purpose of proving (\ref
{eq:aimstep3d}).

We introduce the filtration $\{\cF_n\}_{n\in\N\cup\{0\}}$,
defined by $\cF_n := \gs(\widetilde{\underline\gs}_0, \ldots,
\widetilde{\underline\gs}_n)$, and we note that $\widetilde m$
is a bounded stopping time for this filtration. Let us set
%
%e3.56 ###
%
\begin{equation} \label{eq:gaga}
\gga= \gga(\gep, a) := \sup_{x \in[-\gep, 0]}
\bE_x \biggl[ \frac12 ( 1 + e^{-A \widetilde{\underline\gs
}_{1}} )
e^{B\kappa(\gep, a)}
(\widetilde{\underline\gs}_{1})^{B \kappa(\gep, a)} \biggr] ,
\end{equation}
where we recall that $\bP_x$ denotes the law of the regenerative set
started at $x$, that is, $\bP_x(\widetilde\tau^\ga\in\cdot) :=
\bP(\widetilde\tau^\ga+ x \in\cdot)$.
From (\ref{eq:aimstep3c}) and the regenerative property of $\widetilde
\tau^\ga$,
we obtain
%
%e3.57 ###
%
\begin{equation} \label{eq:superG}
\bE[ G_{n+1} | \cF_n ] \le \gga G_n .
\end{equation}
If $\gga\le1$, this relation shows that the process $\{G_n\}_{n\ge0}$,
with $G_0 := 1$, is a supermartingale.
Since $\widetilde m$ is a bounded stopping time,
from the optional sampling theorem we deduce
that $\bE[G_{\widetilde m}] \le\bE[G_0] = 1$,
which clearly yields (\ref{eq:aimstep3d}).

It only remains to show that indeed $\gga\le1$, provided
$\gep$ and $a$ are small in the usual sense. Observe that $\widetilde
{\underline\gs}_1$, defined in (\ref{eq:tildegs}),
is a discretized version of
the variable $d_{1-\gep} = d_{1-\gep}(\widetilde\tau^\ga)$,
defined in (\ref{eq:gtds})
(recall that $\delta= 1$): more precisely,
$\widetilde{\underline\gs}_1 = \gep\lceil d_{1-\gep} / \gep\rceil$,
therefore
$d_{1-\gep} \le\widetilde{\underline\gs}_1
\le d_{1-\gep} + \gep$.
Setting $\kappa:= \kappa(\gep, a)$ for short
and applying (\ref{eq:D}), we obtain
%
%e3.58 ###
%
\begin{eqnarray} \label{eq:gaga2}
&& \bE_x \biggl[ \frac12 ( 1 + e^{-A \widetilde{\underline\gs
}_{1}} )
e^{B\kappa}
(\widetilde{\underline\gs}_{1})^{B \kappa} \biggr]\nonumber\\
&&\qquad\le \bE_x \biggl[ \frac12 ( 1 +
e^{-A d_{1-\gep}} ) e^{B\kappa}
(d_{1-\gep} + \gep)^{B \kappa} \biggr] \nonumber\\[-8pt]\\[-8pt]
&&\qquad = \frac{\sin(\pi\ga)}{\pi} \int_{1-\gep}^\infty
\biggl[ \frac12 ( 1 +
e^{-A t} ) e^{B\kappa}
(t + \gep)^{B \kappa} \biggr]\nonumber\\
&&\qquad\quad\hspace*{58.2pt}{}\times \frac{((1-\gep)-x)^\ga}
{(t-(1-\gep))^{\ga} (t-x)} \,\dd t .\nonumber
\end{eqnarray}
Plainly, there exists $\kappa_0 > 0$ such that
the integral in (\ref{eq:gaga2}) is finite for $\kappa\in[0,\kappa_0]$,
for every $x \in[-\gep, 0]$, and it is in fact a \textit{continuous function}
of $(x,\kappa) \in[-\gep, 0] \times[0, \kappa_0]$. Furthermore,
the integral
is strictly smaller than $1$ for $\kappa= 0$ and every $x \in[-\gep, 0]$,
as it is clear from the first line of (\ref{eq:gaga2}).
Therefore, by continuity, there exists $\kappa_1 \in(0, \kappa_0)$ such
that the integral in (\ref{eq:gaga2}) is strictly smaller than one for
$(x,\kappa) \in[-\gep, 0] \times[0,\kappa_1]$.
Looking back at (\ref{eq:gaga}), we see that indeed $\gga\le1$
provided $\kappa(\gep, a) \le\kappa_1$. Thanks to (\ref{eq:crucialk2}),
it suffices to take $\gep$ and $a$ small in the usual sense, and the
proof of
$f^{2} \simeq f^{3}$ is complete.

%%%%%%%%%%%%%%%%%%%%%%%%%%%%%%%%%%%%%%%%%%%%%%%%%%%%%%%%%%%%%%%%%%%%%%%%%%%%%%%%%%
%%%%%%%%%%%%%%%%%%%%%%%%%%%%%%%%%%%%%%%%%%%%%%%%%%%%%%%%%%%%%%%%%%%%%%%%%%%%%%%%%%
%%%%%%%%%%%%%%%%%%%%%%%%%%%%%%%%%%%%%%%%%%%%%%%%%%%%%%%%%%%%%%%%%%%%%%%%%%%%%%%%%%
%s3.4 ###
\subsection{Step 4: Inverse coarse-graining of the regenerative set}
\label{sec:step4}

This step is the close analog of Step 1 (cf. Section \ref{sec:step1})
in the continuum set-up, and a straightforward modification of Step 4
in \cite{cfBdH}.
We will therefore be rather concise.

Recall that the function $f^4$ is nothing but the continuum
finite-volume free energy, cf. (\ref{eq:psipsitilde}),
hence according to (\ref{eq:fform}) it corresponds to the Hamiltonian
[recall (\ref{eq:Zcont}) and (\ref{eq:tildeFt})]
%
%e3.59 ###
%
\begin{eqnarray}
\label{eq:H4}
H^4_{t, \gep, \delta} (a,h) :\!&=& \frac1a
\int_0^t \widetilde\gD(u)
\bigl( \dd\gb(u) + h \,\dd u \bigr) \nonumber\\[-8pt]\\[-8pt]
&=& \frac1a
\sum_{k=1}^{\widetilde m} \int_{\widetilde{\underline\gs
}_{k-1}}^{\widetilde{\underline\gs}_{k}}
\widetilde\gD(u) \bigl( \dd\gb(u) + h \,\dd u \bigr) ,\nonumber
\end{eqnarray}
where we have set $\widetilde\gD(u) := \widetilde\gD^\ga(u)$ for short.
As in the third step, we redefine
$\underline{\widetilde\gs}_{\widetilde m} := t$ for simplicity
[otherwise, the $k=\widetilde m$ term in the sum in (\ref{eq:H4})
would require a separate notation], but we will drop this convention
later.

We now rewrite $H^{3}_{t,\gep,\delta}(a,h)$ by introducing the process
%
%e3.60 ###
%
\begin{equation}
\widehat\gD(u) := \sum_{k=1}^{\widetilde m}
\widetilde s_k \ind_{(\widetilde{\underline\gs}_{k-1}, \widetilde
{\underline\gs}_{k}]} (u) ,
\end{equation}
so that by (\ref{eq:H3a}) we can write
%
%e3.61 ###
%
\begin{equation}
\label{eq:H3a-2}
H^{3}_{t,\gep,\delta}(a,h) =
\frac{1}{a} \sum_{k=1}^{\widetilde m}
\int_{\widetilde{\underline\gs}_{k-1}}^{\widetilde{\underline\gs}_{k}}
\widehat\gD(u) \bigl( \dd\gb(u) + h \,\dd u \bigr) .
\end{equation}
Our aim is to show that $f^3 \simeq f^4$, but we
prove only $f^4 \prec f^3$, since the
argument for the opposite inequality
is very similar. We have [recall (\ref{eq:DH})]
%
%e3.62 ###
%
\begin{eqnarray}
&&
a H^{(4,3)}_{t, \gep, \delta} (a, h , \rho)\nonumber\\
&&\qquad=
\rho h \sum_{k=1}^{\widetilde m} \int_{\widetilde{\underline\gs
}_{k-1}}^{\widetilde{\underline\gs}_{k}}
\widehat\gD(u) \,\dd u \\
&&\qquad\quad{}+
\sum_{k=1}^{\widetilde m} \int_{\widetilde{\underline\gs
}_{k-1}}^{\widetilde{\underline\gs}_{k}}
\bigl( \widetilde\gD(u) - \widehat\gD(u) \bigr) \bigl(
\dd\gb(u) + h \,\dd u
\bigr),\nonumber
\end{eqnarray}
and therefore, arguing as in (\ref{eq:4.10}) and (\ref{eq:ctrl1}),
it is sufficient to show that for every choice of $A$ and $B>0$
%
%e3.63 ###
%
\begin{eqnarray}
\label{eq:eb6.1}
&&\limsup_{t \to\infty}
\frac1t \log\bE
\Biggl[
\exp
\Biggl(
A \int_0^t \vert\widetilde\gD(u) - \widehat\gD(u) \vert
\,\dd u \nonumber\\[-8pt]\\[-8pt]
&&\qquad\hspace*{72pt}{} -
B \sum_{k=1}^{\widetilde m} \widetilde s_k (
\widetilde{\underline\gs}_{k}-\widetilde{\underline\gs}_{k-1}
)
\Biggr)
\Biggr] \le 0 ,\nonumber
\end{eqnarray}
provided $\delta$ and $\gep$ are small in the usual sense. Note that $a$
has disappeared.
%Let us now focus on the complement of the union of the excursions of
%$B^\ga$
%of length at least $\gd$ and denote the intersection of such a set
%with $[0,t]$ by $J_{t, \gd}$.

Let us now focus on the union of the excursions of $\widetilde B^\ga$
whose length is
shorter than $\delta$ and denote the intersection of such a set
with $[0,t]$ by $J_{t, \delta}$.
Then, in analogy with (\ref{eq:hihi}), we have the bound
%
%e3.64 ###
%
\begin{equation}
\label{eq:b7a}
\int_0^t \vert\widetilde\gD(u) - \widehat\gD(u) \vert
\,\dd u \le \vert
J_{t, \delta} \vert+ \widetilde m \gep .
\end{equation}
We now integrate out the $\widetilde s$ variables in (\ref{eq:eb6.1})
[recall that they are i.i.d. $B(1/2)$ variables] and
observe that, since $\widetilde{\underline\gs}_{k}-
\widetilde{\underline\gs}_{k-1}\ge\delta$,
for every $\delta>0$ there exists $\gep_0$ such that for
$\gep\le\gep_0$
%
%e3.65 ###
%
\begin{equation}
A \widetilde m \gep + \frac12
\sum_{k=1}^{\widetilde m} \log\biggl( \frac12 + \frac12 \exp
\bigl(-B (
\widetilde{\underline\gs}_{k}-\widetilde{\underline\gs}_{k-1})
\bigr) \biggr) \le 0 .
\end{equation}
Also notice that, by construction, $|J_{t,\delta} \cap(\underline
{\widetilde\gs}_{k-1},
\underline{\widetilde\gs}_{k}] | \le(\delta+ \gep) \le2 \delta$
for all
$k = 1, \ldots, \widetilde m$, hence $|J_{t,\delta}| \le2 \delta
\widetilde m$.
Therefore, it remains to show that
%
%e3.66 ###
%
\begin{eqnarray}
\label{eq:dep2-1}
&&\limsup_{t \to\infty}
\frac1t \log\bE
\Biggl[
\exp
\Biggl(2 A \delta \widetilde m\nonumber\\[-8pt]\\[-8pt]
&&\qquad\hspace*{69.5pt}{} + \frac12
\sum_{k=1}^{\widetilde m} \log\biggl( \frac12 + \frac12 \exp
\bigl(-B (
\widetilde{\underline\gs}_{k}-\widetilde{\underline\gs}_{k-1})
\bigr) \biggr)\Biggr)
\Biggr] \le 0 .\hspace*{-30pt}\nonumber
\end{eqnarray}

At this point, it is practical to go back to
the original definition of $\underline{\widetilde\gs}_{\widetilde m}$
[cf. (\ref{eq:tildegs})]; this produces a
change in the exponent of (\ref{eq:dep2-1}) which is smaller than
$(\log2)/2$ and this is irrelevant for
the estimate we are after. We then rewrite (\ref{eq:dep2-1}) as
%
%e3.67 ###
%
\begin{eqnarray} \label{eq:poll}
\limsup_{t \to\infty}\frac1t \log
\bE[ G_{\widetilde m} ] \le 0 \nonumber\\[-8pt]\\[-8pt]
\eqntext{\mbox{where }
\displaystyle G_n := \prod_{i=1}^n e^{2 A \delta}
\sqrt{\frac12 \bigl(1 + e^{-B (
\widetilde{\underline\gs}_{k}-\widetilde{\underline\gs}_{k-1}
)} \bigr)} .}
\end{eqnarray}
Let us set
%
%e3.68 ###
%
\begin{equation} \label{eq:gagaga}
\gga= \gga(\delta, \gep) = \sup_{x \in[-\gep, 0]}
e^{2 A \delta} \bE_x \Biggl[ \sqrt{\frac12 (1 + e^{-B
\widetilde{\underline\gs}_{1}} )} \Biggr] ,
\end{equation}
and introduce the filtration $\{\cF_n := \gs(\widetilde{\underline
\gs}_0, \ldots,
\widetilde{\underline\gs}_n)\}_{n\in\N}$.
By the regenerative property of $\widetilde\tau^\ga$, we can write
%
%e3.69 ###
%
\begin{equation} \label{eq:superGbis}
\bE[ G_{n+1} | \cF_n ] \le \gga G_n ,
\end{equation}
therefore if $\gga\le1$ the process $\{G_n\}_{n\ge0}$, with
$G_0 := 0$, is a supermartingale. Since $\widetilde m$
is a bounded stopping time, the optional sampling theorem yields
$\bE[G_{\widetilde m}] \le1$, from which (\ref{eq:poll}) follows.
We are left with showing that $\gamma\le1$ if $\delta$ and
$\gep$ are small in the usual sense.

Recall that $d_s = d_{s}(\widetilde\tau) = \inf\{u > s\dvtx u \in
\widetilde\tau^\ga\}$
[cf. (\ref{eq:gtds})] and observe that, by definition,
$\underline{\widetilde\gs}_1 = j \gep$ if and only if
$d_{\delta- \gep} \in((j-1)\gep, j\gep]$
[cf. (\ref{eq:tildegs})]. Therefore, we may write
$\underline{\widetilde\gs}_1 \ge d_{\delta- \gep}
\ge d_{\delta- \gep+ x}$ for $x \le0$, whence
%
%e3.70 ###
%
\begin{eqnarray} \label{eq:sqsqab}\quad
\bE_x \Bigl[ \sqrt{\tfrac12 (1 +
e^{-B \widetilde{\underline\gs}_{1}} )} \Bigr] &\le&
\bE_x \Bigl[ \sqrt{\tfrac12 (1 + e^{-B
d_{\delta- \gep+ x}} )} \Bigr] \nonumber\\[-8pt]\\[-8pt]
&=&
\bE\Bigl[ \sqrt{\tfrac12 (1 + e^{-B
d_{\delta- \gep}} )} \Bigr] .\nonumber
\end{eqnarray}
Looking back at (\ref{eq:gagaga}), we see that $\gga\le1$
if we show that the right-hand side of (\ref{eq:sqsqab})
is less than $\exp(-2A\delta)$, when $\delta$ and $\gep$ are small
in the usual sense. This condition can be simplified by letting
$\gep\searrow0$: since $d_{\delta- \gep} \to d_\delta$, $\bP$-a.s.,
it suffices to show that
%
%e3.71 ###
%
\begin{equation}
\label{eq:andth}
\bE\Bigl[
\sqrt{ \tfrac12 \bigl( 1+\exp(-B d_{\delta}) \bigr) } \Bigr]
< \exp(-2 A \delta) \qquad
\mbox{for all $\delta> 0$ small enough} .\hspace*{-35pt}
\end{equation}
The law of the variable
$d_{\delta}$ is given in (\ref{eq:D}), hence with a change of
variables we may write
%
%e3.72 ###
%
\begin{eqnarray}
\label{eq:fromDL-RV}
&&\frac1\delta\Biggl( 1-
\bE\Biggl[
\sqrt{ \frac12 \bigl( 1+\exp(-B d_{\delta}) \bigr) } \Biggr]
\Biggr) \nonumber\\[-8pt]\\[-8pt]
&&\qquad=\frac{\sin( \pi\ga) }{\pi} \int_0^\infty
\frac1 \delta\Biggl[ 1-
\sqrt{ \frac12 \bigl( 1+\exp\bigl(-B \delta(1+v)\bigr) \bigr) } \Biggr]
\,\frac{\dd v}{v^\ga(1+v)} .\hspace*{-28pt}\nonumber
\end{eqnarray}
Since the term between square brackets in the right-hand side is
positive and
asymptotically equivalent, as $\delta\searrow0$, to $\delta B(1+v)/4$,
Fatou's lemma guarantees that the limit as $\delta\searrow0$ of the expression
in (\ref{eq:fromDL-RV}) is equal to $+\infty$
and this entails that (\ref{eq:andth})
holds.

This concludes the proof of Step 4 and, hence, the proof of
Theorem \ref{th:main-tech}.

\begin{appendix}
%s4 ###
\section{\texorpdfstring{Completing the proof of Proposition 2.2%\protect\ref{th:Kin}
}{Completing the proof of Proposition 2.2}}
\label{sec:contbound}

In this section, we are going to prove (\ref{eq:cru}), that is, for every
$\eta\in(0,\infty)$ there exists $D(\eta) \in(0,\infty)$ such that
%
%e4.1 ###
%
\begin{equation} \label{eq:cruapp}
\bE[ \bbE[ \exp( \eta \Theta_T (\gb,
\widetilde\gD^\ga) ) ] ]
\le D(\eta) e^{ D(\eta) T} \qquad \mbox{for every }
T > 0 .
\end{equation}
We first state some important estimates concerning the regenerative set
$\widetilde\tau^\ga$.

%%%%%%%%%%%%%%%%%%%%%%%%%%%%%%%%%%%%%%%%%%%%%%%%%%%%%%%%%%%

%s4.1 ###
\subsection{Regenerative set, excursions and local time}
\label{sec:regapp}

We recall the basic link between regenerative set and subordinators.
Let $(\gs= \{\gs_t\}_{t\ge0}, \bP)$ denote the stable subordinator
of index $\ga$, that is, the L\'evy process with zero drift,
zero Brownian component and with L\'evy measure given by\vspace*{2pt}
$\Pi(\dd x) := \frac{C}{x^{1+\ga}} \ind_{(0,\infty)}(x) \,\dd x$
with $C>0$. We choose as usual a right-continuous version of $\gs$.
The value of the constant $C$ is quite immaterial (it corresponds to
rescaling time or space by a constant factor) and a useful normalization
it to fix $C$ so that $\int_0^\infty(1-e^{-x}) \Pi(\dd x) = 1$.
In this way, the \textit{L\'evy exponent} of $\sigma$, defined by
$\Phi(\lambda) := - \log\bE[e^{-\gl\sigma_1}]
= \int_0^\infty(1-e^{-\lambda x}) \Pi(\dd x)$,
equals exactly $\lambda^{\ga}$ for all $\gl\ge0$.

If we denote by $\Delta\gs_t := \gs_{t+}-\gs_t$ the size
of the jump of $\gs$ at epoch $t$,
it is well known that $\gs_t = \sum_{s \in(0,t]} \Delta\gs_s$,
that is,
$\gs$ increases only by jumps.
A remarkable property of $\gs$ is its scale invariance: $\{\gs_{ct}\}
_{t\ge0}$ has
the same law as $\{c^{1/\ga}\gs_{t}\}_{t\ge0}$.
We also recall some basic estimates (cf. Theorems 8.2.1 and 8.2.2
in \cite{cfBinGolTeu}):
%
%e4.2 ###
%
\begin{eqnarray} \label{eq:aszero}
\bP(\gs_1 > x) &=& \frac{(\mbox{const.})}{x^\ga} \bigl( 1 + o(1)
\bigr)\qquad
\mbox{as } x \to+\infty, \nonumber\\[-8pt]\\[-8pt]
\bP(\gs_1 < x) &=& \exp\biggl(- \frac{(\mbox{const.}')}{x^{\ga/(1-\ga)}}
\bigl( 1 + o(1) \bigr) \biggr) \qquad
\mbox{as } x \searrow0 .\nonumber
\end{eqnarray}

If we set $\cE:= [0,\infty) \times(0,\infty)$, the random set of points
$\{(t,\Delta\gs_t)\}_{t\in[0,\infty)} \cap\cE$ (note that we only keep
the positive jumps $\Delta\gs_t > 0$) is a \textit{Poisson random
measure} (sometimes simply
called \textit{Poisson process}) on $\cE$ with intensity measure
$\dd t \otimes\Pi(\dd x)$, where of course $\dd t$ denotes the
Lebesgue measure.
The stochastic process $\{\Delta\gs_t\}_{t\in[0,\infty)}$ is called a
\textit{Poisson point process} on $(0,\infty)$ with intensity measure
$\Pi$.

The basic link with regenerative sets
is as follows: the random closed set of $[0,\infty)$
defined as the closure of the image of the process $\gs$,
that is, $\overline{\{\gs_t\}_{t\ge0}}$, is precisely
the $\ga$-stable regenerative set $\widetilde\tau^\ga$ we are
considering. Therefore,
the set of jumps $\{\Delta\gs_t\}_{t \ge0}$ coincides with the set of
widths $\{|I_n|\}_{n\in\N} \cup\{0\}$ of the excursions of
$\widetilde\tau^\ga$.

Let us discuss an application of these results that will be useful later.
If we denote by $L_t := \inf\{u \ge0\dvtx \gs_u > t\}$ the inverse of
$\gs$,
known as the \textit{local time} of $\widetilde\tau^\ga$, we may write
%
%e4.3 ###
%
\begin{eqnarray}\qquad\quad
\sum_{n\in\N\dvtx I_n \subseteq(0,2)}
|I_n|^{1-\gep} &=& \sum_{t \in(0, L_2)}
(\Delta\gs_t)^{1-\gep}= \sum_{t \in(0, L_2)} f(\Delta\gs_t)
\nonumber\\[-8pt]\\[-8pt]
\eqntext{\mbox{where } f(x) := x^{1-\gep} \ind_{[0,2]}(x) ,}
\end{eqnarray}
therefore for $\gl> 0$ we have by Cauchy--Schwarz
\begin{eqnarray*}
&&
\bE\biggl[ \exp\biggl( \gl \sum_{n\in\N\dvtx I_n \subseteq
(0,2)} |I_n|^{1-\gep}\biggr) \biggr] \\
&&\qquad\le \sum_{m \in\N}
\bE\biggl[ \exp\biggl( \gl \sum_{t \in(0, m)} f(\Delta\gs_t)
\biggr) \ind_{\{m-1 < L_2 \le m \}} \biggr] \\
&&\qquad\le \sum_{m \in\N}
\sqrt{ \bE\biggl[ \exp\biggl( 2 \gl \sum_{t \in(0, m)}
f(\Delta\gs_t)
\biggr) \biggr] \bP[ m-1 < L_2 \le m ] } .
\end{eqnarray*}
By the definition of $L$, the scale invariance of $\gs$ and (\ref{eq:aszero}),
we have for some $c>0$
%
%e4.4 ###
%
\begin{eqnarray}
\bP[ m-1 < L_2 \le m ] % \le\bP[ L_2 > m-1 ]
&\le& \bP[ \gs_{m-1} < 2 ] =
\bP\biggl[ \gs_{1} < \frac{2}{(m-1)^{1/\ga}} \biggr]\nonumber\\[-8pt]\\[-8pt]
&\le& e^{- c (m-1)^{1/(1-\ga)}} .\nonumber
\end{eqnarray}
By Campbell's formula for Poisson processes (cf.
equation (3.17) in \cite{cfKin2}) we obtain
%
%e4.5 ###
%
\begin{eqnarray} \label{eq:Cgl}
&&\bE\biggl[ \exp\biggl( 2 \gl \sum_{t \in(0, m)} f(\Delta\gs_t)
\biggr) \biggr] \nonumber\\
&&\qquad= \exp
\biggl( m \int_0^\infty\bigl(e^{2 \gl f(x)} - 1\bigr) \Pi(\dd x)
\biggr)
= e^{C(\gl) m} , \\
\eqntext{\mbox{where }
\displaystyle C(\gl) := \int_0^2 \frac{e^{2 \gl x^{1-\gep}} -
1}{x^{1+\ga}} \,\dd x
< \infty \mbox{ for } 0 < \gep< 1-\ga.}
\end{eqnarray}
From the last relations we then obtain, for some $c_1 \in(0,\infty)$,
%
%e4.6 ###
%
\begin{eqnarray}
\bE\biggl[ \exp\biggl( \gl \sum_{n\in\N\dvtx I_n \subseteq(0,2)}
|I_n|^{1-\gep}
\biggr) \biggr] &\le& \sum_{m \in\N} e^{1/2 (C(\gl) m
- c (m-1)^{1/(1-\ga)})}\nonumber\\[-8pt]\\[-8pt]
&\le& c_1 e^{c_1 (C(\gl))^{1/\ga}} ,\nonumber
\end{eqnarray}
where the last inequality can be checked, for example, by approximating
the sum
with an integral and developing the function
$e^{1/2 [C(\gl) x - c x^{1/(1-\ga)}]}$ around its maximum.

Since $e^{2 \gl y} - 1 \le2 \gl e^{4\gl} y$ for $y \in[0,2]$,
it follows from (\ref{eq:Cgl}) that $C(\gl) \le(\mbox{const.}) e^{5
\gl}$.
By Markov's inequality, we then obtain
%
%e4.7 ###
%
\begin{eqnarray}
\bP\biggl[ \sum_{n\in\N\dvtx I_n \subseteq(0,2)} |I_n|^{1-\gep}
> x \biggr]
&\le& c_1 e^{c_1 (C(\gl))^{1/\ga} - \gl x}\nonumber\\[-8pt]\\[-8pt]
&\le& c_1 e^{c_2 e^{5\gl/\ga} - \gl x} ,\nonumber
\end{eqnarray}
for some $c_2 \in(0,\infty)$. Optimizing over $\gl$ yields, for
every $x > 0$,
%
%e4.8 ###
%
\begin{eqnarray}
\bP\biggl[ \sum_{n\in\N\dvtx I_n \subseteq(0,2)} |I_n|^{1-\gep}
> x \biggr]
&\le& \min\bigl\{ c_1 e^{-(\alpha x/5)
[\log((\alpha x)/(5 c_2)) - 1]} ,
1 \bigr\} \nonumber\\[-8pt]\\[-8pt]
&\le& c_3 e^{- c_3 x} ,\nonumber
\end{eqnarray}
for a suitable $c_3 \in(0,\infty)$. We can finally estimate the
quantity we are interested in:
\begin{eqnarray*}
&& \bE\Biggl[ \exp\Biggl( \gamma \sqrt{T}
\sqrt{\sum_{n\in\N\dvtx I_n \subseteq(0,2)} |I_n|^{1-\gep}} \Biggr)
\Biggr] \\
&&\qquad = \int_0^\infty\bP\Biggl[ \exp\Biggl( \gamma \sqrt{T}
\sqrt{\sum_{n\in\N\dvtx I_n \subseteq(0,2)} |I_n|^{1-\gep}} \Biggr)
> t \Biggr]\,
\dd t \\
&&\qquad = \int_0^\infty\bP\biggl[
\sum_{n\in\N\dvtx I_n \subseteq(0,2)} |I_n|^{1-\gep}
> \frac{(\log t)^2}{\gamma^2 T}
\biggr] \,\dd t \\
&&\qquad \le c_3 \int_0^\infty e^{- c_3 (\log t)^2/(\gamma^2 T)}
\,\dd t\\
&&\qquad= c_3 \int_{-\infty}^\infty e^x e^{- c_3 x^2/(\gamma^2 T)}
\,\dd x\\
&&\qquad\le c_4 \gamma \sqrt{T} e^{c_4 \gamma^2 T} ,
\end{eqnarray*}
for some $c_4 \in(0,\infty)$, by a Gaussian integration. We have thus
proven that,
if $\gep< 1-\ga$, there exists $c_4 \in(0,\infty)$ such that for
all $\gamma,T > 0$
%
%e4.9 ###
%
\begin{equation} \label{eq:bleach}
\bE\Biggl[ \exp\Biggl( \gamma \sqrt{T}
\sqrt{\sum_{n\in\N\dvtx I_n \subseteq(0,2)} |I_n|^{1-\gep}} \Biggr)
\Biggr]
\le c_4 \gamma \sqrt{T} e^{c_4 \gamma^2 T} .
\end{equation}

%%%%%%%%%%%%%%%%%%%%%%%%%%%%%%%%%%%%%%%%%%%%%%%%%%%%%%%%%%%
%s4.2 ###
\subsection{\texorpdfstring{Proof of (\protect\ref{eq:cruapp})}{Proof of (A.1)}}

We recall that
%
%e4.10 ###
%
\begin{equation}\label{eq:Thetaapp}
\Theta_T (\gb, \widetilde\gD^\ga)
:= {\sup_{-1 \le x \le T , 0 \le y \le T + 1}}
| \cH_{0,y;\theta_x \gb}(\widetilde\gD^\ga) | .
\end{equation}
Recalling (\ref{eq:Zcont}), we can write
\[
\cH_{0,y; \theta_x \gb}(\widetilde\gD^\ga) =
-2 \gl\int_{0}^{y} \widetilde\Delta^\ga(u) \,\dd(\theta_x \gb
)(u) -
2 \gl h \int_0^{y} \widetilde\Delta^\ga(u) \,\dd u ,
\]
and note that the second term is bounded in absolute
value by $2 \gl h y$. For the purpose of proving
(\ref{eq:cruapp}), we may therefore focus on the
first term: we set
%
%e4.12 ###
%e4.11 ###
%
\begin{eqnarray} \label{eq:gamma}
\gamma_{x,y}(\gb, \widetilde\gD^\ga)  :\!& = &
\int_{0}^{y} \widetilde\Delta^\ga(u) \,\dd(\theta_x
\gb)(u)\nonumber\\[-8pt]\\[-8pt]
&=&
\int_{x}^{x+y} \widetilde\Delta^\ga(u-x) \,\dd\gb(u) ,\nonumber \\
\label{eq:Gamma}
\Gamma_T (\gb, \widetilde\gD^\ga)
:\!& = & \sup_{(x,y) \in\cS_{T}}
\gamma_{x,y}(\gb, \widetilde\gD^\ga) \nonumber\\[-8pt]\\[-8pt]
\eqntext{\mbox{where }
\cS_T := [-1,T] \times[0,T+1] .}
\end{eqnarray}
We stress that $\Gamma_T$ is defined as the supremum of $\gamma_{x,y}$,
not of $|\gamma_{x,y}|$.
Notice however that, for fixed $\widetilde\gD^\ga$, the process
$\gamma= \{\gamma_{x,y}(\gb,\widetilde\gD^\ga)\}_{x,y}$ under
$\bbP$ is Gaussian and centered,
in particular it has the same law as $-\gamma$.
Since $e^{|x|} \le e^x + e^{-x}$, we may then write
%
%e4.13 ###
%
\begin{equation}\qquad
\bE[ \bbE[ \exp( \eta\Theta_T (\gb, \widetilde
\gD^\ga) ) ] ]
\le 2 e^{2 \gl h (T+1)}
\bE[ \bbE[ \exp( 2 \eta\gl \Gamma_T (\gb,
\widetilde\gD^\ga) ) ] ] .
\end{equation}
Looking back at (\ref{eq:cruapp}), we are left with showing that, for
every $\eta> 0$,
there exists (a possibly different) $D(\eta) \in(0,\infty)$ such that
%
%e4.14 ###
%
\begin{equation} \label{eq:ttoo}
\bE[ \bbE[ \exp( \eta \Gamma_T (\gb, \widetilde
\gD^\ga) ) ] ]
\le D(\eta) e^{D (\eta) T} \qquad \forall T > 0 .
\end{equation}

Let us set $\Gamma_T := \Gamma_T (\gb, \widetilde\gD^\ga)$ for short.
It is convenient to split
%
%e4.15 ###
%
\begin{equation} \label{eq:mysplit}
\bE[ \bbE[ \exp( \eta \Gamma_T ) ] ]
= \bE\bigl[ \exp( \eta \bbE[\Gamma_T ] )
\cdot\bbE\bigl[ \exp\bigl( \eta (\Gamma_T - \bbE[\Gamma_T])
\bigr) \bigr] \bigr] .
\end{equation}
To prove (\ref{eq:ttoo}), we use the powerful tools of the theory of continuity
of Gaussian processes. Let us introduce (for a fixed realization of
$\widetilde\gD^\ga$)
the \textit{canonical metric} associated to the gaussian process $\gamma
$, defined
for $(x,y), (x',y') \in\cS_T = [-1,T] \times[0,T+1]$
by
%
%e4.16 ###
%
\begin{equation} \label{eq:cano}
d((x,y), (x', y')) :=
\sqrt{ \bbE\bigl[ \bigl( \gamma_{x',y'}(\gb,\widetilde\gD^\ga) -
\gamma_{x,y}(\gb,\widetilde\gD^\ga) \bigr)^2 \bigr] } .
\end{equation}
For $\gep> 0$ we define $N_T(\gep) = N_{T,\widetilde\gD^\ga}(\gep
)$ as
the least number of open balls of radius $\gep$ (in the canonical
metric) needed to cover the parameter space
$\cS_T$. The quantity $\log N_T(\gep)$
is called the \textit{metric entropy} of $\gamma$.
It is known (\cite{cfAdl}, Corollary 4.15) that
the finiteness of $\int_0^\infty\sqrt{\log N_T(\gep)} \,\dd\gep$
ensures the existence of a version of the process $\gamma$ which
is continuous in the parameter space.
Moreover, there exists a universal constant $K \in(0,\infty)$ such that
%
%e4.17 ###
%
\begin{equation} \label{eq:mmean}
\bbE[\Gamma_T(\gb, \widetilde\gD^\ga)] \le K \int
_0^\infty
\sqrt{\log N_{T,\widetilde\gD^\ga}(\gep)} \,\dd\gep.
\end{equation}
We show below that, for $\bP$-a.e. realization of $\widetilde\gD^\ga$,
indeed $\int_0^\infty\sqrt{\log N_{T, \widetilde\gD^\ga}(\gep)}%\,
\*\dd\gep< \infty$,
so we may (and will) choose henceforth a continuous
version of the process~$\gamma$.%\looseness=1

To estimate\vspace*{1pt} the right-hand side of (\ref{eq:mysplit}),
let us denote by $\gs_T^2 = \gs_{T,\widetilde\gD^\ga}^2$ the
maximal variance
of the process $\gamma$, that is, $\gs_T^2 :=
\sup_{(x,y) \in\cS_T} \bbE[\gamma_{x,y}(\gb, \widetilde\gD^\ga)^2]$.
Since $\gamma$ is continuous, it follows easily
by Borell's inequality (\cite{cfAdl}, Theorem 2.1) that
% \bbP[ | \Gamma_T - \bbE[\Gamma_T] |
% > t ] \le2 e^{- \frac12 t^2 / \gs_T^2} ,
%therefore
% & \bbE[ \exp( \eta(\Gamma_T - \bbE[\Gamma_T]) ) ]
% = \int_0^\infty\bbP[ \exp( \eta(\Gamma_T - \bbE[
% > t ] \dd t \\
% & = \int_0^\infty\bbP[ \Gamma_T - \bbE[\Gamma_T] >
% \frac{\log t}{\eta} ] \dd t \le1 +
% 2 \int_1^\infty\exp( - \frac12 \frac{(\log t)^2}
% {\eta^2 \gs_T^2} ) \dd t \\
% & = 1 + 2 \int_0^\infty
% \exp( - \frac12 \frac{y^2}
% {\eta^2 \gs_T^2} ) \exp(y) \dd y =
% 1 + \sqrt{2\pi} \gs_T \exp( \frac12 \eta^2
%
\[
\bbE\bigl[ \exp\bigl( \eta(\Gamma_T - \bbE[\Gamma_T]) \bigr) \bigr]
\le C' \gs_T \exp\bigl( \tfrac12 \eta^2 \gs_T^2
\bigr) ,
\]
where $C' \in(0,\infty)$ is an absolute constant.
Now observe that $\gs_T^2$ is uniformly bounded: by (\ref{eq:gamma})
and the isometry property of the Wiener\vadjust{\goodbreak} integral, since $|\widetilde
\gD^\ga(\cdot)|\le1$, we can write
%
%e4.18 ###
%
\begin{eqnarray} \label{eq:gs}
\gs_T^2 :\!&=& \sup_{(x,y) \in\cS_T}
\bbE[ \gamma_{x,y}(\gb,\widetilde\gD^\ga)^2 ] \nonumber\\[-8pt]\\[-8pt]
&=&
\sup_{(x,y) \in\cS_T}
\int_{x}^{x+y} \widetilde\gD^\ga(u-x)^2 \,\dd u \le T+1 .\nonumber
\end{eqnarray}
Looking back at (\ref{eq:mysplit}) and recalling (\ref{eq:mmean}), we
have proven that
there exists $C \in(0,\infty)$ such that
%
%e4.19 ###
%
\begin{eqnarray} \label{eq:appartial}
&&\bE[ \bbE[ \exp( \eta \Gamma_T(\gb, \widetilde
\gD^\ga) ) ] ]\nonumber\\[-8pt]\\[-8pt]
&&\qquad\le C e^{C \eta^2 T} \bE\biggl[
\exp\biggl( K \eta\int_0^\infty\sqrt{\log N_{T, \widetilde\gD
^\ga}(\gep)} \,\dd\gep
\biggr) \biggr].\nonumber
\end{eqnarray}

To complete the proof of (\ref{eq:ttoo}),
it remains to estimate $N_{T, \widetilde\gD^\ga}(\gep)$, which
requires some effort.
For a fixed realization of $\widetilde\gD^\ga$,
we introduce the function $\rho_T\dvtx \R^+ \to\R^+$ defined by
%
%e4.20 ###
%
\begin{equation} \label{eq:rho}
\rho_T(\delta) := \sup_{(x,y), (x',y') \in\cS_T\dvtx |(x,y) -
(x', y')| \le\delta}
d ( (x,y), (x', y') ) ,
\end{equation}
where $|(x,y) - (x', y')|^2 := (x-x')^2 + (y-y')^2$ denotes the
Euclidean norm
and we recall that the canonical metric $d$ is defined in (\ref{eq:cano}).
Note that $\rho_T(\cdot)$ is a nondecreasing function
which is eventually constant: $\rho_T(\delta) = \rho_T(\sqrt
{2}(T+1))$ for
every $\delta\ge\sqrt{2}(T+1)$, simply because $\sqrt{2}(T+1)$ is the
diameter of the space $\cS_{T} = [-1,T] \times[0,T+1]$.

Plainly, for every fixed $\delta>0$, we can cover the square $\cS
_{T}$ with
no more than $(\frac{T+1}{\delta} + 1)^{2}$ open squares of side
$\delta$.
Since the \textit{Euclidean} distance
between a point in a square of side $\delta$ and the center of the
square is at most $\delta/\sqrt{2}$, the corresponding
distance in the canonical metric is at most $\rho_{T}(\delta/\sqrt{2})$,
by the definition of $\rho_{T}$.
Therefore, a square of side $\delta$ can be covered with a ball (in
the canonical metric) of radius
$\rho_{T}(\delta/\sqrt{2})$ centered at the center of the square.
If we set $\gep:= \rho_{T}(\delta/\sqrt{2})$, this means that we
need at most
$(\frac{T+1}{\delta} + 1)^{2}$ balls (in the canonical metric) of
radius $\gep$
to cover the whole parameter space $\cS_{T}$. Put otherwise,
we have shown that for every $\gep> 0$,
%
%e4.21 ###
%
\begin{equation}
N_T(\gep) \le \biggl(1 + \frac{T+1}{\sqrt{2} \rho_T^{-1}(\gep
)} \biggr)^2 ,
\end{equation}
where $\rho_T^{-1}$ is well defined because $\rho_{T}$ is nondecreasing
and \textit{continuous}, as it will be clear below.
Since $N_T(\gep) = 1$ for $\gep> \rho_T((T+1)/\sqrt{2})$
(we can cover $\cS_{T}$ with just one ball), we obtain the estimate
\[
\int_0^\infty\sqrt{\log N_T(\gep)} \,\dd\gep\le
\int_0^{\rho_T((T+1)/\sqrt{2})}
\sqrt{2 \log\biggl( 1 + \frac{T+1}{\sqrt{2} \rho_T^{-1}(\gep
)} \biggr)} \,\dd\gep.
\]
By a change of variables and integrating by parts, we obtain
%
%e4.22 ###
%
\begin{eqnarray} \label{eq:long}\qquad\quad
&&\int_0^\infty\sqrt{\log N_T(\gep)} \,\dd\gep\nonumber\\
&&\qquad\le
\int_0^{({T+1})/{\sqrt{2}}}
\sqrt{2 \log\biggl(1 + \frac{T+1}{\sqrt{2} t} \biggr)} \,\dd\rho
_T(t)\nonumber\\
&&\qquad= \sqrt{2 \log2}
\rho_T \biggl( \frac{T+1}{\sqrt{2}} \biggr)\\
&&\qquad\quad{} + \int_0^{
({T+1})/{\sqrt{2}}}
\frac{\rho_T(t)}{t \sqrt{2 \log(1 + ({T+1})/({\sqrt{2}
t}) )}}
\frac{T+1}{T+1+ \sqrt{2} t} \,\dd t \nonumber\\
&&\qquad\le \sqrt{2} \rho_T \biggl( \frac{T+1}{\sqrt{2}} \biggr)
+
\int_0^{({T+1})/{\sqrt{2}}} \frac{\rho_T(t)}{t \sqrt{2 \log
(1 + ({T+1})/({\sqrt{2} t} ))}}\,
\dd t ,\nonumber
\end{eqnarray}
where in the integration by parts we have used the fact that,
for $\bP$-a.e. realization of $\widetilde\gD^\ga$, we have
$\sqrt{2 \log(1 + \frac{T+1}{\sqrt{2} t})} \rho_{T}(t) \to0$
as $t \to0$,
as we prove below.

To proceed with the estimates, we need to obtain bounds on
$\rho_{T}$, hence, we start
from the definition (\ref{eq:gamma}) of $\gamma_{x,y}(\gb,
\widetilde\gD^{\ga})$.
By the properties the Wiener integral, we can write
\begin{eqnarray*}
&& d ((x,y), (x',y') )^2 \\
&&\qquad=
\bbE\bigl[ \bigl( \gamma_{x',y'}(\gb,\widetilde\gD^\ga)
- \gamma_{x,y}(\gb,\widetilde\gD^\ga) \bigr)^2 \bigr] \\
&&\qquad = \int_{-1}^{2T+1} \bigl( \widetilde\gD^\ga(u-x') \ind
_{[x',x'+y']}(u) -
\widetilde\gD^\ga(u-x) \ind_{[x,x+y]}(u) \bigr)^2 \,\dd u \\
&&\qquad = \int_{-1}^{2T+1} \bigl| \widetilde\gD^\ga(u-x') \ind
_{[x',x'+y']}(u) -
\widetilde\gD^\ga(u-x) \ind_{[x,x+y]}(u) \bigr| \,\dd u ,
\end{eqnarray*}
where the last equality holds simply because $\widetilde\gD^\ga
(\cdot)$ takes values
in $\{0,1\}$. Incidentally, this expression shows that the canonical metric
$d(\cdot, \cdot)$ is continuous on $\cS_{T}$ (because the translation
operator is continuous in $L^{1}$). Therefore, $\rho_{T}(\cdot)$ is a
continuous
function, as we stated before.

By the triangle inequality, we get for $x' \le x$
\begin{eqnarray*}
d ((x,y), (x',y') )^2 &\le&
\int_{-1}^{2T+1} \widetilde\gD^\ga(u-x')
\bigl| \ind_{[x',x'+y']}(u) - \ind_{[x,x+y]}(u) \bigr| \,\dd u \\
&&{}
+ \int_{-1}^{2T+1} | \widetilde\gD^\ga(u-x') -
\widetilde\gD^\ga(u-x) | \ind_{[x,x+y]}(u) \,\dd u \\
&\le& |x'-x| + |(x'+y') - (x+y)| \\
&&{} +
\int_{x}^{2T+1} | \widetilde\gD^\ga(u-x') -
\widetilde\gD^\ga(u-x) | \,\dd u .
\end{eqnarray*}
Recall that
$\widetilde\gD^\ga(s) = \sum_{n\in\N} \widetilde\xi_n \ind
_{I_n}(s)$, where
$\{I_n\}_{n\in\N}$ are the connected components of the open set
$(\widetilde\tau^\ga)^\complement$
and $\{\widetilde\xi_n\}_{n\in\N}$ are i.i.d. Bernoulli variables
of parameter $1/2$.
For every finite interval $I$, we have
the bound $\int_{\R} |\ind_{I}(u-x') - \ind_{I}(u-x)| \,\dd u
\le2 \min\{|I_n|, |x'-x|\}$, whence
%
%e4.23 ###
%
\begin{eqnarray}
&&\int_{x}^{2T+1} | \widetilde\gD^\ga(u-x') -
\widetilde\gD^\ga(u-x) | \,\dd u \nonumber\\[-8pt]\\[-8pt]
&&\qquad\le
\sum_{n\in\N\dvtx I_n \cap(0,2(T+1)) \ne\varnothing} \min\{|I_n|,
\delta\}.\nonumber
\end{eqnarray}
Therefore, recalling definition (\ref{eq:rho}), we can write
%
%e4.24 ###
%
\begin{equation} \label{eq:rho2}
\rho_T(\delta)^2 \le 3 \delta+
\sum_{n\in\N\dvtx I_n \cap(0,2(T+1)) \ne\varnothing} \min\{|I_n|,
\delta\} .
\end{equation}
Observe that the sum in the right-hand side can be rewritten as $\delta
N_{\delta} +
A_{\delta}$, where $N_{\delta}$ is the \textit{number} of excursions $I_{n}$
that intersect $(0,2(T+1))$ with $|I_{n}| > \delta$
and $A_{\delta}$ is the \textit{total area} covered by the excursions $I_{n}$
that intersect $(0,2(T+1))$ with $|I_{n}| \le\delta$. The asymptotic behavior
as $\delta\searrow0$ of $N_{\delta}$ and $A_{\delta}$ is as
follows: there exists a positive
constant $c = c(\ga)$ such that
%
%e4.25 ###
%
\begin{equation} \label{eq:asNA}
\lim_{\delta\searrow0} \delta^{\ga} N_{\delta} =
\lim_{\delta\searrow0} \frac{A_{\delta}}{\delta^{1-\ga}}
= c L_{2(T+1)} ,\qquad
\mbox{$\bP$-a.s.} ,
\end{equation}
where $\{L_{t}\}_{t\ge0}$ is the local time associated to the regenerative
set $\widetilde\tau^\ga$ (whose definition is recalled in
Appendix \ref{sec:regapp}).
The relations in (\ref{eq:asNA})
are proven in \cite{cfRevYor} [cf.~Proposition XII-(2.9)
and Exercise XII-(2.14)] in the special case $\ga= \frac12$, but the
proof is easily extended to the general case.
Looking back at (\ref{eq:rho2}), it follows that, for $\bP$-a.e. realization
of $\widetilde\gD^{\ga}$, we have $\rho_{T}(\delta) \sim\sqrt{2
c} \sqrt{L_{2(T+1)}} \delta^{(1-\ga)/2}$
as $\delta\searrow0$. In particular,
$\sqrt{\log(1 + \frac{T+1}{\sqrt{2} t})} \rho_{T}(t) \to0$
as $t \searrow0$,
a property used in the integration by parts in (\ref{eq:long}).

We are ready to bound the terms in the last line of (\ref{eq:long}).
Note that the first term is easily controlled:
by definition $d((x,y), (x',y')) \le2 \gs_T$,
hence it follows by (\ref{eq:gs}) that
%
%e4.26 ###
%
\begin{equation} \label{eq:firterm}
\sqrt{2} \rho_T \biggl( \frac{T+1}{\sqrt{2}} \biggr)
\le 2 \sqrt{2} \sqrt{T+1} .
\end{equation}
Now observe that from (\ref{eq:rho2}), we have
%
%e4.27 ###
%
\begin{eqnarray}
\rho_{T}(\delta) \le F_{T+1}(\delta) \hspace*{70pt}\nonumber\\[-8pt]\\[-8pt]
\eqntext{\mbox{where }
\displaystyle F_{M}(\delta) := \sqrt{ 3\delta+
\sum_{n\in\N\dvtx I_n \cap(0,2M) \ne\varnothing} \min\{|I_n|, \delta
\}} .}
\end{eqnarray}
By the scale invariance of the regenerative set $\widetilde\tau^\ga$
it follows that, under $\bP$, $\{F_{M}(t)\}_{t\ge0}$ has the same law as
$\{\sqrt{M} F_{1}(\frac{t}{M})\}_{t\ge0}$. Therefore, we can
bound the second term in the last line of (\ref{eq:long}) as follows:
%
%e4.28 ###
%
\begin{eqnarray} \label{eq:secterm}
&&\int_0^{({T+1})/{\sqrt{2}}} \frac{\rho_T(t)}{t \sqrt{2 \log
(1 + ({T+1})/({\sqrt{2} t}) )}}\,
\dd t \nonumber\\
&&\qquad\le
\int_0^{({T+1})/{\sqrt{2}}} \frac{F_{T+1}(t)}
{t \sqrt{2 \log(1 + ({T+1})/({\sqrt{2} t}) )}}\,
\dd t\\
&&\qquad\stackrel{d}{=} \sqrt{T+1} \cM,\nonumber
\end{eqnarray}
where, performing the change of variable $t = (T+1)s$ in the integral,
we have introduced the variable $\cM$ defined by
%
%e4.29 ###
%
\begin{eqnarray}
\cM:\!&=& \int_0^{{1}/{\sqrt{2}}} \frac{F_{1}(s)}
{s \sqrt{2 \log(1 + {1}/({\sqrt{2} s}) )}} \,\dd s\nonumber\\[-8pt]\\[-8pt]
&=& \int_0^{{1}/{\sqrt{2}}} \frac{1}{s}
\sqrt{\frac{3s + \sum_{n\in\N\dvtx I_n \cap(0,2) \ne\varnothing}
\min\{|I_n|, s\}}
{2 \log(1 + {1}/({\sqrt{2} s}) )}} \,\dd s .\nonumber
\end{eqnarray}
We can finally come back to (\ref{eq:appartial}): applying
(\ref{eq:long}), (\ref{eq:firterm}) and (\ref{eq:secterm}) we obtain
%
%e4.30 ###
%
\begin{equation} \label{eq:appartial2}
\bE\biggl[
\exp\biggl( K \eta\int_0^\infty\sqrt{\log N_{T, \widetilde\gD
^\ga}(\gep)} \,\dd\gep
\biggr) \biggr] \le \bE\bigl[ e^{K \eta \sqrt{T+1}
\bigl(2 \sqrt{2} + \cM\bigr)} \bigr] .
\end{equation}

It only remains to estimate the law of $\cM$.
Let us fix an arbitrary $\gep\in(0, 1-\ga)$:
applying the Cauchy--Schwarz inequality, we obtain
%
%e4.31 ###
%
\begin{eqnarray}
\cM &\le& \sqrt{\int_0^{{1}/{\sqrt{2}}}
\frac{1}{2 s^{1-\gep} \log(1 + {1}/({\sqrt{2}
s}) )} \,\dd s}
\nonumber\\[-8pt]\\[-8pt]
&&{}\times \sqrt{\int_0^{{1}/{\sqrt{2}}} \biggl(\frac
{3}{s^\gep} +
\sum_{n\in\N\dvtx I_n \cap(0,2) \ne\varnothing}
\frac{\min\{|I_n|, s\}}{s^{1+\gep}} \biggr) \,\dd s} .\nonumber
\end{eqnarray}
The first integral being finite, we may focus on the second one,
in particular on the sum over the excursions $\{I_{n}\}_{n\in\N}$.
Consider first the excursions such that $|I_{n}|\ge\frac{1}{\sqrt{2}}$,
for which $\min\{|I_n|, s\} = s$:
there are at most $2/(1/\sqrt{2}) + 1 = 2 \sqrt{2} + 1$ such
excursions with
$I_{n} \cap(0,2) \ne\varnothing$, therefore
\[
\int_0^{1/\sqrt{2}} \sum_{n\in\N\dvtx I_n \cap(0,2) \ne\varnothing
,
|I_{n}| \ge{1}/{\sqrt{2}}}
\frac{\min\{|I_n|, s\}}{s^{1+\gep}} \,\dd s \le
\bigl(2\sqrt{2} + 1\bigr) \int_0^{1/\sqrt{2}} \frac{1}{s^{\gep}} \,\dd s
< \infty.
\]
Plainly, also the last excursion $I_{n} \ni2$ gives a finite contribution.
It remains to consider the excursions $I_{n}$ included in $(0,2)$
such that $|I_{n}| < \frac{1}{\sqrt{2}}$, for which we may write
\begin{eqnarray*}
&& \int_0^{{1}/{\sqrt{2}}} \sum_{I_n \subseteq(0,2),
|I_n| < {1}/{\sqrt{2}}}
\frac{\min\{|I_n|, s\}}{s^{1+\gep}} \,\dd s \\
&&\qquad= \sum_{I_n \subseteq
(0,2),
|I_n| < {1}/{\sqrt{2}}}
\biggl( \int_0^{|I_n|} \frac{1}{s^{\gep}} \,\dd s
+ \int_{|I_n|}^{{1}/{\sqrt{2}}}
\frac{|I_n|}{s^{1+\gep}} \,\dd s \biggr) \\
&&\qquad =
\sum_{I_n \subseteq(0,2), |I_n| < {1}/{\sqrt{2}}}
\biggl( \frac{|I_n|^{1-\gep}}{1-\gep} + \frac{1}{\gep} |I_n|
\biggl( \frac{1}{|I_n|^\gep} - \bigl(\sqrt{2}\bigr)^\gep\biggr)
\biggr)\\
&&\qquad\le \frac{1}{\gep(1-\gep)} \sum_{I_n \subseteq(0,2)}
|I_n|^{1-\gep} .
\end{eqnarray*}
We have thus shown that there exist constants $0 < a,b < \infty$
(depending on $\gep$) such that
%
%e4.32 ###
%
\begin{equation} \label{eq:estM}
\cM\le a + b
\sqrt{\sum_{n \in\N\dvtx I_n \subseteq(0,2)} |I_n|^{1-\gep}} .
\end{equation}

We can finally conclude the proof of (\ref{eq:ttoo}). From
(\ref{eq:appartial}), (\ref{eq:appartial2}) and (\ref{eq:estM})
it follows that equation (\ref{eq:ttoo}) is proven once we
show that for every $C >0$ there exists
$D = D(C) \in(0,\infty)$ such that for every $T > 0$
%
%e4.33 ###
%
\begin{equation}
\bE\Biggl[ \exp\Biggl( C \sqrt{T}
\sqrt{\sum_{n\in\N\dvtx I_n \subseteq(0,2)}
|I_n|^{1-\gep}} \Biggr) \Biggr] \le D \exp(D T) .
\end{equation}
But this is a direct consequence of equation (\ref{eq:bleach}).

%%%%%%%%%%%%%%%%%%%%%%%%%%%%%%%%%%%%%%%%%%%%%%%%%%%%%%%%%%%%%%%%%%%%%%%%%%%%%%%%
%%%%%%%%%%%%%%%%%%%%%%%%%%%%%%%%%%%%%%%%%%%%%%%%%%%%%%%%%%%%%%%%%%%%%%%%%%%%%%%%
%%%%%%%%%%%%%%%%%%%%%%%%%%%%%%%%%%%%%%%%%%%%%%%%%%%%%%%%%%%%%%%%%%%%%%%%%%%%%%%%

%s5 ###
\section{\texorpdfstring{Proof of Lemma 3.5%\protect\ref{th:crucialk}
}{Proof of Lemma 3.5}}
\label{sec:crucialk}

We recall that $\tau= \{\tau_{n}\}_{n\in\N}$ and $\widetilde\tau
^\ga$ denote,
respectively, the renewal process and the regenerative set, both defined
under the law $\bP$. For $x \ge0$, we denote by $\bP_x$
the law of the sets $\tau$ and $\widetilde\tau^\ga$ started at $x$,
that is, $\bP_x(\tau\in\cdot) := \bP(\tau+ x = \{\tau_{n} + x\}
_{n\in\N} \in\cdot)$ and analogously for $\widetilde\tau^\ga$.
For the definition of the
vectors $\Sigma:= (m; s_{1}, \ldots, s_{m}; \underline\gs_{1},
\ldots, \underline\gs_{m})$ and
$\widetilde\Sigma:= (\widetilde m; \widetilde s_{1}, \ldots,
\widetilde s_{\widetilde m};
\widetilde{\underline\gs}_{1}, \ldots, \widetilde{\underline\gs
}_{\widetilde m})$, we refer
to Section \ref{sec:step3}.

In this section, we fix $\delta= 1$.
We have to estimate the Radon--Nikodym density $\frac{\dd\widetilde
\Sigma}{\dd\Sigma}$
%$\frac{\dd\widetilde\Sigma}{\dd\Sigma}$
of the laws of $\widetilde\Sigma$ and $\Sigma$ [which does not depend
on the sign variables; see the explanation between (\ref{eq:H3a})
and (\ref{eq:H2modbis})], namely the quantity
%
%e5.1 ###
%
\begin{equation} \label{eq:appratio}
\frac{\dd\widetilde\Sigma}{\dd\Sigma} (l; x_1, \ldots, x_l)
=
\frac{\bP( (\widetilde m; \widetilde{\underline\gs}_{1},
\ldots,
\widetilde{\underline\gs}_{m}) = (l; x_1, \ldots, x_l) )}
{\bP( (m; \underline\gs_{1}, \ldots, \underline\gs_{m})
= (l; x_1, \ldots, x_l) )} .
\end{equation}
Note that by construction
$({\underline\gs}_{i+1} - {\underline\gs}_i) \in[\delta,\infty)
\cap\gep\N$,
and since $\delta= 1$ we assume that $x_{i+1} - x_i \in[1,\infty)
\cap\gep\N$.
Using the regenerative property of $\widetilde\tau^\ga$ and the renewal
property of $\tau$, the ratio in (\ref{eq:appratio})
can be estimated in terms of the probability
of the first coarse-grained returns of $\widetilde\tau^\ga$ and
$\tau$:
%
%e5.2 ###
%
\begin{equation} \label{eq:boundSigma0}
\prod_{i=1}^{l} c(x_{i} - x_{i-1}) \le
\frac{\dd\widetilde\Sigma}{\dd\Sigma} (l; x_1, \ldots, x_l)
\le
\prod_{i=1}^{l} C(x_{i} - x_{i-1}) ,
\end{equation}
where we set for convenience $x_0 := 0$ and we have introduced, for $z
\in
[1, \infty) \cap\gep\N$,
%
%e5.3 ###
%
\begin{equation}
C(z) := \sup_{y, \widetilde y \in(0, \gep]}
\frac{\bP_{\widetilde y} ( \inf\{ u > 1\dvtx u \in\widetilde
\tau^\ga\}
\in (z, z+\gep] )}
{\bP_{{y}/{a^2}} ( \inf\{ i > {1}/{a^2}\dvtx i \in\tau
\}
\in ({z}/{a^2}, ({z+\gep})/{a^2}] )} ,
\end{equation}
and $c(z)$ is defined analogously, replacing the supremum (over $y$
and $\widetilde y$) by the infimum (over the same variables and range).
For the purpose of proving Lem\-ma~\ref{th:crucialk},
it is actually more convenient to give a slightly different estimate
than~(\ref{eq:boundSigma0}), namely
%
%e5.4 ###
%
\begin{equation}
\exp\Biggl( -\sum_{i=1}^l G(x_i - x_{i-1}) \Biggr) \le
\frac{\dd\widetilde\Sigma}{\dd\Sigma} (l; x_1, \ldots, x_l)
\le
\exp\Biggl( \sum_{i=1}^l G(x_i - x_{i-1}) \Biggr) ,\hspace*{-34pt}
\end{equation}
where $G(z) = G_{\gep, a}(z)$ is defined, always for $z \in
[1, \infty) \cap\gep\N$, by
%
%e5.5 ###
%
\begin{eqnarray} \label{eq:parapeo}
G_{\gep, a}(z) &:=& \sup_{y, \widetilde y \in(0, \gep]}
\biggl| \log\biggl(
\biggl({\bP_{{y}/{a^2}} \biggl( \inf\biggl\{ i > \frac{1}{a^2}\dvtx i
\in\tau\biggr\}
\in \biggl(\frac{z}{a^2}, \frac{z+\gep}{a^2}\biggr] \biggr)}\biggr)\hspace*{-30pt}\nonumber\\[-8pt]\\[-8pt]
&&\qquad\quad\hspace*{40.8pt}{}\times
\bigl({\bP_{\widetilde y} \bigl( \inf\{ u > 1\dvtx u \in\widetilde\tau
^\ga\}
\in (z, z+\gep] \bigr)}\bigr)^{-1} \biggr) \biggr|.\hspace*{-30pt}\nonumber
\end{eqnarray}
Recalling the statement of Lemma \ref{th:crucialk}, we are left with
showing that
%
%e5.6 ###
%
\begin{equation} \label{eq:crucialk2bis}
G_{\gep, a}(z) \le \kappa(\gep, a)
( \log z + 1 ) \qquad \mbox{with }
\lim_{\gep\to0} \limsup_{a \to0} \kappa(\gep, a)
= 0 .
\end{equation}

We claim that the rescaled renewal process $a^{2} \tau= \{a^{2} \tau
_{n}\}_{n\in\N}$,
viewed as a random closed subset of $[0,\infty)$,
converges in distribution as $a \to0$ toward the regenerative set
$\widetilde\tau^\ga$,
where we equip the family of closed subsets of $[0,\infty)$ with
the topology of Matheron, as described in \cite{cfFFM}.
To check this claim,
we recall from Appendix \ref{sec:regapp} that $\widetilde\tau^\ga$
is the closure of the
image of the (stable) subordinator with L\'evy exponent $\Phi(\gl) :=
\lambda^\ga$.
If we denote by $\{N_t\}_{t\ge0}$ a standard Poisson process on $\R$
of rate $\gamma> 0$,
independent of all the processes considered so far, the random set
$a^2 \tau$ can be viewed as the image of the subordinator $\{a^{2}
\tau_{N_t}\}_{t \ge0}$,
whose L\'evy exponent is given by
%
%e5.7 ###
%
\begin{eqnarray}
\Phi_a(\gl) :\!&=& - \log\bE[e^{-\gl a^2 \tau_{N_1}} ]
= \gamma( 1 - \bE[e^{-\gl a^2 \tau_1} ] )
\nonumber\\[-8pt]\\[-8pt]
&=& \gamma\sum_{n\in\N} ( 1-e^{-\lambda a^2 n} ) K(n).\nonumber
\end{eqnarray}
If we fix $\gamma= \gamma(a)$ so that
$\Phi_a(1) = 1$, as prescribed by Proposition (1.14) in \cite{cfFFM},
it follows easily by our assumption
(\ref{eq:K}) that $\lim_{a\to0} \Phi_a(\gl) = \Phi(\gl) = \gl
^\ga$
for every $\gl\ge0$. By Proposition (3.9) in \cite{cfFFM},
the pointwise convergence of the L\'evy exponents entails the
convergence in distribution of the corresponding regenerative sets,
which proves the claim.

From the convergence in distribution of $a^2 \tau$ toward $\widetilde
\tau^\ga$
it follows that the numerator in the right-hand side of (\ref
{eq:parapeo}) converges as $a \to0$
toward the denominator with $\widetilde y$ replaced by $y$,
for all fixed $\gep\in(0,1)$, $z \in[1,\infty) \cap\gep\N$ and
$y \in(0, \gep]$.
In the following lemma, we provide a quantitative control
on this convergence, as a function of $z$ and $y$.
\begin{lemma} \label{lem:pippo}
Fix $\gep\in(0, 1/3)$. There exists
$\zeta_\gep(a) > 0$ with $\lim_{a\to0} \zeta_\gep(a) = 0$ such that
%
%e5.8 ###
%
\begin{eqnarray}\quad \label{eq:toprove1}
\bigl(1 - \zeta_\gep(a) \bigr) z^{-\zeta_\gep(a)} &\le&
\frac{\bP_{{y}/{a^2}} ( \inf\{ i > {1}/{a^2}\dvtx i
\in\tau\}
\in ({z}/{a^2}, ({z+\gep})/{a^2}] )}
{\bP_{y} ( \inf\{ u > 1\dvtx u \in\widetilde\tau^\ga\}
\in (z, z+\gep] )} \nonumber\\[-8pt]\\[-8pt]
&\le&
\bigl( 1 + \zeta_\gep(a) \bigr) z^{\zeta_\gep(a)} ,\nonumber
\end{eqnarray}
for all $a \in(0,a_0)$ (with $a_0 > 0$), $y \in[0,1/3]$ and
$z \in[1,\infty) \cap\gep\N$.
\end{lemma}

We point out that Lemma \ref{lem:pippo} is proved below through explicit
estimates, without reference to the convergence in distribution
of $a^2 \tau$ toward $\widetilde\tau^\ga$ stated above.

We now apply (\ref{eq:toprove1}) to (\ref{eq:parapeo}):
since $|{\log}(1+x)| \le2 |x|$ for $x$ small, for small $a$ we obtain
%
%e5.9 ###
%
\begin{eqnarray} \label{eq:hehe}
G_{\gep, a}(z) &\le& 2 \zeta_\gep(a) ( \log z + 1
)\nonumber\\[-8pt]\\[-8pt]
&&{} +
\sup_{y, \widetilde y \in(0, \gep]} \biggl| \log\biggl(
\frac{\bP_{ y} ( \inf\{ u > 1\dvtx u \in\widetilde\tau^\ga\}
\in (z, z+\gep] )}
{\bP_{\widetilde y} ( \inf\{ u > 1\dvtx u \in\widetilde\tau
^\ga\}
\in (z, z+\gep] )} \biggr) \biggr| .\nonumber
\end{eqnarray}
Recalling the definition (\ref{eq:gtds}) of $d_t(\widetilde\tau^\ga)$
and applying (\ref{eq:D}), for $z \in[1,\infty) \cap\gep\N$ we
can write
%
%e5.10 ###
%
\begin{eqnarray} \label{eq:pint0}
&&\bP_{y} \bigl( \inf\{ u > 1\dvtx u \in\widetilde\tau^\ga\}
\in (z, z+\gep] \bigr)
\nonumber\\[-8pt]\\[-8pt]
&&\qquad= \frac{\sin(\pi\ga)}{\pi} \int_{z}^{z+\gep}
\frac{(1 - y)^\ga}{(t-1)^{\ga} (t - y)} \,\dd t .\nonumber
\end{eqnarray}
From this explicit expression it is easy to check that
the second term in the right-hand side of (\ref{eq:hehe})
vanishes as $\gep\to0$
\textit{uniformly in $z \in[1,\infty) \cap\gep\N$},
hence (\ref{eq:crucialk2bis}) holds true.

%%%%%%%%%%%%%%%%%%%%%%%%%%%%%%%%%%%%%%%%%%%%%%%%%%%%%%%%%%%%%%%%%%%%%
%s5.1 ###
\subsection{\texorpdfstring{Proof of Lemma \protect\ref{lem:pippo}}{Proof of Lemma B.1}}

We have already obtained in (\ref{eq:pint0}) an explicit expression
for the denominator in (\ref{eq:toprove1}). It is however more convenient
to give an alternative expression: recalling again
the definition (\ref{eq:gtds}) of the variable $d_t(\widetilde\tau
^\ga)$
and applying (\ref{eq:joint}), we can rewrite the denominator in
(\ref{eq:toprove1}) as
%
%e5.11 ###
%
\begin{equation} \label{eq:pint1}
I(y,z) :=
\frac{\ga \sin(\pi\ga)}{\pi}
\int_{y}^{1} \dd s \int_{z}^{z+\gep} \dd t
\frac{1}{(s - y)^{1-\ga} (t - s)^{1+\ga}} .
\end{equation}
Recalling that $K(n) := \bP(\tau_1 = n)$ and setting $U(n) := \bP( n
\in\tau)$,
we can rewrite the numerator in (\ref{eq:toprove1}) using the renewal
property as
%
%e5.12 ###
%
\begin{eqnarray} \label{eq:pint2}
J_a(y,z) :\!&=& \mathop{\sum_{{y}/{a^2} \le k \le
{1}/{a^2}}}_{{z}/{a^2} < l \le({z+\gep})/{a^2}}
U\biggl(k - \frac{y}{a^2}\biggr)
K(l - k)\nonumber\\[-8pt]\\[-8pt]
&=& \mathop{\sum_{s \in[y,1] \cap a^2 \N}}_{t \in(z,
z+\gep] \cap a^2 \N}
U \biggl( \frac{1}{a^2}(s-y) \biggr)
K \biggl( \frac{1}{a^2}(t -s) \biggr) .\nonumber
\end{eqnarray}
We now use \cite{cfDon}, Theorem B, coupled with our basic
assumption on the inter-arrival distribution (\ref{eq:K}), to see that
% K(\ell) \stackrel{\ell\to\infty}\sim\frac{L(\ell)}{\ell^{1+
%with $L(\cdot)$ slowly varying, hence by
%
%e5.13 ###
%
\begin{equation}
\label{eq:green}
U(\ell) \stackrel{\ell\to\infty}\sim
\frac{\ga \sin(\pi\ga)}{\pi}
\frac{1}{L(\ell) \ell^{1-\ga}} .
\end{equation}
Using the asymptotic relations (\ref{eq:K}) and (\ref{eq:green}) and
a Riemann sum argument (with some careful handling of the slowly
varying functions, see the details below), one can check that
(\ref{eq:pint2}) converges toward (\ref{eq:pint1}) as $a \to0$,
for all \textit{fixed} $\gep\in(0,1/3)$,
$z \in[1, \infty) \cap\gep\N$ and $y \in(0, \gep]$.
However to obtain (\ref{eq:toprove1}), a more attentive estimate is required.
We set $n := 1/a^2$ for notational convenience, so that,
with some abuse of notation, we can rewrite (\ref{eq:pint2}) as
%
%e5.14 ###
%
\begin{eqnarray} \label{eq:finalaim}
J_n(y,z) :\!&=& \mathop{\sum_{n y \le k \le n}}_{nz < l \le n(z+\gep)}
U(k - ny)
K(l - k) \nonumber\\[-8pt]\\[-8pt]
&=& \mathop{\sum_{s \in[y,1] \cap{1}/{n}\N}}_{t
\in(z, z+\gep] \cap{1}/{n}\N}
U \bigl( n(s-y) \bigr)
K \bigl( n(t -s) \bigr) .\nonumber
\end{eqnarray}
We can now rephrase (\ref{eq:toprove1}) in the
following way: for every fixed $\gep\in(0,1/3)$, there
exist $\zeta_\gep(n) > 0$, with $\lim_{n\to\infty} \zeta_\gep(n)
= 0$,
and $n_0 \in\N$ such that
%
%e5.15 ###
%
\begin{equation} \label{eq:toprove2}
\bigl(1 - \zeta_\gep(n) \bigr) z^{-\zeta_\gep(n)} \le
\frac{J_n(y,z)}{I(y,z)} \le
\bigl( 1 + \zeta_\gep(n) \bigr) z^{\zeta_\gep(n)} ,
\end{equation}
for all $n\ge n_0$, $y \in[0,1/3]$ and $z \in[1,\infty) \cap\gep\N$.
We recall that $I(y,z)$ is defined in (\ref{eq:pint1}).
For convenience, we divide the rest of the proof in three steps.

%%%%%%%%%%%%%%%%%%%%%%%%%%%%%%%%%%%%%%%%%%%%%%%%%%%%%%%%

%s5.1.1 ###
\subsubsection*{Step 1}

We first show that the terms in (\ref{eq:finalaim}) with $k \le ny +
\sqrt{n}$,
that is,
%
%e5.16 ###
%
\begin{equation} \label{eq:finalaim1}
A_n(y,z) := \mathop{\sum_{n y \le k \le ny + \sqrt{n}}}_{nz < l
\le n(z+\gep)}
U(k - ny) K(l - k)
\end{equation}
give a negligible contribution to (\ref{eq:toprove2}).

By paying a positive constant, we can replace $K(\cdot)$
and $U(\cdot)$ by their asymptotic behaviors; cf. (\ref{eq:K})
and (\ref{eq:green}). Note that $k \le ny + \sqrt{n} \le n/2$ for large
$n$, because $y \le1/3$, and therefore $n(z-1/2) \le(l-k) \le n(z+1/3)$,
because $\gep\le1/3$, for all $l,k$ in the range of summation.
We thus obtain the upper bound
%
%e5.17 ###
%
\begin{equation} \label{eq:boundA1}
A_n(y,z) \le C_1 \sum_{0 < h \le\sqrt{n}}
\frac{1}{L(h) h^{1-\ga}}
\sum_{n(z-{1}/{2}) < m \le n(z+{1}/{3})} \frac
{L(m)}{m^{1+\ga}} ,
\end{equation}
for some absolute constant $C_1>0$. We now show that,
for some absolute constant $C_2 > 0$ (not depending on $z$),
we can write $L(m) \le C_2 L(nz)$ for every $m$ in the range
of summation. To this purpose, we recall
the representation theorem of slowly varying functions:
%
%e5.18 ###
%
\begin{eqnarray} \label{eq:RepSV}
L(x) = a(x) \exp\biggl( \int_1^x \frac{b(t)}{t} \,\dd t
\biggr) \nonumber\\[-8pt]\\[-8pt]
\eqntext{\mbox{with } \displaystyle\lim_{x\to\infty} a(x) \in(0,\infty)
\mbox{ and } \lim_{x\to\infty} b(x) = 0 ;}
\end{eqnarray}
see Theorem 1.3.1 in \cite{cfBinGolTeu}.
Setting $\gamma_n := \sup_{x \ge n/2} |b(x)|$,
we have $\lim_{n\to\infty} \gamma_n = 0$ and
for $m \in\{n(z-1/2), n(z+1/3)\}$ we can write for $z\ge1$
%
%e5.19 ###
%
\begin{eqnarray} \label{eq:LL}
\frac{L(m)}{L(nz)} &\le& \frac{a(m)}{a(nz)}
\exp\biggl( \gamma_n \int_{n(z-1/2)}^{n(z+1/3)}
\frac{1}{t} \,\dd t \biggr)\nonumber\\[-8pt]\\[-8pt]
&\le& \frac{\sup_{k \ge n/2} a(k)}{\inf_{k \ge n}a(k)}
\exp\biggl( \gamma_n \log\frac{z+{1}/{3}}{z-{1}/{2}}\biggr) .\nonumber
\end{eqnarray}
Since $z\ge1$, it is clear that the right-hand side of (\ref{eq:LL})
is bounded from above by some absolute constant $C_2$
(in fact, it even converges to $1$ as $n\to\infty$).
From (\ref{eq:boundA1}), we then obtain
%
%e5.20 ###
%
\begin{eqnarray} \label{eq:intermA}\quad
A_n(y,z) & \le & C_2 L(nz) \sum_{0 < h \le\sqrt{n}}
\frac{1}{L(h) h^{1-\ga}}
\sum_{n(z-{1}/{2}) < m \le n(z+{1}/{3})} \frac{1}{m^{1+\ga
}} \nonumber\\[-8pt]\\[-8pt]
& \le & C_3 \frac{L(nz)}{n^\ga z^{1+\ga}}
\sum_{0 < h \le\sqrt{n}} \frac{1}{L(h) h^{1-\ga}}
\le C_4 \frac{L(nz)}{n^\ga z^{1+\ga}}
\frac{n^{\ga/2}}{\ga L(\sqrt n)} ,\nonumber
\end{eqnarray}
where $C_3, C_4$ are absolute positive constant and
the last inequality is a classical result (Proposition 1.5.8
in \cite{cfBinGolTeu}). Using again the representation
(\ref{eq:RepSV}), in analogy with (\ref{eq:LL}), we can write
%
%e5.21 ###
%
\begin{eqnarray} \label{eq:LL2}
\frac{L(nz)}{L(\sqrt n)} &\le& \frac{a(nz)}{a(\sqrt n)}
\exp\biggl( \gamma_n \int_{\sqrt n}^{nz} \frac{1}{t} \,\dd t\biggr)
\le C_5 \exp\biggl( \gamma_n \log\frac{nz}{\sqrt n}
\biggr)\nonumber\\[-8pt]\\[-8pt]
&=& C_5 n^{\gamma_n/2} z^{\gamma_n} ,\nonumber
\end{eqnarray}
for some absolute constant $C_5$. Coming back to (\ref{eq:intermA}),
we have shown that there exists absolute constants $C_6$ and $n_0$
such that for all
$n\ge n_0$, $z \in[1,\infty) \cap\gep\N$ and $y \in[0,1/3]$
%
%e5.22 ###
%
\begin{equation} \label{eq:boundA2}
A_n(y,z) \le \frac{C_6}{n^{(\ga- \gamma_n)/2}}
\frac{z^{\gamma_n}}{z^{1+\ga}} .
\end{equation}

Let us now look back at the integral $I(y,z)$, defined in (\ref{eq:pint1}).
It is easy to check that for every fixed
$\gep\in(0,1/3)$ there exists an absolute constant $C_7 = C_7(\gep)
> 0$
such that
%
%e5.23 ###
%
\begin{equation} \label{eq:boundI0}
I(y,z) \ge \frac{C_7}{z^{1+\ga}} ,
\end{equation}
for all $y \in[0,1/3]$ and $z \in[1,\infty) \cap\gep\N$.
If we set $\zeta'(n) := \max\{\gamma_n,
C_6 /\break (C_7 n^{(\ga- \gamma_n)/2})\}$,
we have $\lim_{n\to\infty} \zeta'(n) = 0$ and from (\ref{eq:boundA2})
and (\ref{eq:boundI0}) we have shown that for every fixed
$\gep\in(0,1/3)$ there exists $n_0 \in\N$ such that for $n\ge n_0$
we have
%
%e5.24 ###
%
\begin{equation} \label{eq:toprove2.1}
\frac{A_n(y,z)}{I(y,z)} \le \zeta'(n) z^{\zeta'(n)} ,
\end{equation}
for all $z \in[1,\infty) \cap\gep\N$ and $y \in[0,1/3]$.
This completes the first step.

%%%%%%%%%%%%%%%%%%%%%%%%%%%%%%%%%%%%%%%%%%%%%%%%%%%%%%%%

%s5.1.2 ###
\subsubsection*{Step 2}

We now consider the terms in (\ref{eq:finalaim}) with
$k > ny + \sqrt{n}$, or equivalently $s > y + \frac{1}{\sqrt n}$,
that is, we introduce the quantity
%
%e5.25 ###
%
\begin{equation} \label{eq:nel}
B_n(y,z) :=
\mathop{\sum_{s \in(y + {1}/{\sqrt n},1] \cap{1}/{n}\N
}}_{t \in(z, z+\gep] \cap{1}/{n}\N}
U \bigl( n(s-y) \bigr)
K \bigl( n(t -s) \bigr) ,
\end{equation}
and we observe that $J_n(y,z) = A_n(y,z) + B_n(y,z)$, see
(\ref{eq:finalaim}) and (\ref{eq:finalaim1}).
Our aim is to prove (\ref{eq:toprove2}): in view of relation
(\ref{eq:toprove2.1}), it remains
to show that for every fixed $\gep\in(0,1/3)$ there
exist $\zeta''(n) > 0$, with $\lim_{n\to\infty} \zeta''(n) = 0$,
and $n_0 \in\N$ such that
%
%e5.26 ###
%
\begin{equation} \label{eq:toprove3}
\bigl(1 - \zeta''(n) \bigr) z^{-\zeta''(n)} \le
\frac{B_n(y,z)}{I(y,z)} \le
\bigl( 1 + \zeta''(n) \bigr) z^{\zeta''(n)} ,
\end{equation}
for all $n\ge n_0$, $y \in[0,1/3]$ and
$z \in[1,\infty) \cap\gep\N= \{1, 1+\gep, 1+2\gep, \ldots\}$.
In this step, we prove that (\ref{eq:toprove3}) holds
for $z \in[1+\gep, \infty) \cap\gep\N$, that
is we exclude the case $z = 1$, that will be considered separately
in the third step.

By construction, the arguments of the functions $U(\cdot)$
and $K(\cdot)$ appearing in (\ref{eq:nel})
tend to $\infty$ as $n\to\infty$ uniformly in the range of
summation: in fact $n(s-y) \ge\sqrt n$ and $n(t -s) \ge\gep n$,
because we assume that $z \ge1+\gep$.
We can therefore replace $U(\cdot)$ and $K(\cdot)$ by their asymptotic
behaviors, given in (\ref{eq:K}) and (\ref{eq:green}), by committing
an asymptotically negligible error: more precisely, we can write
%
%e5.27 ###
%
\begin{eqnarray} \label{eq:mezzo}\qquad
B_n(y,z) &=& \bigl(1+o(1) \bigr) \frac{C_\ga}{n^2}
\mathop{\sum_{s \in(y + {1}/{\sqrt n},1] \cap{1}/{n}\N
}}_{t \in(z, z+\gep] \cap{1}/{n}\N}
\biggl[ \frac{L ( n(s-y) )}{L ( n(t -s) )} \biggr]\nonumber\\[-8pt]\\[-8pt]
&&\hspace*{139pt}{}\times\frac{1}{(s - y)^{1-\ga} (t - s)^{1+\ga}} ,\nonumber
\end{eqnarray}
where we set $C_\ga:= \ga\sin(\pi\ga)/\pi$ for short
and where, here and in the sequel, $o(1)$~denotes a quantity
(possibly depending on $\gep$ and varying from place to place)
that vanishes as $n \to\infty$ \textit{uniformly in $y \in[0,1/3]$
and in $z \in[1+\gep, \infty) \cap\gep\N$}.

We now estimate the ratio in square brackets in the right-hand side of
(\ref{eq:mezzo}).
Recalling the representation theorem of slowly varying functions
[see (\ref{eq:RepSV})] uniformly for $s,t$ in the range of summation, we
can write
%
%e5.28 ###
%
\begin{equation} \label{eq:ratio}
\frac{L ( n(s-y) )}{L ( n( t -s) )} =
\bigl(1+o(1)\bigr) \exp\biggl( \int_{n(t -s)}^{n(s-y)}
\frac{b(x)}{x} \,\dd x \biggr) ,
\end{equation}
with the convention $\int_\gb^\gamma(\cdots) := - \int_\gamma^\gb
(\cdots)$ if $\gb>\gamma$.
Let us set
%
%e5.29 ###
%
\begin{equation} \label{eq:etan}
\eta_n := {\sup_{x \ge\min\{\sqrt n, \gep n\}}} |b(x)| ,
\end{equation}
so that $\eta_n \to0$ as $n\to\infty$.
Uniformly for $s,t$ in the range of summation, we can write
%
%e5.30 ###
%
\begin{eqnarray} \label{eq:ratioub}
\biggl| \int_{n(t-s)}^{n(s-y)} \frac{b(x)}{x} \,\dd x \biggr| &\le&
\eta_n \biggl| \int_{n(t-s)}^{n(s-y)} \frac{1}{x} \,\dd x
\biggr|\nonumber\\[-8pt]\\[-8pt]
&\le& \eta_n \bigl( |{\log}(t -s)| + |{\log}(s-y)| \bigr) .\nonumber
\end{eqnarray}
In the range of summation of (\ref{eq:mezzo}),
we have $0 < (s-y) \le1$, hence $|{\log}(s-y)|
= -\log(s-y)$, and $\gep\le(t - s) \le z+\gep$, whence
$|{\log}(t-s)| \le-\log\gep+ \log(z+\gep)$ (recall that $\gep< 1 < z$).
Coming back to (\ref{eq:mezzo}), from (\ref{eq:ratio})
and (\ref{eq:ratioub}) we obtain the upper bound
%
%e5.31 ###
%
\begin{eqnarray} \label{eq:del}
B_n(y,z) &\le& \bigl(1+o(1) \bigr) \frac{(z+\gep)^{\eta
_n}}{\gep^{\eta_n}}\nonumber\\[-8pt]\\[-8pt]
&&{}\times\biggl[ \frac{C_\ga}{n^2}
\mathop{\sum_{s \in(y + {1}/{\sqrt n},1] \cap{1}/{n}\N
}}_{t \in(z, z+\gep] \cap{1}/{n}\N}
\frac{1}{(s - y)^{1-\ga+\eta_n} (t - s)^{1+\ga}} \biggr]
,\nonumber
\end{eqnarray}
as well as the corresponding lower bound
%
%e5.32 ###
%
\begin{eqnarray} \label{eq:del2}
B_n(y,z) &\ge& \bigl(1+o(1) \bigr) \frac{\gep^{\eta_n}}{(z+\gep
)^{\eta_n}}\nonumber\\[-8pt]\\[-8pt]
&&{}\times\biggl[ \frac{C_\ga}{n^2}
\mathop{\sum_{s \in(y + {1}/{\sqrt n},1] \cap{1}/{n}\N
}}_{t \in(z, z+\gep] \cap{1}/{n}\N}
\frac{1}{(s - y)^{1-\ga-\eta_n} (t - s)^{1+\ga}} \biggr].\nonumber
\end{eqnarray}

Observe that we can write
$\frac{(z+\gep)^{\eta_n}}{\gep^{\eta_n}} = c_{\gep, z, n}
z^{\eta_n}$,
with $c_{\gep, z, n} = (\frac{1}{\gep} +
\frac{1}{z})^{\eta_n} \to1$ as $n\to\infty$ (for fixed $\gep$) uniformly
in $z \in[1,\infty)$. We can therefore
incorporate $c_{\gep, z, n}$ in the $(1+o(1))$ term
in (\ref{eq:del}) and (\ref{eq:del2}).
Recalling that we aim at proving (\ref{eq:toprove3}),
it remains to show that
for every fixed $\gep\in(0,1/3)$ the terms in
square brackets in the right-hand sides of
(\ref{eq:del}) and (\ref{eq:del2}), divided by the integral
$I(y,z)$ defined in (\ref{eq:pint1}), converge to $1$ as $n\to\infty$
uniformly in $y \in[0,1/3]$ and in $z \in[1+\gep, \infty)$.

Since the summand in the right-hand side
of (\ref{eq:del}) is decreasing in $t$, we can replace the sum over
$t$ by
an integral over a slightly shifted domain, getting the following upper bound
on the term in square brackets in the right-hand side of (\ref{eq:del}):
%
%e5.33 ###
%
\begin{equation} \label{eq:cammin}
[ \cdots]_{(\fontsize{8.36}{10.36}\selectfont{\mbox{\ref{eq:del}}})} \le \int_{z-
1/n}^{z+\gep}
\biggl( \frac{C_\ga}{n} \sum_{s \in(y + {1}/{\sqrt n},1] \cap
{1}/{n}\N}
\frac{1}{(s - y)^{1-\ga+\eta_n} (t - s)^{1+\ga}}
\biggr)\, \dd t .\hspace*{-39pt}
\end{equation}
By direct computation one sees that the term in the right-hand side of this
relation, as a function of $s$, is decreasing in $(0,s_0)$ and
increasing in $(s_0, \infty)$, where
$s_0 = \frac{(1-\ga+\eta_n) t + (1+\ga) y}{2+\eta_n}$.
The precise\vspace*{1pt} value of $s_0$ is actually immaterial: the important
point is that each term in the sum in (\ref{eq:cammin}) can be bounded
from above by an integral over $[s-\frac1n, s]$ (if $s \le s_0$)
or over $[s, s+\frac1n]$ (if $s \ge s_0$). Therefore, we
get an upper bound replacing the sum by
an integral over a slightly enlarged domain:
%
%e5.34 ###
%
\begin{eqnarray} \label{eq:cammin2}\quad
[ \cdots]_{(\fontsize{8.36}{10.36}\selectfont{\mbox{\ref{eq:del}}})} &\le& \frac{\ga \sin(\pi
\ga)}{\pi}\nonumber\\[-8pt]\\[-8pt]
&&{}\times
\int_{y + {1}/{\sqrt n} - 1/n}^{1 + 1/n} \,\dd s
\int_{z- 1/n}^{z+\gep} \dd t
\frac{1}{(s - y)^{1-\ga+\eta_n} (t - s)^{1+\ga}} .\nonumber
\end{eqnarray}
With almost identical arguments one obtains the following lower bound
on the term in square brackets in the right-hand side of (\ref{eq:del2}):
%
%e5.35 ###
%
\begin{eqnarray} \label{eq:cammin3}\quad\qquad
[ \cdots]_{(\fontsize{8.36}{10.36}\selectfont{\mbox{\ref{eq:del2}}})} &\ge& \frac{\ga \sin
(\pi\ga)}{\pi}\nonumber\\[-8pt]\\[-8pt]
&&{}\times
\int_{y + {1}/{\sqrt n} + 1/n}^{1 - 1/n} \,\dd s
\int_{z}^{z + \gep+ 1/n} \,\dd t
\frac{1}{(s - y)^{1-\ga-\eta_n} (t - s)^{1+\ga}} .\nonumber
\end{eqnarray}
One can now check directly that, for every fixed $\gep\in(0,1/3)$,
the ratio between the right-hand side of (\ref{eq:cammin2}) and the integral
$I(y,z)$ defined in (\ref{eq:pint1}) converges to $1$ as $n\to\infty$,
uniformly in $y \in[0,1/3]$ and in $z \in[1+\gep, \infty)$.
Since an analogous statement holds for the right-hand side of (\ref
{eq:cammin3}),
the second step is complete.

%%%%%%%%%%%%%%%%%%%%%%%%%%%%%%%%%%%%%%%%%%%%%%%%%%%%%%%%%

%s5.1.3 ###
\subsubsection*{Step 3}

To complete the proof of Lemma \ref{lem:pippo}, it only remains
to prove that equation (\ref{eq:toprove3}) holds true also for
$z=1$. More explicitly, we have to show that as $n \to\infty$
%
%e5.36 ###
%
\begin{equation} \label{eq:lastrat}
\frac{B_n(y,1)}{I(y,1)} \longrightarrow 1 ,
\end{equation}
uniformly in $y \in[0,1/3]$. We recall that
%
%e5.37 ###
%
\begin{equation} \label{eq:vita}
B_n(y,1) :=
\mathop{\sum_{s \in(y + {1}/{\sqrt n},1] \cap{1}/{n}\N
}}_{t \in(1, 1+\gep] \cap{1}/{n}\N}
U \bigl( n(s-y) \bigr) K \bigl( n(t -s) \bigr),
\end{equation}
while the integral $I(y,z)$ is defined in (\ref{eq:pint1}).

We only sketch the proof of (\ref{eq:lastrat}), because the
arguments are very similar to those used in the preceding
steps. Note that we cannot immediately replace $K(\cdot)$ by its asymptotic
behavior, because its argument $n(t -s)$ can
take small values. It is therefore convenient
to restrict the sum in (\ref{eq:vita}) to
$t \in(1 + 1/\sqrt{n}, 1+\gep]$. For this restricted
sum, call it $B'_n(y,1)$,
one can write a formula analogous to (\ref{eq:mezzo}):
then, arguing as in the second step
(with several simplifications), one shows that
(\ref{eq:lastrat}) holds true with $B_n$ replaced by $B'_n$.
It remains to deal with $B_n - B'_n$, that is, to control
the terms in (\ref{eq:vita}) with $t \le1 + 1/\sqrt{n}$.
In this case one can replace $K(\cdot)$
by its asymptotic behavior by paying a positive constant:
arguing as in the first step, one can show that
$(B_n(y,1) - B'_n(y,1))/I(y,1) \to0$ as $n\to\infty$,
uniformly in $y \in[0,1/3]$.
This completes the proof of (\ref{eq:lastrat}) and of Lemma \ref{lem:pippo}.
\end{appendix}

%%%%%%%%%%%%%%%%%%%%%%%%%%%%%%%%%%%%%%%%%%%%%%%%%%%%%%%%%
\section*{Acknowledgments}
We are very grateful to E. Bolthausen, F. Toninelli, L.~Zambotti
and O. Zeitouni for several useful discussions.

% imsref loaded by lrinkeviciute, 2010-07-22 08:06:44
%

%
\printaddresses


\begin{thebibliography}{26}

%b1 ###
\bibitem{cfAdl}
%
\begin{bbook}[vtex]
\bauthor{\bsnm{Adler},~\bfnm{Robert~J.}\binits{R.~J.}}
(\byear{1990}).
\btitle{An Introduction to Continuity, Extrema, and Related Topics for General
{G}aussian Processes}.
\bpublisher{IMS}, \baddress{Hayward, CA}.
\bid{mr={1088478}}
\end{bbook}
%
\endbibitem

%b2 ###
\bibitem{cfAsm}
%
\begin{bbook}[vtex]
\bauthor{\bsnm{Asmussen},~\bfnm{S{\o}ren}\binits{S.}}
(\byear{2003}).
\btitle{Applied Probability and Queues: Stochastic Modelling and Applied Probability},
\bedition{2nd} ed.
\bseries{Applications of Mathematics (New York)}
\bvolume{51}.
\bpublisher{Springer}, \baddress{New York}.
\bid{mr={1978607}}
\end{bbook}
%
\endbibitem

%b3 ###
\bibitem{cfBPY}
%
\begin{bincollection}[mr]
\bauthor{\bsnm{Barlow},~\bfnm{Martin}\binits{M.}},
\bauthor{\bsnm{Pitman},~\bfnm{Jim}\binits{J.}} \AND
\bauthor{\bsnm{Yor},~\bfnm{Marc}\binits{M.}}
(\byear{1989}).
\btitle{Une extension multidimensionnelle de la loi de l'arc sinus}.
In \bbooktitle{S\'eminaire de {P}robabilit\'es, {XXIII}}.
\bseries{Lecture Notes in Math.}
\bvolume{1372}
\bpages{294--314}.
\bpublisher{Springer}, \baddress{Berlin}.
\bid{doi={10.1007/BFb0083980}, mr={1022918}}
\end{bincollection}
%
\endbibitem

%b4 ###
\bibitem{cfBinGolTeu}
%
\begin{bbook}[mr]
\bauthor{\bsnm{Bingham},~\bfnm{N.~H.}\binits{N.~H.}},
\bauthor{\bsnm{Goldie},~\bfnm{C.~M.}\binits{C.~M.}} \AND
\bauthor{\bsnm{Teugels},~\bfnm{J.~L.}\binits{J.~L.}}
(\byear{1989}).
\btitle{Regular Variation}.
\bseries{Encyclopedia of Mathematics and Its Applications}
\bvolume{27}.
\bpublisher{Cambridge Univ. Press}, \baddress{Cambridge}.
\bid{mr={1015093}}
\end{bbook}
%
\endbibitem

%%b5 ###
%%
%(\byear{2010}).
%%

%b6 ###
\bibitem{cfBG}
%
\begin{barticle}[vtex]
\bauthor{\bsnm{Bodineau},~\bfnm{Thierry}\binits{T.}} \AND
\bauthor{\bsnm{Giacomin},~\bfnm{Giambattista}\binits{G.}}
(\byear{2004}).
\btitle{On the localization transition of random copolymers near selective
interfaces}.
\bjournal{J. Stat. Phys.}
\bvolume{117}
\bpages{801--818}.
\bid{doi={10.1007/s10955-004-5705-7}, mr={2107896}}
\end{barticle}
%
\endbibitem

%b7 ###
\bibitem{cfBGLT}
%
\begin{barticle}[mr]
\bauthor{\bsnm{Bodineau},~\bfnm{Thierry}\binits{T.}},
\bauthor{\bsnm{Giacomin},~\bfnm{Giambattista}\binits{G.}},
\bauthor{\bsnm{Lacoin},~\bfnm{Hubert}\binits{H.}} \AND
\bauthor{\bsnm{Toninelli},~\bfnm{Fabio~Lucio}\binits{F.~L.}}
(\byear{2008}).
\btitle{Copolymers at selective interfaces: New bounds on the phase diagram}.
\bjournal{J. Stat. Phys.}
\bvolume{132}
\bpages{603--626}.
\bid{doi={10.1007/s10955-008-9579-y}, mr={2429695}}
\end{barticle}
%
\endbibitem

%b8 ###
\bibitem{cfBCT}
%
\begin{barticle}[mr]
\bauthor{\bsnm{Bolthausen},~\bfnm{Erwin}\binits{E.}},
\bauthor{\bsnm{Caravenna},~\bfnm{Francesco}\binits{F.}} \AND
\bauthor{\bparticle{de~}\bsnm{Tili{\`e}re},~\bfnm{B{\'e}atrice}\binits{B.}}
(\byear{2009}).
\btitle{The quenched critical point of a diluted disordered polymer model}.
\bjournal{Stochastic Process. Appl.}
\bvolume{119}
\bpages{1479--1504}.
\bid{doi={10.1016/j.spa.2008.07.008}, mr={2513116}}
\end{barticle}
%
\endbibitem

%b9 ###
\bibitem{cfBdH}
%
\begin{barticle}[mr]
\bauthor{\bsnm{Bolthausen},~\bfnm{Erwin}\binits{E.}} \AND
\bauthor{\bparticle{den }\bsnm{Hollander},~\bfnm{Frank}\binits{F.}}
(\byear{1997}).
\btitle{Localization transition for a polymer near an interface}.
\bjournal{Ann. Probab.}
\bvolume{25}
\bpages{1334--1366}.
\bid{doi={10.1214/aop/1024404516}, mr={1457622}}
\end{barticle}
%
\endbibitem

%b10 ###
\bibitem{cfCGG}
%
\begin{barticle}[mr]
\bauthor{\bsnm{Caravenna},~\bfnm{Francesco}\binits{F.}},
\bauthor{\bsnm{Giacomin},~\bfnm{Giambattista}\binits{G.}} \AND
\bauthor{\bsnm{Gubinelli},~\bfnm{Massimiliano}\binits{M.}}
(\byear{2006}).
\btitle{A numerical approach to copolymers at selective interfaces}.
\bjournal{J. Stat. Phys.}
\bvolume{122}
\bpages{799--832}.
\bid{doi={10.1007/s10955-005-8081-z}, mr={2213950}}
\end{barticle}
%
\endbibitem

%b11 ###
\bibitem{cfCY}
%
\begin{barticle}[mr]
\bauthor{\bsnm{Comets},~\bfnm{Francis}\binits{F.}} \AND
\bauthor{\bsnm{Yoshida},~\bfnm{Nobuo}\binits{N.}}
(\byear{2005}).
\btitle{Brownian directed polymers in random environment}.
\bjournal{Comm. Math. Phys.}
\bvolume{254}
\bpages{257--287}.
\bid{doi={10.1007/s00220-004-1203-7}, mr={2117626}}
\end{barticle}
%
\endbibitem

%b15 ###
\bibitem{cfdH}
%
\begin{bbook}[vtex]
\bauthor{\bparticle{den }\bsnm{Hollander},~\bfnm{Frank}\binits{F.}}
(\byear{2009}).
\btitle{Random Polymers}.
\bseries{Lecture Notes in Math.}
\bvolume{1974}.
\bpublisher{Springer}, \baddress{Berlin}.
%2007}.
\bid{doi={10.1007/978-3-642-00333-2}, mr={2504175}}
\end{bbook}
%
\endbibitem

%b12 ###
\bibitem{cfDon}
%
\begin{barticle}[mr]
\bauthor{\bsnm{Doney},~\bfnm{R.~A.}\binits{R.~A.}}
(\byear{1997}).
\btitle{One-sided local large deviation and renewal theorems in the
case of
infinite mean}.
\bjournal{Probab. Theory Related Fields}
\bvolume{107}
\bpages{451--465}.
\bid{doi={10.1007/s004400050093}, mr={1440141}}
\end{barticle}
%
\endbibitem

%b13 ###
\bibitem{cfFFM}
%
\begin{barticle}[mr]
\bauthor{\bsnm{Fitzsimmons},~\bfnm{P.~J.}\binits{P.~J.}},
\bauthor{\bsnm{Fristedt},~\bfnm{Bert}\binits{B.}} \AND
\bauthor{\bsnm{Maisonneuve},~\bfnm{B.}\binits{B.}}
(\byear{1985}).
\btitle{Intersections and limits of regenerative sets}.
\bjournal{Z. Wahrsch. Verw. Gebiete}
\bvolume{70}
\bpages{157--173}.
\bid{doi={10.1007/BF02451426}, mr={799144}}
\end{barticle}
%
\endbibitem

%b14 ###
\bibitem{cfGHLO}
%
\begin{barticle}[vtex]
\bauthor{\bsnm{Garel},~\bfnm{T.}\binits{T.}},
\bauthor{\bsnm{Huse},~\bfnm{D.~A.}\binits{D.~A.}},
\bauthor{\bsnm{Leibler},~\bfnm{S.}\binits{S.}} \AND
\bauthor{\bsnm{Orland},~\bfnm{H.}\binits{H.}}
(\byear{1989}).
\btitle{Localization transition of random chains at interfaces}.
\bjournal{Europhys. Lett.}
\bvolume{8}
\bpages{9--13}.
\end{barticle}
%
\endbibitem

%b16 ###
\bibitem{cfGB}
%
\begin{bmisc}[vtex]
\bauthor{\bsnm{Giacomin},~\bfnm{G.}\binits{G.}}
(\byear{2004}).
\btitle{Localization phenomena in random polymer models}.
\bnote{Unpublished lecture notes.
Universit\'{e} Paris Diderot (Paris 7)
(available on the web-page of the
author)}.
\end{bmisc}
%
\endbibitem

%b17 ###
\bibitem{cfBook}
%
\begin{bbook}[mr]
\bauthor{\bsnm{Giacomin},~\bfnm{Giambattista}\binits{G.}}
(\byear{2007}).
\btitle{Random Polymer Models}.
\bpublisher{Imperial College Press}, \baddress{London}.
\bid{mr={2380992}}
\end{bbook}
%
\endbibitem

%b18 ###
\bibitem{cfGT}
%
\begin{barticle}[mr]
\bauthor{\bsnm{Giacomin},~\bfnm{Giambattista}\binits{G.}} \AND
\bauthor{\bsnm{Toninelli},~\bfnm{Fabio~Lucio}\binits{F.~L.}}
(\byear{2005}).
\btitle{Estimates on path delocalization for copolymers at selective
interfaces}.
\bjournal{Probab. Theory Related Fields}
\bvolume{133}
\bpages{464--482}.
\bid{doi={10.1007/s00440-005-0439-2}, mr={2197110}}
\end{barticle}
%
\endbibitem

%b20 ###
\bibitem{cfKin}
%
\begin{barticle}[vtex]
\bauthor{\bsnm{Kingman},~\bfnm{J.~F.~C.}\binits{J.~F.~C.}}
(\byear{1973}).
\btitle{Subadditive ergodic theory}.
\bjournal{Ann. Probab.}
\bvolume{1}
\bpages{883--909}.
%Frank Spitzer and J. M. Hammersley, and a reply by the author}.
\bid{mr={0356192}}%
\end{barticle}%
%
\endbibitem%

%b19 ###
\bibitem{cfKin2}
%
\begin{bbook}[vtex]
\bauthor{\bsnm{Kingman},~\bfnm{J.~F.~C.}\binits{J.~F.~C.}}
(\byear{1993}).
\btitle{Poisson Processes}.
\bseries{Oxford Studies in Probability}
\bvolume{3}.
\bpublisher{Oxford Univ. Press}, \baddress{New York}.
\bid{mr={1207584}}
\end{bbook}
%
\endbibitem

%b21 ###
\bibitem{cfKMT}
%
\begin{barticle}[mr]
\bauthor{\bsnm{Koml{\'o}s},~\bfnm{J.}\binits{J.}},
\bauthor{\bsnm{Major},~\bfnm{P.}\binits{P.}} \AND
\bauthor{\bsnm{Tusn{\'a}dy},~\bfnm{G.}\binits{G.}}
(\byear{1976}).
\btitle{An approximation of partial sums of independent {RV}'s, and the sample
{DF}. {II}}.
\bjournal{Z. Wahrsch. Verw. Gebiete}
\bvolume{34}
\bpages{33--58}.
\bid{mr={0402883}}%
\end{barticle}%
%
\endbibitem%

%b22 ###
\bibitem{cfNicolas}
%
\begin{barticle}[vtex]
\bauthor{\bsnm{Petrelis},~\bfnm{Nicolas}\binits{N.}}
(\byear{2009}).
\btitle{Copolymer at selective interfaces and pinning potentials: Weak coupling
limits}.
\bjournal{Ann. Inst. H. Poincar\'e Probab. Statist.}
\bvolume{45}
\bpages{175--200}.
\bid{doi={10.1214/07-AIHP160}, mr={2500234}}
\end{barticle}
%
\endbibitem

%b23 ###
\bibitem{cfRevYor}
%
\begin{bbook}[mr]
\bauthor{\bsnm{Revuz},~\bfnm{Daniel}\binits{D.}} \AND
\bauthor{\bsnm{Yor},~\bfnm{Marc}\binits{M.}}
(\byear{1999}).
\btitle{Continuous Martingales and {B}rownian Motion},
\bedition{3rd} ed.
\bseries{Grundlehren der Mathematischen Wissenschaften [Fundamental Principles
of Mathematical Sciences]}
\bvolume{293}.
\bpublisher{Springer}, \baddress{Berlin}.
\bid{mr={1725357}}
\end{bbook}
%
\endbibitem

%b24 ###
\bibitem{cfSinai}
%
\begin{barticle}[mr]
\bauthor{\bsnm{Sina{\u\i}},~\bfnm{Ya.~G.}\binits{Y.~G.}}
(\byear{1993}).
\btitle{A random walk with a random potential}.
\bjournal{Theory Probab. Appl.}
\bvolume{38}
\bpages{382--385}.
\bid{doi={10.1137/1138036}, mr={1317991}}
\end{barticle}
%
\endbibitem

%b25 ###
\bibitem{cfTAAP}
%
\begin{barticle}[mr]
\bauthor{\bsnm{Toninelli},~\bfnm{Fabio~Lucio}\binits{F.~L.}}
(\byear{2008}).
\btitle{Disordered pinning models and copolymers: Beyond annealed bounds}.
\bjournal{Ann. Appl. Probab.}
\bvolume{18}
\bpages{1569--1587}.
\bid{doi={10.1214/07-AAP496}, mr={2434181}}
\end{barticle}
%
\endbibitem

%b26 ###
\bibitem{cfTcg}
%
\begin{barticle}[mr]
\bauthor{\bsnm{Toninelli},~\bfnm{Fabio~Lucio}\binits{F.~L.}}
(\byear{2009}).
\btitle{Coarse graining, fractional moments and the critical slope of random
copolymers}.
\bjournal{Electron. J. Probab.}
\bvolume{14}
\bpages{no. 20, 531--547}.
\bid{mr={2480552}}
\end{barticle}
%
\endbibitem

\end{thebibliography}
\end{document}